# ROCKET: <u>Ro</u>bust <u>C</u>onfidence Intervals via <u>Ke</u>ndall's <u>T</u>au for Transelliptical Graphical Models

Rina Foygel Barber and Mladen Kolar

September 1, 2017

**Abstract**

Understanding complex relationships between random variables is of fundamental importance in high-dimensional statistics, with numerous applications in biological and social sciences. Undirected graphical models are often used to represent dependencies between random variables, where an edge between two random variables is drawn if they are conditionally dependent given all the other measured variables. A large body of literature exists on methods that estimate the structure of an undirected graphical model, however, little is known about the distributional properties of the estimators beyond the Gaussian setting. In this paper, we focus on inference for edge parameters in a high-dimensional transelliptical model, which generalizes Gaussian and nonparanormal graphical models. We propose ROCKET, a novel procedure for estimating parameters in the latent inverse covariance matrix. We establish asymptotic normality of ROCKET in an ultra high-dimensional setting under mild assumptions, without relying on oracle model selection results. ROCKET requires the same number of samples that are known to be necessary for obtaining a $\sqrt{n}$ consistent estimator of an element in the precision matrix under a Gaussian model. Hence, it is an optimal estimator under a much larger family of distributions. The result hinges on a tight control of the sparse spectral norm of the non-parametric Kendall's tau estimator of the correlation matrix, which is of independent interest. Empirically, ROCKET outperforms the nonparanormal and Gaussian models in terms of achieving accurate inference on simulated data. We also compare the three methods on real data (daily stock returns), and find that the ROCKET estimator is the only method whose behavior across subsamples agrees with the distribution predicted by the theory.

## 1 Introduction

Probabilistic graphical models [Lauritzen, 1996] have been widely used to explore complex system and aid scientific discovery in areas ranging from biology and neuroscience to financial modeling and social media analysis. An undirected graphical model consists of a graph $G = (V, E)$, where $V = \{1, \ldots, p\}$ is the set of vertices and $E$ is the set of edges, and a $p$-dimensional random vector $X = (X_1, \ldots, X_p)^\top$ that is Markov with respect to $G$. In particular, we have that $X_a$ and $X_b$ are conditionally independent given the remaining variables $\{X_c \mid c \in \{1, \ldots, p\} \backslash \{a, b\}\}$ if and only if $\{a, b\} \notin E$. One of the central questions in high-dimensional statistics is estimation of the undirected graph $G$ given $n$ independent realizations of $X$, as well as quantifying uncertainty of the estimator.

In this paper we focus on (asymptotic) inference for elements in the latent inverse covariance matrix under the semiparametric elliptical copula model [Embrechts et al., 2003, Klüppelberg et al., 2008], also known as the transelliptical model [Liu et al., 2012b]. Let $X_1, \ldots, X_n$ be $n$ independent copies of the random vector $X$ that follows a transelliptical distribution,

$$X \sim \mathsf{TE}(\Sigma, \xi; f_1, \ldots, f_p), \tag{1.1}$$

where $\Sigma \in \mathbb{R}^p$ is a correlation matrix (that is, $\Sigma_{jj} = 1$ for $j = 1, \ldots, p$), $\xi \in \mathbb{R}$ is a nonnegative random variable with $\mathbb{P}\{\xi = 0\} = 0$, and $f_1, \ldots, f_p$ are univariate, strictly increasing functions. Recall that $X$ follows a transelliptical distribution if the marginal transformation $(f_1(X_1), \ldots, f_p(X_p))$ of $X$ follows a (centered) elliptically contoured distribution with covariance matrix $\Sigma$ [Fang et al., 1990]. Let $\Omega = \Sigma^{-1}$ be the inverse covariance matrix, also known as the precision matrix. Under a Gaussian model, nonzero elements in $\Omega$ correspond to pairs of variables that are conditionally dependent, that is, form an edge in the graph $G$; under an elliptical model, nonzero elements in $\Omega$ correspond to variables that are conditionally correlated (but in general it is possible to have $\Omega_{ab} = 0$ where $f_a(X_a)$ and $f_b(X_b)$ are conditionally uncorrelated, but not conditionally independent). Under the model in (1.1), we construct an estimator for a fixed element of the precision matrix, $\Omega_{ab}$, that is asymptotically normal. Furthermore, we construct a confidence interval for the unknown parameter $\Omega_{ab}$ that is valid and robust to model selection mistakes. Finally, we construct a uniformly valid hypothesis test for the presence of an edge in the graphical model.

Our main theoretical result establishes that given initial estimates of the regression coefficients for $(f_a(X_a), f_b(X_b))$ on $(f_j(X_j))_{j \neq a,b}$, one can obtain a $\sqrt{n}$-consistent and asymptotically normal estimator for $\Omega_{ab}$. These initial estimators need to converge at a sufficiently fast rate (see Section 3). In particular, we note that we do not require strict sparsity in these regressions, and allow for an error rate that is achievable by known methods such as a nonconvex Lasso [Loh and Wainwright, 2013] (see Section 3.1). To achieve $\sqrt{n}$-consistent rate, our estimator requires the same scaling for the sample size $n$ as in the Gaussian case; this sample size scaling is minimax optimal [Ren et al., 2013].

Given accurate initial estimates, in order to construct the asymptotically normal estimator, we prove a key result: that the vector $\text{sign}(X_i - X_{i'})$ is subgaussian at the scale $C(\Sigma)$ (the condition number of $\Sigma$), with dependence on the dimension $p$ coming only through $C(\Sigma)$ (this problem was initially posed by Han and Liu [2013], where subgaussianity was proved for some special cases). This result allows us to construct an asymptotically normal estimator by combining the initial regression coefficient estimates with the Kendall's tau rank correlation matrix. In particular, the subgaussianity result allows us to establish a new concentration result on the operator norm of the Kendall's tau correlation matrix that holds with exponentially high probability. This result allows us to uniformly control deviations of quadratic forms involving the Kendall's tau correlation matrix over approximately sparse vectors. These results are of independent interest and could be used to extend recent results of Mitra and Zhang [2014], Wegkamp and Zhao [2013] and Han and Liu [2013] to the elliptical copula setting. Furthermore, subgaussianity of $\text{sign}(X_i - X_{i'})$, which in turn leads to a bound on the error of the Kendall's tau estimate of $\Sigma$ in the sparse spectral norm, allows us to study properties of penalized rank regression in high-dimensions.

We base our confidence intervals and hypothesis tests on the asymptotically normal estimator of the element $\Omega_{ab}$ (see Section 2). We point out that our results hold under milder conditions than those required in Ren et al. [2013], which treats the special case of Gaussian graphical models. Most notably, we give a $\sqrt{n}$-consistent estimator for elements in the precision matrix without requiring strong parametric assumptions.

## 1.1 Relationship To Literature

Our work contributes to several areas. First, we contribute to the growing literature on graphical model selection in high dimensions. There is extensive literature on the Gaussian graphical model, where it is assumed that $X \sim N(0, \Sigma)$, in which case the edge set $E$ of the graph $G$ is encoded by the non-zero elements of the precision matrix $\Omega$ [Meinshausen and Bühlmann, 2006, Yuan and Lin, 2007, Rothman et al., 2008, Friedman et al., 2008, d'Aspremont et al., 2008, Fan et al., 2009, Lam and Fan, 2009, Yuan, 2010, Cai et al., 2011, Liu and Wang, 2012, Zhao and Liu, 2014]. Learning structure of the Ising model based on the penalized pseudo-likelihood was studied in Höfling and Tibshirani [2009], Ravikumar et al. [2010] and Xue et al. [2012]. More recently, Yang et al. [2013] studied estimation of graphical models under the assumption that each of the nodes' conditional distribution belongs to an exponential family distribution. See also Guo et al. [2011a], Guo et al. [2011b], Lee and Hastie [2012], Cheng et al. [2013], Yang et al. [2012] and Yang et al. [2014] who studied mixed graphical models, where the nodes' conditional distributions are not necessarily all from the same family (for instance, there may be continuous-valued nodes as well as discrete-valued nodes). The parametric Gaussian assumption was relaxed in Liu et al. [2009], where graph estimation was studied under a Gaussian copula model. More recently, Liu et al. [2012a], Xue and Zou [2012] and Liu et al. [2012b] show that the graph can be recovered in the Gaussian and elliptical semiparametric model class under the same conditions on the sample size $n$, number of nodes $p$ and the maximum node degree in the graph $k$ as if the estimation was done under the Gaussian assumption. In our paper, we construct a novel $\sqrt{n}$-consistent estimator of an element in the precision matrix without requiring oracle model selection properties.

Second, we contribute to the literature on high-dimensional inference. Recently, there has been much interest on performing valid statistical inference in the high-dimensional setting. Zhang and Zhang [2013], Belloni et al. [2013a], Belloni et al. [2013d], van de Geer et al. [2014], Javanmard and Montanari [2014], Javanmard and Montanari [2013], and Farrell [2013] developed methods for construction of confidence intervals for low dimensional parameters in high-dimensional linear and generalized linear models, as well as hypothesis tests. These methods construct honest, uniformly valid confidence intervals and hypothesis test based on the $\ell_1$-penalized estimator in the first stage. Similar results were obtained in the context of the $\ell_1$-penalized least absolute deviation and quantile regression [Belloni et al., 2013c,b]. Lockhart et al. [2014] study significance of the input variables that enter the model along the lasso path. Lee et al. [2013] and Taylor et al. [2014] perform post-selection inference conditional on the selected model. Liu [2013], Ren et al. [2013] and Chen et al. [2013] construct $\sqrt{n}$-consistent estimators for elements of the precision matrix $\Omega$ under a Gaussian assumption. We extend these results to perform valid inference under semiparametric ellitical copula models. In a recent independent work, Gu et al. [2015] propose a procedure for inference under a nonparanormal model. We will provide a detailed comparison in Section 3 and Section 5.

## 1.2 Notation

Let $[n]$ denote the set $\{1, \ldots, n\}$ and let $\mathbb{I}\{\cdot\}$ denote the indicator function. For a vector $a \in \mathbb{R}^d$, we let $\text{supp}(a) = \{j \ : \ a_j \neq 0\}$ be the support set, and let $||a||_q$, for $q \in [1, \infty)$, be the $\ell_q$-norm defined as $||a||_q = (\sum_{i \in [n]} |a_i|^q)^{1/q}$ with the usual extensions for $q \in \{0, \infty\}$, that is, $||a||_0 = |\text{supp}(a)|$ and $||a||_\infty = \max_{i \in [n]} |a_i|$.

For a matrix $A \in \mathbb{R}^{n_1 \times n_2}$, for sets $S \subset [n_1]$ and $T \subset [n_2]$, we write $A_{ST}$ to denote the $|S| \times |T|$ submatrix of $A$ obtained by extracting the appropriate rows and columns. The sets $S$ and/or $T$ can be replaced by single indices, for example, for $S \subset [n_1]$ and $j \in [n_2]$, $A_{Sj}$ is a $|S|$-length vector. If $A \in \mathbb{R}^{n \times n}$ is a square matrix, for any $T \subset [n]$ we may write $A_T$ to denote the square submatrix $A_{TT}$.

For a matrix $A \in \mathbb{R}^{n_1 \times n_2}$, we use the notation $\text{vec}(A)$ to denote the vector in $\mathbb{R}^{n_1 n_2}$ formed by stacking the columns of $A$. We denote the Frobenius norm of $A$ by $||A||_\mathsf{F}^2 = \sum_{i \in [n_1], j \in [n_2]} A_{ij}^2$, and the operator norm (spectral norm) by $||A||_{\mathsf{op}}$, that is, the largest singular value of $A$. The norms $||A||_1$ and $||A||_\infty$ are applied entrywise, with $||A||_1 = \sum_{ij} |A_{ij}|$ and $||A||_\infty = \max_{ij} |A_{ij}|$. We write $\mathsf{C}(A)$ to denote the condition number of $A$, that is, the ratio between the largest and smallest singular values. For two matrices $A \in \mathbb{R}^{n \times m}$ and $B \in \mathbb{R}^{r \times s}$, $A \otimes B \in \mathbb{R}^{nr \times ms}$ denotes the Kronecker product, with $(A \otimes B)_{ik,jl} = A_{ij} B_{kl}$. For two matrices of the same size, $A, B \in \mathbb{R}^{n \times m}$, $A \circ B \in \mathbb{R}^{n \times m}$ denotes the Hadamard product (that is, the entrywise product), with $(A \circ B)_{ij} = A_{ij} B_{ij}$. Kronecker products and Hadamard products are defined also for vectors, by treating a vector as a matrix with one column.

Throughout, $\Phi(\cdot)$ denotes the cumulative distribution function of the standard normal distribution, that is, $\Phi(t) = \mathbb{P}\{N(0,1) \leq t\}$.

## 2 Preliminaries and method

Before introducing our method, we begin with some preliminary definitions and properties of the transelliptical distribution, and related models.

**Gaussian and nonparanormal graphical models**  Suppose that $X = (X_1, \ldots, X_p)$ follows a multivariate normal distribution, $X \sim N(\mu, \Sigma)$. A Gaussian graphical model represents the structure of the covariance matrix $\Sigma$ with a graph, where an edge between nodes $a$ and $b$ indicates that $\Omega_{ab} \neq 0$, where $\Omega = \Sigma^{-1}$ is the precision (inverse covariance) matrix. This model can be generalized by allowing for arbitrary marginal transformations on the variables $X_1, \ldots, X_p$. Liu et al. [2009] study the resulting distribution, the nonparanormal model (also known as a Gaussian copula), where we write $X \sim \mathsf{NPN}(\Sigma; f_1, \ldots, f_p)$, if the marginally transformed vector $(f_1(X_1), \ldots, f_p(X_p))$ follows a (centered) multivariate normal distribution,

$$(f_1(X_1), \ldots, f_p(X_p)) \sim N(0, \Sigma) \ .$$

The sparse structure of the underlying graphical model, representing the sparsity pattern in $\Omega = \Sigma^{-1}$, can then be recovered using similar methods as in the Gaussian case. Note that the Gaussian model is a special case of the nonparanormal model (by setting $f_1, \ldots, f_p$ each to be the identity function, or to be linear functions if we would like a nonzero mean).

**Elliptical and transelliptical graphical models**  The elliptical model is a generalization of the Gaussian graphical model that allows for heavier-tailed dependence between variables. The random vector $X = (X_1, \ldots, X_p)$ follows an elliptical distribution with the mean vector $\mu \in \mathbb{R}^p$, covariance matrix $\Sigma \in \mathbb{R}^{p \times p}$, and a random variable (the "radius") $\xi \geqslant 0$, denoted by $X \sim \mathsf{E}(\mu, \Sigma, \xi)$, if we can write $X = \mu + \xi \cdot A \cdot U$, where $AA^\top = \Sigma$ is a Cholesky decomposition of $\Sigma$, and where $U \in \mathbb{R}^p$ is a unit vector drawn uniformly at random (independently from the radius $\xi$). Note that the level sets of this distribution are given by ellipses, centered at $\mu$ and with shape determined by $\Sigma$. The Gaussian model is a special case of the elliptical model (by taking $\xi \sim \chi_p$).

The transelliptical model (also known as an elliptical copula) combines the elliptical distribution with marginal transformations, much as the nonparanormal distribution applies marginal transformations to a multivariate Gaussian. For a random vector $X \in \mathbb{R}^p$ we write

$$X \sim \mathsf{TE}(\Sigma, \xi; f_1, \ldots, f_p)$$

to denote that the marginally transformed vector $(f_1(X_1), \ldots, f_p(X_p))$ follows a centered elliptical distribution, specifically,

$$(f_1(X_1), \ldots, f_p(X_p)) \sim \mathsf{E}(0, \Sigma, \xi) \ .$$

Here the marginal transformation functions $f_1, \ldots, f_p$ are assumed to be strictly increasing. Note that the Gaussian, nonparanormal, and elliptical models are each special cases of this model.

**Pearson's rho and Kendall's tau**  From this point on, we assume for each distribution that $\mu = 0$ and that $\Sigma$ is a correlation matrix (that is, diagonal elements are equal to one, $\Sigma_{aa} = 1$). In the case of the Gaussian distribution $X \sim N(0, \Sigma)$, the entries of $\Sigma$ are the (population-level) Pearson's correlation coefficients for each pair of variables, which in this case we can also write as $\Sigma_{ab} = \mathbb{E}[X_a X_b]$. In this setting, we can estimate $\Sigma$ with the sample covariance.

In the nonparanormal setting, $X \sim \mathsf{NPN}(\Sigma; f_1, \ldots, f_p)$, it is no longer the case that $\Sigma_{ab}$ is equal to the (population-level) correlation $\mathrm{Corr}(X_a, X_b)$, due to the marginal transformations. However, we can estimate $f_1, \ldots, f_p$ by performing marginal empirical transformations of each $X_a$ to the standard normal distribution. After taking these empirical transformations, $\Sigma$ can again be estimated via the empirical covariances. Similarly, for the elliptical model $X \sim \mathsf{E}(0, \Sigma, \xi)$, after rescaling so that $\mathbb{E}[\xi^2] = p$ we also have $\Sigma_{ab} = \mathbb{E}[X_a X_b]$. We can therefore again estimate $\Sigma$ via the empirical covariance.

For the transelliptical distribution, in contrast, this is no longer possible. Taking scaling $\mathbb{E}[\xi^2] = p$ for simplicity, we generalize the calculations above to have $\Sigma_{ab} = \mathbb{E}[f_a(X_a) f_b(X_b)]$. Therefore, if we can estimate the marginal transformations $f_1, \ldots, f_p$, then we can estimate $\Sigma$ using the empirical covariance of the transformed data. However, unlike the nonparanormal model, estimating $f_1, \ldots, f_p$ is not straightforward. The reason is that, for the elliptical distribution $\mathsf{E}(0, \Sigma, \xi)$, the marginal distributions are not known unless the distribution of the radius $\xi$ is known. Therefore, marginally for each $X_a$, we cannot estimate $f_a$ because we do not know what should be the marginal distribution after transformation, that is, what should be the marginal distribution of $f_a(X_a)$. (In contrast, in the nonparanormal model, $f_a(X_a)$ is marginally normal.)

As an alternative, Liu et al. [2012b] use the Kendall rank correlation coefficient (Kendall's tau). At the population level, it is given by

$$\tau_{ab} := \tau(X_a, X_b) = \mathbb{E}\left[\text{sign}(X_a - X'_a) \cdot \text{sign}(X_b - X'_b)\right],$$

where $X'$ is an i.i.d. copy of $X$. Unlike Pearson's rho, the Kendall's tau coefficient is invariant to marginal transformations: since $f_a, f_b$ are strictly increasing functions, we see that

$$\text{sign}(f_a(X_a) - f_a(X'_a)) \cdot \text{sign}(f_b(X_b) - f_b(X'_b)) = \text{sign}(X_a - X'_a) \cdot \text{sign}(X_b - X'_b).$$

At the sample level, Kendall's tau can be estimated by taking a U-statistic comparing each pair of distinct observations:

$$\widehat{\tau}_{ab} = \frac{1}{\binom{n}{2}} \sum_{1 \leqslant i < i' \leqslant n} \text{sign}(X_{ia} - X_{i'a}) \cdot \text{sign}(X_{ib} - X_{i'b}). \tag{2.1}$$

When $X$ follows an elliptical distribution, Theorem 2 of Lindskog et al. [2003] gives us the following relationship between Kendall's tau and the Pearson's rho coefficients given by the covariance matrix $\Sigma$:

$$\Sigma_{ab} = \sin\left(\frac{\pi}{2} \tau_{ab}\right) \text{ for each } a, b \in [p]. \tag{2.2}$$

Since Kendall's tau is invariant to marginal transformations, this identity holds for the transelliptical family as well. For this reason, Liu et al. [2012b] estimate the covariance matrix $\Sigma$ by

$$\widehat{\Sigma}_{ab} = \sin\left(\frac{\pi}{2} \widehat{\tau}_{ab}\right). \tag{2.3}$$

Note, however, that $\widehat{\Sigma}$ is not necessarily positive semidefinite.

While Spearman's rho, like Kendall's tau, is also invariant to marginal transformations, Liu et al. [2012b] comment that there is no equivalence between $\Sigma$ and the population-level Spearman's rho values (analogous to (2.2) for Kendall's tau) which holds uniformly across the entire elliptical (or transelliptical) family. Therefore, this type of estimator as in (2.3) could only be carried out with Kendall's tau.

For the remainder of this paper, $\widehat{\Sigma}$ denotes the estimate given here in (2.3). The matrix of the Kendall's tau coefficients is denoted as $T$, with entries $T_{ab} := \tau_{ab}$, and $\widehat{T}$ denotes its empirical estimate (with entries as in (2.1)).

**Comparing models: tail dependence** It is clear that, compared to a Gaussian graphical model, the nonparanormal model allows for data that may be extremely heavy-tailed (in the marginal distributions). A more subtle consideration is the question of tail dependence between two or more of the variables. In particular, the nonparanormal model does not allow for tail dependence between two variables to be any stronger than in the Gaussian distribution itself. Specifically, consider pairwise $\alpha$-tail dependence between $X_a$ and $X_b$, given by

$$\mathsf{Tail}_\alpha(X_a, X_b) := \text{Corr}\left(\mathbb{1}\left\{X_a \geqslant q^{X_a}_\alpha\right\}, \mathbb{1}\left\{X_b \geqslant q^{X_b}_\alpha\right\}\right),$$

where $q^{X_a}_\alpha$ is the $\alpha$-quantile of the marginal distribution of $X_a$, and same for $X_b$. Taking $\alpha \to 1$, this is a measure of the correlation between the extreme right tail of $X_a$ and the extreme right tail of $X_b$. (Of course, we can also consider the left tail of the distribution of $X_a$ and/or $X_b$.)

Note that marginal transformations of each variable do not affect this measure, since the quantiles $q_\alpha^{X_a}, q_\alpha^{X_b}$ take these transformations into account. In particular, the nonparanormal distribution has the same tail correlations $\mathsf{Tail}_\alpha(X_a, X_b)$ as the multivariate Gaussian distribution (with the same $\Sigma$). In contrast, an elliptical or transelliptical model can exhibit much higher tail correlations. Since real data often exhibits heavy tail dependence between variables, the flexible transelliptical model may be a better fit in many applications.

We demonstrate this behavior with a simple example in Figure 1. Take

$$X = (X_1, X_2) \sim \mathsf{E}(0, \Sigma, \xi) \text{ with } \Sigma = \begin{pmatrix} 1 & 1/\sqrt{2} \\ 1/\sqrt{2} & 1 \end{pmatrix}, \quad (2.4)$$

where $\xi \sim \chi_2 \cdot \sqrt{d}/\chi_d$ for $d \in \{0.1, 1, 5, 10, \infty\}$, corresponding to a multivariate t-distribution with $d$ degrees of freedom (note that $d = \infty$ is equivalent to taking $X \sim N(0, \Sigma)$). Note that at $\alpha = 0.5$, the relevant quantiles are $q_\alpha^{X_1} = q_\alpha^{X_2} = 0$, and so the tail correlation $\mathsf{Tail}_\alpha(X_1, X_2)$ is equal to the Kendall's tau coefficient $\tau(X_1, X_2) = \frac{2}{\pi} \arcsin(\Sigma_{12}) = 0.5$ at any value of $d$. Figure 1 shows that, as $\alpha \to 1$, the tail correlation decreases towards zero for the normal distribution ($d = \infty$) but grows for low values of $d$.

Therefore, the shift from a nonparanormal to a transelliptical model is important, since it allows us to model variables with high tail dependence, that is, high dependence between their "extreme events".

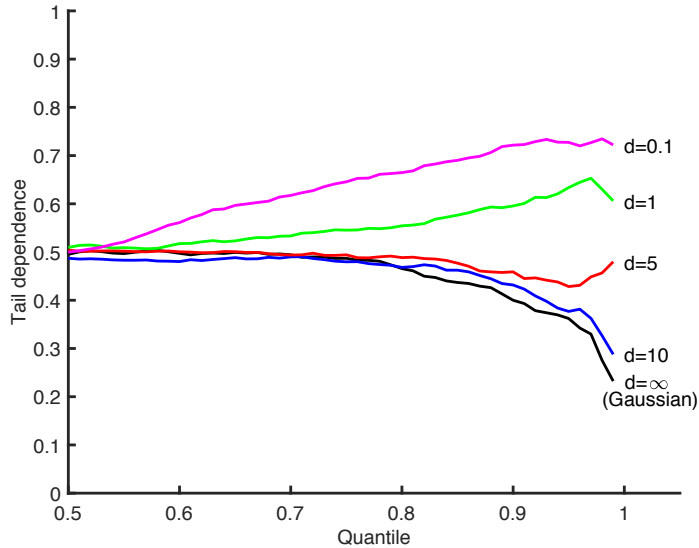

Figure 1: Tail dependence for normal and elliptical distributions on $\mathbb{R}^2$. Data is generated as in (2.4). The figure displays $\mathsf{Tail}_\alpha(X_1, X_2)$, estimated empirically from a sample size $n = 20000$.

## 2.1 ROCKET: an asymptotically normal estimator

Suppose that our data points $X_i$ are drawn i.i.d. from a transelliptical distribution with covariance matrix $\Sigma$. We would like to perform inference on a particular entry of the precision matrix $\Omega = \Sigma^{-1}$, specifically, we are interested in producing a confidence interval for $\Omega_{ab}$ where $a \neq b \in \{1, \ldots, p\}$ is a prespecified node pair.

To move towards constructing a confidence interval, we introduce a few definitions and calculations. First, let $I = \{1, \ldots, p\} \setminus \{a, b\}$, and observe that by block-wise matrix inversion, we can calculate the

$\{a, b\} \times \{a, b\}$ sub-block of $\Omega$ as follows:

$$\Omega_{ab,ab} = \left(\Sigma_{ab,ab} - \Sigma_{ab,I}\Sigma_I^{-1}\Sigma_{I,ab}\right)^{-1}. \tag{2.5}$$

Define $\gamma_a = \Sigma_I^{-1}\Sigma_{Ia}$ and $\gamma_b = \Sigma_I^{-1}\Sigma_{Ib}$. In the nonparanormal graphical model setting, these are the regression coefficients when $f_a(X_a)$ or $f_b(X_b)$ is regressed on $\{f_j(X_j) : j \in I\}$; in the linear model setting, this idea has been used in Sun and Zhang [2012b] and Belloni et al. [2013a]. We then have

$$\Sigma_{ab,I}\Sigma_I^{-1}\Sigma_{I,ab} = (\gamma_a\ \gamma_b)^\top \Sigma_{I,ab} = \Sigma_{I,ab}^\top (\gamma_a\ \gamma_b) = (\gamma_a\ \gamma_b)^\top \Sigma_I (\gamma_a\ \gamma_b).$$

We can therefore rewrite (2.5) as follows (this somewhat redundant formulation will allow for a favorable cancellation of error terms later on):

$$\Theta := (\Omega_{ab,ab})^{-1} = \Sigma_{ab,ab} - (\gamma_a\ \gamma_b)^\top \Sigma_{I,ab} - \Sigma_{I,ab}^\top(\gamma_a\ \gamma_b) + (\gamma_a\ \gamma_b)^\top \Sigma_I(\gamma_a\ \gamma_b). \tag{2.6}$$

We abuse notation and index the entries of $\Theta$ with the indices $a$ and $b$, that is, we denote $\Theta$ as lying in $\mathbb{R}^{\{a,b\}\times\{a,b\}}$ rather than $\mathbb{R}^{2\times 2}$.

Next, we define an oracle estimator of $\Theta$, defined by plugging the *true* values of $\gamma_a$ and $\gamma_b$ and the *empirical* estimate of $\Sigma$ (given in (2.3)) into (2.6) above:

$$\widetilde{\Theta} = \widehat{\Sigma}_{ab,ab} - (\gamma_a\ \gamma_b)^\top \widehat{\Sigma}_{I,ab} - \widehat{\Sigma}_{I,ab}^\top(\gamma_a\ \gamma_b) + (\gamma_a\ \gamma_b)^\top \widehat{\Sigma}_I(\gamma_a\ \gamma_b). \tag{2.7}$$

Later on (in Theorem 4.1), we will show that due to standard results on the theory of U-statistics, this oracle estimator is asymptotically normal. If $\widetilde{\Theta}$ were known, then, we would have achieved our goal for inference in this model, as $\widetilde{\Omega}_{ab} = \left(\widetilde{\Theta}^{-1}\right)_{ab}$ weakly converges to a Normal random variable centered at $\Omega_{ab}$ with variance that scales as $\mathcal{O}(1/n)$ (we calculate this variance later).

Of course, in practice we do not know the true values of $\gamma_a$ and $\gamma_b$, and must instead use some available estimators, denoted by $\breve{\gamma}_a$ and $\breve{\gamma}_b$ (we discuss how to obtain these preliminary estimates later on). Given the estimators of the regression vectors, we then define our estimator of $\Theta$ as follows:

$$\breve{\Theta} = \widehat{\Sigma}_{ab,ab} - (\breve{\gamma}_a\ \breve{\gamma}_b)^\top \widehat{\Sigma}_{I,ab} - \widehat{\Sigma}_{I,ab}^\top(\breve{\gamma}_a\ \breve{\gamma}_b) + (\breve{\gamma}_a\ \breve{\gamma}_b)^\top \widehat{\Sigma}_I(\breve{\gamma}_a\ \breve{\gamma}_b). \tag{2.8}$$

Since we are interested in $\Omega_{ab}$ rather than in the matrix $\Theta$, as a final step we define our estimator

$$\breve{\Omega}_{ab} = \left(\breve{\Theta}^{-1}\right)_{ab}. \tag{2.9}$$

In order to make inference about $\Omega_{ab}$, we approximate the distribution of $\breve{\Omega}_{ab}$, which is a function of $\breve{\Theta}$. We first treat the distribution of the corresponding entry in the oracle estimator $\widetilde{\Theta}$. To do so, let $u, v \in \mathbb{R}^{p_n}$ be the vectors with entries

$$u_a = 1, u_b = 0, u_I = -\gamma_a \text{ and } v_a = 0, v_b = 1, v_I = -\gamma_b,$$

and observe from (2.6) and (2.7) that $\Theta_{ab} = u^\top \sin\left(\frac{\pi}{2}T\right)v$ while $\widetilde{\Theta}_{ab} = u^\top \sin\left(\frac{\pi}{2}\widehat{T}\right)v$. Now, taking a linear approximation to $\sin(\cdot)$, we can write

$$\widetilde{\Theta}_{ab} - \Theta_{ab} \approx \left\langle uv^\top \circ \frac{\pi}{2}\cos\left(\frac{\pi}{2}\widehat{T}\right), \widehat{T} - T\right\rangle$$
$$= \frac{1}{\binom{n}{2}}\sum_{i<i'}\left\langle uv^\top \circ \frac{\pi}{2}\cos\left(\frac{\pi}{2}\widehat{T}\right), \text{sign}(X_i - X_{i'})\text{sign}(X_i - X_{i'})^\top - T\right\rangle.$$

To study the variability of this error, we consider the kernel

$$g(X, X') = \text{sign}(X - X')^\top \left(uv^\top \circ \cos\left(\frac{\pi}{2}T\right)\right)\text{sign}(X - X').$$

We will see later on that understanding the behavior of this kernel will allow us to characterize the distribution of the oracle estimator $\widetilde{\Theta}_{ab}$, and from there, our empirical estimator $\breve{\Theta}_{ab}$ and ultimately $\breve{\Omega}_{ab}$. Of course, $g(X, X')$ itself depends on unknown quantities, namely $u$, $v$, and $T$, so we replace these with their estimates in our empirical version of the kernel: define the (random) kernel

$$\breve{g}(X, X') = \text{sign}(X - X')^\top \left(\breve{u}\breve{v}^\top \circ \cos\left(\frac{\pi}{2}\widehat{T}\right)\right) \text{sign}(X - X'),$$

where

$$\breve{u}_a = 1, \breve{u}_b = 0, \breve{u}_I = -\breve{\gamma}_a \text{ and } \breve{v}_a = 0, \breve{v}_b = 1, \breve{v}_I = -\breve{\gamma}_b.$$

(Note that we have defined $\breve{u}$ and $\breve{v}$ so that $\breve{\Theta}_{ab} = \breve{u}^\top \widehat{\Sigma} \breve{v}$.) Then define

$$\breve{S}_{ab} = \frac{\pi}{\det(\breve{\Theta})} \cdot \sqrt{\frac{1}{n}\sum_i \left(\frac{1}{n-1}\sum_{i' \neq i} \breve{g}(X_i, X_{i'}) - \text{mean}(\breve{g})\right)^2}$$

where $\text{mean}(\breve{g}) = \binom{n}{2}^{-1} \sum_{i < i'} \breve{g}(X_i, X_{i'})$. We will see later on that $\breve{S}_{ab}^2/n$ estimates the variance of $\breve{\Theta}_{ab}$ and that the expression above arises naturally from the theory of U-statistics.

Our main result, Theorem 3.5 below, will prove that $\sqrt{n} \cdot \frac{\breve{\Omega}_{ab} - \Omega_{ab}}{\breve{S}_{ab}}$ follows a distribution that is approximately standard normal. Therefore, an approximate $(1-\alpha)$-confidence interval for $\Omega_{ab}$ is given by

$$\breve{\Omega}_{ab} \pm z_{\alpha/2} \cdot \frac{\breve{S}_{ab}}{\sqrt{n}}, \tag{2.10}$$

where $z_{\alpha/2}$ is the appropriate quantile of the normal distribution, that is, $\mathbb{P}\left\{N(0,1) > z_{\alpha/2}\right\} = \alpha/2$.

**Notation for fixed vs random quantities** From this point on, as much as possible throughout the main body of the paper, quantities that depend on the data and depend on the initial estimates $\breve{\gamma}_a, \breve{\gamma}_b$ are denoted with a "check" accent, for example, $\breve{\Theta}$. Quantities that depend on the data, but do not depend on $\breve{\gamma}_a, \breve{\gamma}_b$, are denoted with a "hat" accent, for example, $\widehat{\Sigma}$. Any quantities with neither a "hat" nor a "check" are population quantities, that is, they are not random. Two important exceptions are the data itself, $X_1, \ldots, X_n$, and the oracle estimator, $\widetilde{\Theta}$, which is of course data-dependent (but does not depend on $\breve{\gamma}_a, \breve{\gamma}_b$).

## 3 Main results

In this section, we give a theoretical result showing that the confidence interval constructed in (2.10) has asymptotically the correct coverage probability, as long as we have reasonably accurate estimators of $\gamma_a = \Sigma_I^{-1}\Sigma_{Ia}$ and $\gamma_b = \Sigma_I^{-1}\Sigma_{Ib}$. Our asymptotic result considers a problem whose dimension $p_n \geq 2$ grows with the sample size $n$. We also allow for the sparsity level in the true inverse covariance matrix $\Omega \in \mathbb{R}^{p_n \times p_n}$ to grow.[1] We use $k_n$ to denote an approximate bound on the sparsity in each column of $\Omega$ (details given below).

We begin by stating several assumptions on the distribution of the data and on the initial estimators $\breve{\gamma}_a$ and $\breve{\gamma}_b$. All of the constants appearing in these assumptions should be interpreted as values that do not depend on the dimensions $(n, p_n, k_n)$ of the problem.

**Assumption 3.1.** *The data points $X_1, \ldots, X_n \in \mathbb{R}^{p_n}$ are i.i.d. draws from a transelliptical distribution, $X_i \overset{iid}{\sim} \mathsf{TE}(\Sigma, \xi; f_1, \ldots, f_{p_n})$, where $f_1, \ldots, f_{p_n}$ are any strictly monotone functions, $\xi \geq 0$ is any*

---

[1] While $\Sigma$, $\Omega$, etc, all depend on the sample size $n$ since the dimension of the problem grows, we abuse notation and do not write $\Sigma_n$, $\Omega_n$, etc; the dependence on $n$ is implicit.

*random variable with* $\mathbb{P}\{\xi = 0\} = 0$, *and the covariance matrix* $\Sigma \in \mathbb{R}^{p_n \times p_n}$ *is positive definite, with* $\mathrm{diag}(\Sigma) = \mathbf{1}$ *and bounded condition number,*

$$\mathsf{C}(\Sigma) = \frac{\lambda_{\max}(\Sigma)}{\lambda_{\min}(\Sigma)} \leqslant C_{\mathsf{cov}} ,$$

*for some constant* $C_{\mathsf{cov}}$.

**Assumption 3.2.** *The $a$-th and $b$-th columns of the true inverse covariance $\Omega$, denoted by $\Omega_a$ and $\Omega_b$, are approximately $k_n$-sparse, with*

$$||\Omega_a||_1 \vee ||\Omega_b||_1 \leqslant C_{\mathsf{sparse}} \sqrt{k_n} ,$$

*for some constant* $C_{\mathsf{sparse}}$.

**Assumption 3.3.** *For some constant $C_{\mathsf{est}}$ and for some $\delta_n > 0$, with probability at least $1 - \delta_n$, for each $c = a, b$, the preliminary estimate $\breve{\gamma}_c$ of the vector $\gamma_c$ satisfies*

$$||\breve{\gamma}_c - \gamma_c||_2 \leqslant C_{\mathsf{est}} \sqrt{\frac{k_n \log(p_n)}{n}}, \ ||\breve{\gamma}_c - \gamma_c||_1 \leqslant C_{\mathsf{est}} \sqrt{\frac{k_n^2 \log(p_n)}{n}} . \qquad (3.1)$$

**Assumption 3.4.** *Define the kernel $h(X, X') = \mathrm{sign}(X - X') \otimes \mathrm{sign}(X - X') \in \mathbb{R}^{p_n^2}$ and let $h_1(X) = \mathbb{E}[h(X, X') \mid X]$. Define the total variance $\Sigma_h = \mathsf{Var}(h(X, X'))$ and the conditional $\Sigma_{h_1} = \mathsf{Var}(h_1(X))$, where $X, X' \stackrel{iid}{\sim} \mathsf{TE}(\Sigma, \xi; f_1, \ldots, f_{p_n})$. Then for some constant $C_{\mathsf{kernel}} > 0$,[2]*

$$C_{\mathsf{kernel}} \cdot \Sigma_h \preceq \Sigma_{h_1} \preceq \Sigma_h .$$

Assumption 3.1 assumes that the smallest and largest eigenvalues of the correlation matrix $\Sigma$ are bounded away from zero and infinity, respectively. This assumption is commonly assumed in the literature on learning structure of probabilistic graphical models [Ravikumar et al., 2011, Liu et al., 2009, 2012a]. Assumption 3.2 does not require that the precision matrix $\Omega$ be exactly sparse, which is commonly assumed in the literature on exact graph recovery [see, for example, Ravikumar et al., 2011], but only requires that rows $\Omega_a$ and $\Omega_b$ have an $\ell_1$ norm that does not grow too fast. Note that if $\Omega_c$, for $c = a, b$, is $k_n$-sparse vector, then

$$||\Omega_c||_1 \leqslant \sqrt{k_n}||\Omega_c||_2 \leqslant \sqrt{k_n} \lambda_{\max}(\Omega) \leqslant C_{\mathsf{cov}} \sqrt{k_n}$$

and we could then set $C_{\mathsf{sparse}} = C_{\mathsf{cov}}$. Assumption 3.3 is a high-level condition, which assumes existence of initial estimators of $\gamma_a$ and $\gamma_b$ that converge at a fast enough rate. In the next section, we will see that Assumption 3.1 together with a stronger version of Assumption 3.2 are sufficient for Assumption 3.3 to be satisfied with a specific estimator that is efficient to compute. Finally, Assumption 3.4 is imposed to allow for estimation of the asymptotic variance $\breve{\Omega}_{ab}$. While Assumption 3.1 depends only on the correlation matrix $\Sigma$ without reference to the distribution of the radius $\xi$, Assumption 3.4 depends on both $\Sigma$ and $\xi$, and therefore cannot be derived as a consequence of the choice of $\Sigma$.

We now state our main result.

**Theorem 3.5.** *Under Assumptions 3.1, 3.2, 3.2, and 3.4, there exists a constant $C_{\mathsf{converge}}$, depending on $C_{\mathsf{cov}}, C_{\mathsf{sparse}}, C_{\mathsf{est}}, C_{\mathsf{kernel}}$ but not on the dimensions $(n, p_n, k_n)$ of the problem, such that*

$$\sup_{t \in \mathbb{R}} \left| \mathbb{P}\left\{ \sqrt{n} \cdot \frac{\breve{\Omega}_{ab} - \Omega_{ab}}{\breve{S}_{ab}} \leqslant t \right\} - \Phi(t) \right| \leqslant C_{\mathsf{converge}} \cdot \sqrt{\frac{k_n^2 \log^2(p_n)}{n}} + \frac{1}{p_n} + \delta_n .$$

---
[2] Here we use the positive semidefinite ordering on matrices, that is, $A \succeq B$ if $A - B \succeq 0$. Note that the second part of the inequality, $\Sigma_{h_1} \preceq \Sigma_h$, is always true by the law of total variance.

We note that the result holds uniformly over a large class of data generating processes satisfying Assumptions 3.1, 3.2, 3.3, and 3.4, which are relatively weak assumptions compared to much of the sparse estimation and inference literature; we emphasize that the result holds without requiring exact model selection or oracle properties, which hold only for restrictive sequences of data generating processes. For example, we do not require the "beta-min" condition (that is, a lower bound on $|\Omega_{ab}|$ for all true edges) or any incoherence conditions [Bühlmann and van de Geer, 2011], which may be implausible in practice. Instead of requiring perfect model selection, we only require estimation consistency as given in Assumption 3.3; our weaker assumptions would not be sufficient to guarantee model selection consistency.

As an immediate corollary, we see that the confidence interval constructed in (2.10) is asymptotically correct:

**Corollary 3.6.** *Under the assumptions and notation of Theorem 3.5, the $(1-\alpha)$-confidence interval constructed in (2.10) fails to cover the true parameter $\Omega_{ab}$ with probability no higher than*

$$\alpha + 2\left[C_{\mathsf{converge}} \cdot \sqrt{\frac{k_n^2 \log^2(p_n)}{n}} + \frac{1}{p_n} + \delta_n\right].$$

Again this result holds uniformly over a large class of data generating distributions.

Theorem 3.5 is striking as it shows that we can form an asymptotically normal estimator of $\Omega_{ab}$ under the transelliptical distribution family with the sample complexity $n = \Omega\left(k_n^2 \log^2(p_n)\right)$. This sample size requirement was shown to be optimal for obtaining an asymptotically normal estimator of an element in a precision matrix from multivariate normal data [Ren et al., 2013]. More precisely, let[3]

$$\mathcal{G}_0(c_0, c_1, k_n) = \left\{ \begin{array}{l} \Omega = (\Omega_{ab})_{a,b \in [p_n]} \ : \ \max_{a \in [p_n]} \sum_{b \neq a} \mathbb{I}\{\Omega_{ab} \neq 0\} \leqslant k_n, \\ \text{and } c_0 \leqslant \lambda_{\min}(\Omega) \leqslant \lambda_{\max}(\Omega) \leqslant c_1. \end{array} \right\}$$

where $c_0, c_1 > 0$ are constants. Then Theorem 1 in Ren et al. [2013] proves

$$\inf_{a,b} \inf_{\breve{\Omega}_{ab}} \sup_{\mathcal{G}_0(c_0,c_1,k_n)} \mathbb{P}\left\{\left|\breve{\Omega}_{ab} - \Omega_{ab}\right| \geqslant \epsilon_0 \left(n^{-1} k_n \log(p_n) \vee n^{-1/2}\right)\right\} \geqslant \epsilon_0$$

and, therefore, our estimator is rate optimal in terms of the sample size scaling. (Above, the infimum is taken over any estimator $\breve{\Omega}_{ab}$ which is a measurable function of the data.) We can also consider a related optimality question: whether the confidence interval we produce has the optimal (that is, lowest possible) width, given the desired coverage level. In the Gaussian setting, Ren et al. [2013]'s method produces an interval which has asymptotically minimal length at the given sample size, due to the fact that the variance of their estimator matches the Fisher information. Our ROCKET method does not enjoy this theoretical property, but empirically we observe that our confidence intervals are only slightly wider than those produced by Ren et al. [2013]'s method, for Gaussian data.

At this point, it is also worth mentioning the result of Gu et al. [2015], who study inference under Gaussian copula graphical models. They base their inference procedure on decorrelating a pseudo score function for the parameter of interest and showing that it is normally distributed. Their main result, stated in Theorem 4.10, requires the sample size to satisfy

$$k_n^3 M_n^6 \left(\frac{\log(p_n)}{n}\right)^{3/2} + k_n^2 M_n^3 \frac{\log(p_n)}{n} = o\left(n^{-1/2}\right)$$

where $M_n = \max_{a \in [p_n]} \sum_{b \in [p_n]} |\Omega_{ab}|$. As $M_n$ can be potentially as large as $\sqrt{k_n}$, it is immediately clear that our result achieves much better scaling on the sample size.

---

[3]In their work, the constants $c_0, c_1$ are instead denoted by a constant $M \geqslant c_0^{-1} \vee c_1$; we use different notation here to distinguish from the $M$ used in Gu et al. [2015] which plays a very different role, and which we denote by $M_n$ as it is not necessarily constant.

## 3.1 Initial estimators

The validity of our inference method relies in part on the accuracy of the initial estimators $\breve{\gamma}_a$ and $\breve{\gamma}_b$, which are assumed to satisfy error bounds with high probability as stated in Assumption 3.3—that is, with high probability, we have

$$||\breve{\gamma}_c - \gamma_c||_2 \leq C_{\mathsf{est}}\sqrt{\frac{k_n \log(p_n)}{n}}, \ ||\breve{\gamma}_c - \gamma_c||_1 \leq C_{\mathsf{est}}\sqrt{\frac{k_n^2 \log(p_n)}{n}} \ ,$$

for $c = a, b$, where $C_{\mathsf{est}}$ is some constant. Below, we will prove that these required error rates can be obtained, under an additional sparsity assumption, by the Lasso estimators

$$\breve{\gamma}_c = \underset{\gamma \in \mathbb{R}^I ; ||\gamma||_1 \leq C_{\mathsf{cov}}\sqrt{2k_n}}{\operatorname{argmin}} \left\{ \frac{1}{2} \gamma^\top \hat{\Sigma}_I \gamma - \gamma^\top \hat{\Sigma}_{Ic} + \lambda ||\gamma||_1 \right\} \tag{3.2}$$

for each $c = a, b$, when the penalty parameter $\lambda$ is chosen appropriately. In fact, these optimization problems may not be convex, because $\hat{\Sigma}_I$ will not necessarily be positive semidefinite.

We now turn to proving that any local minima for (3.2) for $c = a, b$ will satisfy the required error rates of Assumption 3.3. To proceed, we will use the theoretical results of Loh and Wainwright [2013], which gives a theory for local minimizers of nonconvex regularized objective functions. In particular, any local minimizers of the two optimization problems will satisfy requirements of Assumption 3.3 and, therefore, we only need to be able to run optimization algorithms that find local minima. We specialize their main result to our setting.

**Theorem 3.7** (Adapted from Loh and Wainwright [2013, Theorem 1]). *Consider any $n, p \geq 1$, any $A \in \mathbb{R}^{p \times p}$ and $z \in \mathbb{R}^p$, and any $k$-sparse $x^\star \in \mathbb{R}^p$ with $||x^\star||_1 \leq R$. Suppose that $A$ satisfies restricted strong convexity conditions*

$$v^\top A v \geq \alpha_1 ||v||_2^2 - \tau_1 ||v||_1^2 \cdot \frac{\log(p)}{n} \ . \tag{3.3}$$

*If*

$$n \geq \frac{16R^2 \tau_1 \max\{\alpha_1, \tau_1\} \log(p)}{\alpha_1^2} \tag{3.4}$$

*and*

$$\max\left\{ 4||Ax^\star - z||_\infty, 4\alpha_1 \sqrt{\frac{\log(p)}{n}} \right\} \leq \lambda \leq \frac{\alpha_1}{6R} \tag{3.5}$$

*then for any $\breve{x}$ that is a local minimum of the objective function $\frac{1}{2} x^\top A x - x^\top z + \lambda ||x||_1$ over the set $\{x \in \mathbb{R}^d : ||x||_1 \leq R\}$, it holds that*

$$||\breve{x} - x^\star||_2 \leq \frac{1.5\lambda\sqrt{k}}{\alpha_1} \ and \ ||\breve{x} - x^\star||_1 \leq \frac{6\lambda k}{\alpha_1} \ .$$

We apply Loh and Wainwright [2013]'s results, Theorem 3.7, to our problem of estimating $\gamma_a$ and $\gamma_b$ in a setting where we assume exact sparsity. (It is likely that similar results would hold for approximate sparsity, but here we use exact sparsity to fit the assumptions of this existing theorem.)

**Corollary 3.8.** *Suppose that Assumption 3.1 holds. Assume additionally that the columns $\Omega_a, \Omega_b$ of the true inverse covariance $\Omega = \Sigma^{-1}$ are $k_n$-sparse. Then there exist constants $C_{\mathsf{sample}}, C_{\mathsf{Lasso}}$, depending on $C_{\mathsf{cov}}$ but not on $(n, k_n, p_n)$, such that if $n \geq C_{\mathsf{sample}} k_n \log(p_n)$ then, with probability at least $1 - \frac{1}{2p_n}$, any local minimizer $\breve{\gamma}_a$ of the objective function*

$$\frac{1}{2} \gamma^\top \hat{\Sigma}_I \gamma - \gamma^\top \hat{\Sigma}_{Ia} + \lambda ||\gamma||_1$$

*over the set* $\{\gamma \in \mathbb{R}^I : ||\gamma||_1 \leq C_{\mathsf{cov}}\sqrt{2k_n}\}$ *satisfies*

$$||\breve{\gamma}_a - \gamma_a||_2 \leq 3\sqrt{2}C_{\mathsf{cov}}\lambda\sqrt{k_n} \text{ and } ||\breve{\gamma}_a - \gamma_a||_1 \leq 24C_{\mathsf{cov}}\lambda\sqrt{k_n} .$$

*where we choose* $\lambda = C_{\mathsf{Lasso}} \cdot \sqrt{\frac{\log(p_n)}{n}}$. *The same holds for estimating* $\gamma_b$.

Using this corollary, we see that a local minimizer of (3.2), $\breve{\gamma}_c$, satisfies Assumption 3.3 with $\delta_n = \frac{1}{p_n}$ and $C_{\mathsf{est}} = 24C_{\mathsf{cov}}C_{\mathsf{Lasso}}$. We remark that in practice, the constant $C_{\mathsf{Lasso}}$ suggested by the theory is in general unknown, but choosing $\lambda$ to be a small multiple of $\sqrt{\log(p_n)/n}$ generally performs well—for instance $\lambda = 2.1\sqrt{\log(p_n)/n}$ as we use in our simulations, where the choice of the constant 2.1 ensures that the penalty term dominates the variance of the elements of the objective function's derivative, that is, the elements of $\widehat{\Sigma}_I\gamma - \widehat{\Sigma}_{Ic}$, at the true solution $\gamma = \gamma_c$.

To prove that this corollary follows from Loh and Wainwright [2013]'s result (Theorem 3.7), it is sufficient to check that the restricted strong convexity condition (3.3) holds with high probability for the matrix $\widehat{\Sigma}_I$, and then compute the necessary values for $\lambda$ and the other parameters of Theorem 3.7. The proof is technical and relies on novel results on concentration of the Kendall's tau correlation matrix. Details are given in Appendix E.

We have provided sufficient condition for a local minimizer of (3.2) to satisfy Assumption 3.3, however, many other estimators can be used as initial estimators. For example, one could use the Dantzig selector [Candès and Tao, 2007]. Potential benefits of the Dantzig selector over the optimization program in (3.2) are twofold. First, the optimization program is convex even when $\widehat{\Sigma}_I$ is not positive semi-definite. Second, one does not need to know an upper bound $R$ on the $\ell_1$ norm of $\Omega_c$ for $c = a, b$. Using the techniques similar to those used to prove Corollary 3.8, we can also prove that Assumption 3.3 holds when the Dantzig selector is used as an initial estimator. For large problems, however, Dantzig selector type methods are computationally much slower than Lasso type methods; in our empirical results, we implement the Lasso rather than the Dantzig selector since we study graphs with as many as 1000 nodes.

In practice, we have found that in simulations, using the Lasso for model selection, and then refitting without a penalty, leads to better empirical performance. Specifically, for each $c = a, b$, we first fit

$$\breve{\gamma}_c^{\mathsf{Lasso}} = \underset{\gamma \in \mathbb{R}^I}{\mathsf{argmin}} \left\{ \frac{1}{2}\gamma^\top\widehat{\Sigma}_I\gamma - \gamma^\top\widehat{\Sigma}_{Ia} + \lambda||\gamma||_1 \right\};$$

or, more precisely, find a local minimum of this nonconvex optimization problem over the ball $\{\gamma : ||\gamma||_1 \leq R\}$ for a large radius $R$. (In practice, every iteration will stay inside this ball; therefore, as long as we see convergence in our iterative algorithm for solving this nonconvex Lasso, we do not concern ourselves with this theoretical boundedness constraint.)

We then extract the combined support of the two solutions, $\breve{J} = \mathsf{supp}(\breve{\gamma}_a^{\mathsf{Lasso}}) \cup \mathsf{supp}(\breve{\gamma}_b^{\mathsf{Lasso}})$, and refit the coefficients using least-squares:

$$\breve{\gamma}_c = \left(\widehat{\Sigma}_{\breve{J}}\right)^{-1}\widehat{\Sigma}_{\breve{J}c} \text{ for } c = a, b .$$

Following the work of Belloni and Chernozhukov [2013] or Sun and Zhang [2012a], it can be shown that the refitted estimators also satisfy the Assumption 3.3; in practice, refitting improves the accuracy of these preliminary estimators by reducing shrinkage bias.

Finally, we remark that if we would like to perform inference for all $\binom{p_n}{2}$ potential edges, then we require $2 \cdot \binom{p_n}{2} \sim p_n^2$ many initial estimators to be computed; this is of course quite computationally demanding. However, Ren et al. [2013] propose a simple modification that significantly reduces computation time: for each node $a$ we can first regress $X_a$ on all the other variables; call this solution $\breve{\gamma}_a^{\mathsf{all}}$. Next for any $b \neq a$, if $(\breve{\gamma}_a^{\mathsf{all}})_b = 0$, then this solution $\breve{\gamma}_a^{\mathsf{all}}$ is already optimal for regressing node $a$ on nodes $I = [p_n]\backslash\{a, b\}$; this will be the case for most nodes $b$ due to sparsity. With this modification, the actual number of regressions required is far smaller—if each node $a$ forms edges with at most $k_n$ other nodes (that is, $\breve{\gamma}_a^{\text{all}}$ is $k_n$-sparse), then we will require only $p_n(k_n + 1)$ many regressions in total to form all of the initial estimators.

## 4 Main technical tools

In this section, we outline the proof of Theorem 3.5 (Section 4.1) and state the key technical result that establishes the sign-subgaussianity property of a vector $X$ following a transelliptical distribution (Section 4.2). We also illustrate an application of this technical result to establishing a bound on $\widehat{\Sigma} - \Sigma$ (Section 4.3).

### 4.1 Sketch of proof for main result

The proof of Theorem 3.5 has two key steps. First, in Step 1, we prove that the distribution of $\widetilde{\Theta}_{ab}$, the oracle estimator of $\Theta_{ab}$, is asymptotically normal, with

$$\sqrt{n} \cdot \frac{\widetilde{\Theta}_{ab} - \Theta_{ab}}{S_{ab} \det(\Theta)} \to N(0, 1)$$

where $S_{ab}$ is the asymptotic variance of $\widetilde{\Omega}_{ab}$. (Explicit form of $S_{ab}$ is given in the proof of Theorem 4.1.) Next, in Step 2, we prove that the difference between the estimator and the oracle estimator, $\breve{\Theta} - \widetilde{\Theta}$, converges to zero at a fast rate, and that the variance estimator $\breve{S}_{ab}$ converges to $S_{ab}$ at a fast rate. Combining these steps we prove that $\breve{\Omega}_{ab}$ is an asymptotically normal estimator of $\Omega_{ab}$. The detailed proofs for each step are found in Appendix D. Here, we outline the main results for each step.

Step 1 establishes the Berry-Esseen type bound for the centered and normalized oracle estimator $\sqrt{n} \cdot \frac{\widetilde{\Theta}_{ab} - \Theta_{ab}}{S_{ab} \cdot \det(\Theta)}$. We approximate the oracle estimator $\widetilde{\Theta}_{ab}$ by a linear function of the Kendall's tau statistic $\widehat{T}$, which is a U-statistic of the data. We prove that the variance of the linear approximation is bounded away from zero and apply existing results on convergence of U-statistics. The following result is proved in Appendix D

**Theorem 4.1.** *Suppose that Assumptions 3.1, 3.2, and 3.4 hold. Then there exist constants $C_{\text{normal}}, C_{\text{variance}}$ depending on $C_{\text{cov}}, C_{\text{sparse}}, C_{\text{kernel}}$ but not on $(n, p_n, k_n)$, such that*

$$\sup_{t \in \mathbb{R}} \left| \mathbb{P}\left\{ \sqrt{n} \cdot \frac{\widetilde{\Theta}_{ab} - \Theta_{ab}}{S_{ab} \cdot \det(\Theta)} \leqslant t \right\} - \Phi(t) \right| \leqslant C_{\text{normal}} \cdot \frac{k_n \log(p_n)}{\sqrt{n}} + \frac{1}{2p_n},$$

*where $S_{ab}$ is defined in the proof and satisfies $S_{ab} \cdot \det(\Theta) \geqslant C_{\text{variance}} > 0$.*

Step 2 contains the main challenge of this problem, since it requires strong results on the concentration properties of the Kendall's tau estimator $\widehat{\Sigma}$ of the covariance matrix $\Sigma$. The main ingredient for this step is a new result on "sign-subgaussianity", that is, proving that the signs vector $\text{sign}(X_i - X_{i'})$ is subgaussian for i.i.d. observations $X_i, X_{i'}$. Our results on sign-subgaussianity are discussed in Section 4.2 and their application to concentration of $\widehat{\Sigma}$ around $\Sigma$ is given in Section 4.3. Using these tools, we are able to prove the following theorem (proved in Appendix D):

**Theorem 4.2.** *Suppose that Assumptions 3.1, 3.2, and 3.3 hold. Then there exists a constant $C_{\text{oracle}}$, depending on $C_{\text{cov}}, C_{\text{sparse}}, C_{\text{est}}$ but not on $(n, p_n, k_n)$, such that, if[4] $n \geqslant 15 k_n \log(p_n)$, then, with probability at least $1 - \frac{1}{2p_n} - \delta_n$,*

$$||\breve{\Theta} - \widetilde{\Theta}||_\infty \leqslant C_{\text{oracle}} \cdot \frac{k_n \log(p_n)}{n}$$

---
[4]Note that the additional condition $n \geqslant 15 k_n \log(p_n)$ can be assumed to hold in our main result Theorem 3.5, since if this inequality does not hold, then the claim in Theorem 3.5 is trivial.

*and*

$$\left|\breve{S}_{ab} \cdot \det(\breve{\Theta}) - S_{ab} \cdot \det(\Theta)\right| \leq C_{\text{oracle}} \cdot \sqrt{\frac{k_n^2 \log(p_n)}{n}}.$$

## 4.2 Sign-subgaussian random vectors

Recall the definition of a subgaussian random vector:

**Definition 4.3.** *A random vector $X \in \mathbb{R}^p$ is $C$-subgaussian if, for any fixed vector $v \in \mathbb{R}^p$, it holds that $\mathbb{E}\left[e^{v^\top X}\right] \leq e^{C \cdot ||v||_2^2/2}$.*

For graphical models where the data points $X_i$ come from a subgaussian distribution, the sample covariance matrix $\frac{1}{n}\sum_i (X_i - \overline{X})(X_i - \overline{X})^\top$, with $\overline{X} = \frac{1}{n}\sum_i X_i$, is known to concentrate near the population covariance, as measured by different norms. For example, elementwise convergence of the sample covariance to the population covariance, that is, convergence in $||\cdot||_\infty$, is sufficient to establish rates of convergence for the graphical Lasso, CLIME or graphical Dantzig selector for estimating the sparse inverse covariance [Ravikumar et al., 2011, Cai et al., 2011, Yuan, 2010]. Similar results can be obtained also for the transelliptical family, since $||\widehat{T} - T||_\infty \leq C\sqrt{\log(p)/n}$ and hence $||\widehat{\Sigma} - \Sigma||_\infty \leq C\sqrt{\log(p)/n}$, as was shown in Liu et al. [2012a] and Liu et al. [2012b]. However, in order to construct asymptotically normal estimators for the elements of the precision matrix, stronger results are needed about the convergence of the sample covariance to the population covariance [Ren et al., 2013]. In particular, a result on convergence in spectral norm, uniformly over all sparse submatrices, is required. One can relate the convergence in the elementwise $\ell_\infty$ norm to (sparse) spectral norm convergence, however, this would lead to suboptimal sample size. One way to obtain a tight bound on the (sparse) spectral norm convergence is by utilizing subgaussianity of the data points $X_i$. This is exactly what we proceed to establish.

Recall from (2.3) the Kendall's tau estimator of the covariance,

$$\widehat{\Sigma} = \sin\left(\frac{\pi}{2}\widehat{T}\right) \text{ where } \widehat{T} = \frac{1}{\binom{n}{2}} \sum_{i<i'} \text{sign}(X_i - X_{i'})\text{sign}(X_i - X_{i'})^\top.$$

Therefore, it is crucial to determine whether the vector $\text{sign}(X_i - X_{i'})$ is itself subgaussian, at a scale that does not depend heavily on the ambient dimension $p_n$.[5] Using past results on elliptical distributions, we can reduce to a simpler case using the arguments of Lindskog et al. [2003] (proved in Appendix F):

**Lemma 4.4.** *Let $X, X' \stackrel{iid}{\sim} \mathsf{TE}(\Sigma, \xi; f_1, \ldots, f_p)$. Suppose that $\Sigma$ is positive definite, and that $\xi > 0$ with probability $1$. Then $\text{sign}(X - X')$ is equal in distribution to $\text{sign}(Z)$, where $Z \sim N(0, \Sigma)$.*

Previous work has shown that a Gaussian random vector $Z \sim N(0, \Sigma)$ is "sign-subgaussian", that is, $\text{sign}(Z)$ is subgaussian with variance proxy that depends on $p_n$ only through $\mathsf{C}(\Sigma)$, for special cases when the covariance $\Sigma$ is identity or equicorrelation matrix [Han and Liu, 2013]. However, a result for general covariance structures was previously unknown.

In the following lemma, we resolve this question, proving that Gaussian vectors are sign-subgaussian (recall $\mathsf{C}(\Sigma)$ is the condition number of $\Sigma$):

**Lemma 4.5.** *Let $Z \sim N(\mu, \Sigma)$. Then $\text{sign}(Z) - \mathbb{E}\left[\text{sign}(Z)\right]$ is $\mathsf{C}(\Sigma)$-subgaussian.*

This lemma is the primary tool for our main results in this paper—specifically, it is the key ingredient to the proof of Theorem 4.2, which bounds the errors $\breve{\Theta} - \widetilde{\Theta}$ and $\breve{S}_{ab} \cdot \det(\breve{\Theta}) - S_{ab} \cdot \det(\Theta)$. Lemma 4.5 is proved in Appendix C. We also use this result in establishing results in the following section.

---

[5] Note that $v^\top \text{sign}(X_i - X_{i'})$ is obviously subgaussian for any distribution on $X$, as it is a sum of subgaussian random variables (since $\text{sign}(\cdot)$ is bounded), however, its scale could grow linearly with $p_n$.

## 4.3 Deterministic and probabilistic bounds on $\widehat{\Sigma} - \Sigma$

Lemma 4.5 is instrumental in obtaining probabilistic bounds on $\widehat{\Sigma} - \Sigma$. Results given in this section are crucial for establishing Theorem 4.2 and Corollary 3.8.

Let $\mathcal{S}_k$ be the set of $k$-sparse vectors in the unit ball,

$$\mathcal{S}_k = \{u \in \mathbb{R}^p : ||u||_2 \leq 1, ||u||_0 \leq k\} \ ,$$

and abusing notation, let $||\cdot||_{\mathcal{S}_k}$ denote the sparse spectral norm for matrices, that is, $||M||_{\mathcal{S}_k} = \max_{u,v \in \mathcal{S}_k} u^\top M v$.

The following lemma provides a bound on the error in Kendall's tau, that is, on $\widehat{T} - T$, in this sparse sectral norm (with the proof given in Appendix F).

**Lemma 4.6.** *Suppose that $k \geq 1$ and $\delta \in (0, 1)$ satisfy $\log(2/\delta) + 2k\log(12p) \leq n$. Then with probability at least $1 - \delta$ it holds that*

$$||\widehat{T} - T||_{\mathcal{S}_k} \leq 32(1 + \sqrt{5})\mathsf{C}(\Sigma) \cdot \sqrt{\frac{\log(2/\delta) + 2k\log(12p)}{n}} \ .$$

Next, we relate $\widehat{\Sigma}$ to $\widehat{T}$, with the following deterministic bound on the sparse spectral norm of the error of the covariance estimator $\widehat{\Sigma}$, which is proven in Appendix F:

**Lemma 4.7.** *The following bound holds deterministically: for any $k \geq 1$,*

$$||\widehat{\Sigma} - \Sigma||_{\mathcal{S}_k} \leq \frac{\pi^2}{8} \cdot k ||\widehat{T} - T||_\infty^2 + 2\pi ||\widehat{T} - T||_{\mathcal{S}_k} \ . \tag{4.1}$$

A result in de la Pena and Giné [1999, Theorem 4.1.8] bounds $||\widehat{T} - T||_\infty$ with high probability (details of this bound are given in Appendix D). Combining the bound on $||\widehat{T} - T||_\infty$ with Lemmas 4.6 and 4.7, we immediately obtain the following corollary:

**Corollary 4.8.** *Take any $\delta_1, \delta_2 \in (0, 1)$ and any $k \geq 1$ such that $\log(2/\delta_2) + 2k\log(12p) \leq n$. Then, with probability at least $1 - \delta_1 - \delta_2$, the following bound on $\widehat{\Sigma} - \Sigma$ holds:*

$$||\widehat{\Sigma} - \Sigma||_{\mathcal{S}_k} \leq$$

$$\frac{\pi^2}{8} \cdot k \cdot \frac{4\log\left(2\binom{p}{2}/\delta_1\right)}{n} + 2\pi \cdot 32(1 + \sqrt{5})\mathsf{C}(\Sigma) \cdot \sqrt{\frac{\log(2/\delta_2) + 2k\log(12p)}{n}} \ . \tag{4.2}$$

Finally, we use a result based on the work of Sun and Zhang [2012a], in order to extend this sparse spectral norm bound to a bound holding for all approximately sparse vectors $u$ and $v$:

**Lemma 4.9** (Based on Proposition 5 of Sun and Zhang [2012a])**.** *For any fixed matrix $M \in \mathbb{R}^{p \times p}$ and vectors $u, v \in \mathbb{R}^p$, and any $k \geq 1$,*

$$|u^\top M v| \leq \left(||u||_2 + ||u||_1/\sqrt{k}\right) \cdot \left(||v||_2 + ||v||_1/\sqrt{k}\right) \cdot ||M||_{\mathcal{S}_k}.$$

Results of Lemma 4.6 and Corollary 4.8 can be compared to Theorem 2 in Mitra and Zhang [2014], which proves essentially the same result for the Kendall's tau estimate of $\Sigma$, but only for the nonparanormal (Gaussian copula) model; their technique does not extend immediately to the transelliptical model. When $C(\Sigma) = O(1)$, we extend their result to the transelliptical model and, as a special case, this provides an alternative proof for their result on the Gaussian copula model. We note that their result does not depend on the condition number of the covariance matrix, but only on the maximum eigenvalue. However, in the context of graphical models it is commonly assumed that the smallest eigenvalue is a constant. Furthermore, our results in Lemma 4.6 and Corollary 4.8 can also be compared with Theorem 4.10 of Han and Liu [2013], which give similar bounds on the spectral norm of sparse submatrices of $\widehat{T} - T$ and $\widehat{\Sigma} - \Sigma$, but with a sign-subgaussianity assumption on the distribution. We rigorously establish the same bounds for all well-conditioned covariance matrices, without explicitly making the sign-subgaussian assumption.

# 5 Simulation studies

In this section, we illustrate finite sample properties of ROCKET described in Section 2 on simulated data. (A real data experiment, and some additional simulations, are presented in Appendix A.)

We use ROCKET to construct confidence intervals for edge parameters and report empirical coverage probabilities as well as the length of constructed intervals. For comparison, we also construct confidence intervals using the procedure of Ren et al. [2013], which is based on the Pearson correlation matrix, a nonparanormal estimator of the correlation matrix (NPN) proposed in Liu et al. [2009], and the pseudo score procedure of Gu et al. [2015]. For the first two methods, we use the plugin estimate of the correlation matrix together with (2.8) to estimate $\Omega_{ab}$. Recall that Liu et al. [2009] estimate the correlation matrix based on the marginal transformation of the observed data. Let

$$\widetilde{F}_a(x) = \begin{cases} \delta_n & \text{if } \widehat{F}_a(x) < \delta_n \\ \widehat{F}_a(x) & \text{if } \delta_n \leqslant \widehat{F}_a(x) \leqslant 1 - \delta_n \\ 1 - \delta_n & \text{if } \widehat{F}_a(x) > 1 - \delta_n, \end{cases}$$

where $\widehat{F}_a(x) = n^{-1}\sum_i \mathbb{1}\{X_{ia} < x\}$ is the empirical CDF of $X_a$ and $\delta_n = \left(4n^{1/4}\sqrt{\pi \log(n)}\right)^{-1}$. The correlation matrix $\widehat{\Sigma} = \left(\widehat{\Sigma}_{ab}\right)_{ab}$ is then estimated as $\widehat{\Sigma}_{ab} = \widehat{\text{Corr}}\left(\Phi\left(\widetilde{F}_a(X_{ia})\right), \Phi\left(\widetilde{F}_b(X_{ib})\right)\right)$. Asymptotic variance of estimators of $\Omega_{ab}$ based on the Pearson or nonparanormal correlation matrix is obtained as $\check{S}_{ab}^2 = n^{-1}\left(\check{\Omega}_{aa}\check{\Omega}_{bb} + \check{\Omega}_{ab}^2\right)$. Gu et al. [2015] estimate $\Omega_{ab}$ as

$$\check{\Omega}_{ab}^{\text{PS}} = \frac{\widehat{\Omega}_{ab}\left(\left(\widehat{\Omega}\widehat{\Sigma}\right)_{ab} + \left(\widehat{\Sigma}\widehat{\Omega}\right)_{ab}\right) - \left(\widehat{\Omega}\widehat{\Sigma}\widehat{\Omega}\right)_{ab}}{\left(\widehat{\Omega}\widehat{\Sigma}\right)_{ab} + \left(\widehat{\Sigma}\widehat{\Omega}\right)_{ab} - 1},$$

where $\widehat{\Sigma}$ is the Kendall's tau estimator of the covariance matrix in (2.3) and $\widehat{\Omega}$ is an initial estimator of the precision matrix. Under suitable conditions, $\widehat{\Omega}_{ab}^{\text{PS}}$ is asymptotically normal with the asymptotic variance that can be consistently estimated as in Corollary 4.12 of Gu et al. [2015]. Gu et al. [2015] suggest using the CLIME estimator [Cai et al., 2011] to construct $\widehat{\Omega}$, however, we find that empirically the method performs better using lasso-with-refitting to estimate each row of $\Omega$, similar to Sun and Zhang [2012a]. For all simulations, we set the tuning parameter $\lambda = 2.1\sqrt{\log(p_n)/n}$, as suggested by our theory—this constant is large enough so that the penalty dominates the variance of each element of the score. All computations are carried out in Matlab.

**Simulation 1.** We generate data from the model $X \sim \mathsf{E}(0, \Sigma, \xi)$, where $\xi$ follows a $t$-distribution with 5 degrees of freedom. The inverse covariance matrix $\Omega$ encodes a grid where each node is connected to its four nearest neighbors with the nonzero elements of $\Omega^0$ equal to $\omega = 0.24$. Diagonal element of $\Omega^0$ are equal to 1. Let $(\Omega^0)^{-1} = \Sigma^0$. Then set $\Sigma = (\text{diag}(\Sigma^0))^{-1/2}\Sigma^0(\text{diag}(\Sigma^0))^{-1/2}$ and $\Omega = \Sigma^{-1}$. (Additional simulations in Appendix B show the same experiment on a chain graph structure.)

We take a grid of size $30 \times 30$ (so that $p_n = 900$) and take sample size $n = 400$. Figure 2 shows quantile-quantile (Q-Q) plots based on 1000 independent realizations of the test statistic error, $\sqrt{n} \cdot \frac{\check{\Omega}_{ab} - \Omega_{ab}}{\check{S}_{ab}}$, for the four methods together with the reference line showing quantiles of the standard normal distribution. From this figure, we observe that the quantiles of the test statistic error $\sqrt{n}\frac{\check{\Omega}_{ab} - \Omega_{ab}}{\check{S}_{ab}}$ based on ROCKET is closest to the quantiles of the standard normal random variable. We further quantify these results in Table 1, which reports empirical coverage and width of the confidence intervals based on $\sqrt{n}\frac{\check{\Omega}_{ab} - \Omega_{ab}}{\check{S}_{ab}}$. From the table, we can observe that the coverage of the confidence intervals based on ROCKET and the pseudo score are closest to nominal coverage of 95%. The three node pairs displayed in this figure and table, namely $\omega_{(2,2),(2,3)}, \omega_{(2,2),(3,3)}, \omega_{(2,2),(10,10)}$, correspond to a true edge, a non-edge between nearby nodes that is therefore easy to mistake for an edge, and a non-edge between distant nodes, respectively.

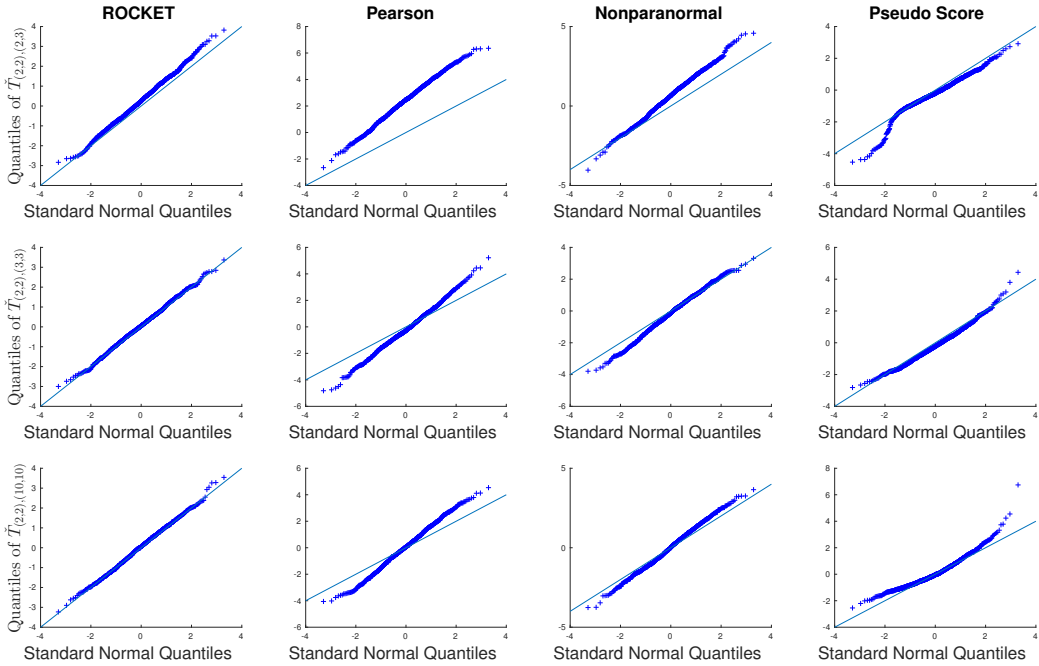

Figure 2: Simulation 1 (transelliptical data). Q-Q plot of $\sqrt{n} \cdot \frac{\breve{\Omega}_{ab} - \Omega_{ab}}{\breve{S}_{ab}}$ when $\Omega$ corresponds to a grid graph structure. Row 1 corresponds to an edge, row 2 to a close non-edge, and row 3 to a far non-edge.

|  | ROCKET | Pearson | NPN | Pseudo Score |
| --- | --- | --- | --- | --- |
| $\omega_{(2,2),(2,3)} = 0.37$ | 94.6 (0.51) | 36.6 (0.88) | 82.4 (0.48) | 92.2 (0.52) |
| $\omega_{(2,2),(3,3)} = 0$ | 94.3 (0.53) | 81.0 (0.86) | 88.3 (0.47) | 94.8 (0.50) |
| $\omega_{(2,2),(10,10)} = 0$ | 94.9 (0.56) | 78.3 (0.88) | 89.1 (0.48) | 95.5 (0.53) |

Table 1: Simulation 1 (transelliptical data). Percent empirical coverage (average length) of 95% confidence intervals based on 1000 independent simulation runs.

These results are not surprising, since neither the Pearson nor the nonparanormal correlation matrix consistently estimate the true $\Sigma$. In contrast, both ROCKET and the pseudo score method are able to construct a test statistic $\sqrt{n} \cdot \frac{\breve{\Omega}_{ab}}{\breve{S}_{ab}}$ that is asymptotically distributed as a normal random variable. The asymptotic distribution provides a good approximation to the finite sample distribution of $\sqrt{n} \cdot \frac{\breve{\Omega}_{ab} - \Omega_{ab}}{\breve{S}_{ab}}$.

**Simulation 2.** We illustrate performance of ROCKET when data are generated from a normal and nonparanormal distribution. We consider $\Omega$ corresponding to a grid as in Simulation 1, and generate $n = 400$ samples from $N(0, \Omega^{-1})$ and $\mathsf{NPN}(\Omega^{-1}; \widetilde{f}_1, \ldots, \widetilde{f}_p)$, where $\widetilde{f}_j = f_{\mathrm{mod}(j-1,5)+1}$ with $f_1(x) = x$, $f_2(x) = \mathrm{sign}(x)\sqrt{|x|}$, $f_3(x) = x^3$, $f_4(x) = \Phi(x)$, $f_5(x) = \exp(x)$.

Table 2 summarizes results from the simulation. We observe that when data are multivariate normal all methods perform well, with ROCKET and the pseudo score having slightly wider intervals, but with similar coverage. When data are generated from a nonparanormal distribution, using the Pearson correlation in (2.8) results in confidence intervals that do not have nominal coverage due to the bias. In this setting, nonparanormal estimator, ROCKET and the pseudo score still have the correct nominal coverage. Note however that when Kendall's tau is equal to zero, Pearson correlation is also equal to zero, and coverage for Pearson improves. See, for example, coverage for $\omega_{(2,2),(3,3)}$ and $\omega_{(2,2),(10,10)}$.

|  | | ROCKET | Pearson | NPN | Pseudo Score |
| --- | --- | --- | --- | --- | --- |
| Gaussian | $\omega_{(2,2),(2,3)} = 0.37$ | 95.6 (0.33) | 94.1 (0.32) | 94.4 (0.32) | 95.9 (0.33) |
| | $\omega_{(2,2),(3,3)} = 0$ | 95.4 (0.35) | 95.9 (0.34) | 95.3 (0.34) | 95.0 (0.35) |
| | $\omega_{(2,2),(10,10)} = 0$ | 96.0 (0.36) | 95.3 (0.35) | 95.6 (0.35) | 94.8 (0.36) |
| Transf. Gaussian | $\omega_{(2,2),(2,3)} = 0.37$ | 95.1 (0.35) | 12.0 (0.33) | 94.9 (0.33) | 95.3 (0.35) |
| | $\omega_{(2,2),(3,3)} = 0$ | 95.3 (0.35) | 93.1 (0.32) | 94.9 (0.34) | 94.7 (0.34) |
| | $\omega_{(2,2),(10,10)} = 0$ | 93.7 (0.36) | 93.2 (0.32) | 94.2 (0.35) | 94.8 (0.39) |

Table 2: Simulation 2 (Gaussian and nonparanormal data). Percent empirical coverage (average length) of 95% confidence intervals based on 1000 independent simulation runs.

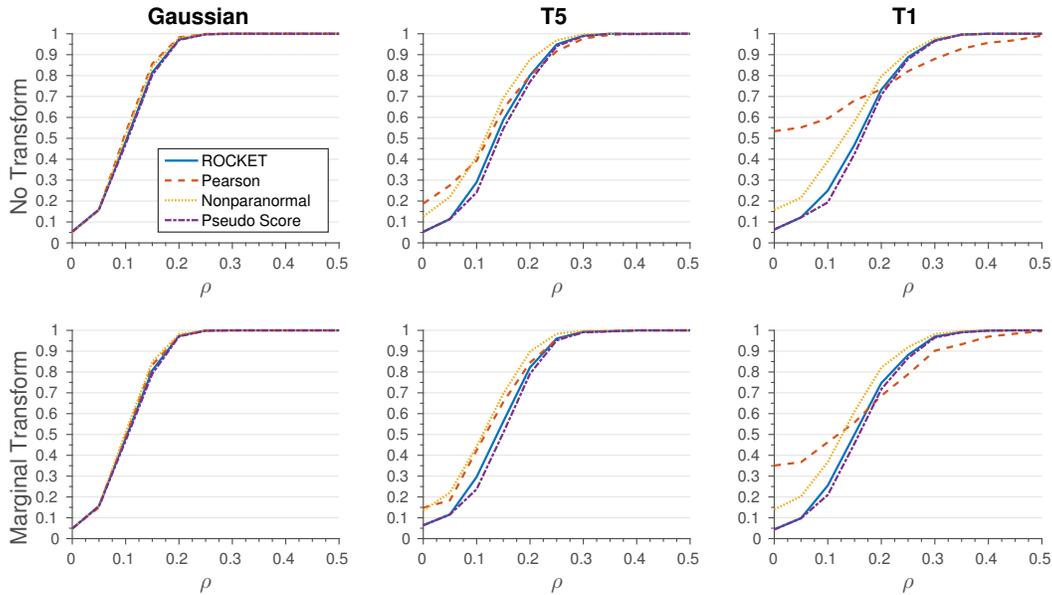

Figure 3: Simulation 3. Power plots for simulated data generated from a Gaussian distribution, and from a multivariate t distribution with 5 d.f. or with 1 d.f.

**Simulation 3.** In this simulation, we illustrate the power of a test based on the statistic $\sqrt{n} \cdot \frac{\check{\Omega}_{ab}}{\check{S}_{ab}}$ to reject the null hypothesis $H_{0,ab} : \Omega_{ab} = 0$. Samples are generated from the $\mathsf{E}(0, \Sigma, \xi)$ with $\xi$ having $\chi_{p_n}$, $t_5$, and $t_1$ distribution and the covariance matrix is of the form $\Sigma = I_P + E$ where $E_{12} = E_{21} = \rho$ and all other entries zero, with $p_n = 1000$ and $n = 400$. Note that $\xi \sim \chi_{p_n}$ implies that $X$ is multivariate normal. We also consider marginal transformation of $X$ as described in Simulation 2. Figure 3 plots empirical power curves based on 1000 independent simulation runs for different settings. When the data follow a normal distribution all methods have similar power. For other distributions, tests based on Pearson and nonparanormal correlation do not have correct coverage and are shown for illustrative purpose only.

# 6  Discussion

We have proposed a novel procedure ROCKET for inference on elements of the latent inverse correlation matrix under high-dimensional elliptical copula models. Our paper has established a surprising result, which states that ROCKET produces an asymptotically normal estimator for an element of the inverse correlation matrix in an elliptical copula model with the same sample complexity that is required to obtain an asymptotically normal estimator for an element in the precision matrix under a multivariate normal distribution. Furthermore, this sample complexity is optimal [Ren et al., 2013]. The result is surprising as the family of elliptical copula models is much larger than the family of multivariate normal distributions. For example, it contains distributions with heavy tail dependence as discussed in Section 2. ROCKET achieves the optimal requirement on the sample size without knowledge of the marginal transformation. Our result is also of significant practical importance. Since normal distribution is only a convenient mathematical approximation to data generating process, we recommend using ROCKET whenever making inference about inverse correlation matrix, instead of methods that heavily rely on normality. From simulation studies, even when data are generated from a normal distribution, ROCKET does not lose power compared to procedures that were specifically developed for inference under normality.

The main technical tool developed in the paper establishes that the sign of normal random vector, taken elementwise, is itself a sub-Gaussian random variable with the sub-Gaussian parameter depending on the condition number of the covariance matrix $\Sigma$ (but not on the dimension $p_n$). Based on this result, we were able to establish a tight tail bound on the deviation of sparse eigenvalues of the Kendall's tau matrix $\widehat{T}$. This result is of independent interest and it would allow us to improve a number of recent results on sparse principal component analysis, factor models and estimation of structured covariance matrices [Mitra and Zhang, 2014, Han and Liu, 2013, Fan et al., 2014]. The sharpest result on the nonparametric estimation of correlation matrices in spectral norm under a Gaussian copula model was established in [Mitra and Zhang, 2014]. Our results establish a similar result for the family of elliptical copula models and provide an alternative proof for the Gaussian copula model.

# Acknowledgments

This work is partially supported by an IBM Corporation Faculty Research Fund at the University of Chicago Booth School of Business, and an Alfred P. Sloan Fellowship. This work was completed in part with resources provided by the University of Chicago Research Computing Center.

# References

Richard G. Baraniuk, Mark A. Davenport, , and Michael B. Wakin. A simple proof of the restricted isometry property for random matrices. *Constructive Approximation*, 28(3):253–263, Jan 2008. ISSN 1432-0940. doi: 10.1007/s00365-007-9003-x. URL http://dx.doi.org/10.1007/s00365-007-9003-x.

Alexandre Belloni and Victor Chernozhukov. Least squares after model selection in high-dimensional sparse models. *Bernoulli'*, 19(2):521–547, May 2013. ISSN 1350-7265. doi: 10.3150/11-bej410. URL http://dx.doi.org/10.3150/11-BEJ410.

Alexandre Belloni, Victor Chernozhukov, and Christian B. Hansen. Inference on treatment effects after selection amongst high-dimensional controls. *Rev. Econ. Stud.*, 81(2):608–650, Nov 2013a. ISSN 1467-937X. doi: 10.1093/restud/rdt044. URL http://dx.doi.org/10.1093/restud/rdt044.

Alexandre Belloni, Victor Chernozhukov, and Kengo Kato. Robust inference in high-dimensional approximately sparse quantile regression models. *arXiv preprint arXiv:1312.7186*, December 2013b.

Alexandre Belloni, Victor Chernozhukov, and Kengo Kato. Uniform post selection inference for lad regression models. *arXiv preprint arXiv:1304.0282*, 2013c. URL http://arxiv.org/abs/1304.0282.

Alexandre Belloni, Victor Chernozhukov, and Ying Wei. Honest confidence regions for logistic regression with a large number of controls. *arXiv preprint arXiv:1304.3969*, 2013d. URL http://arxiv.org/abs/1304.3969.

Peter Bühlmann and Sara A. van de Geer. *Statistics for high-dimensional data*. Springer Series in Statistics. Springer, Heidelberg, 2011. ISBN 978-3-642-20191-2. doi: 10.1007/978-3-642-20192-9. URL http://dx.doi.org/10.1007/978-3-642-20192-9. Methods, theory and applications.

T. Tony Cai, W. Liu, and X. Luo. A constrained $\ell_1$ minimization approach to sparse precision matrix estimation. *J. Am. Stat. Assoc.*, 106(494):594–607, 2011.

Herman Callaert and Paul Janssen. The Berry-Esseen theorem for $U$-statistics. *Ann. Stat.*, 6(2): 417–421, 1978. ISSN 0090-5364.

Emmanuel J. Candès and T. Tao. The dantzig selector: Statistical estimation when $p$ is much larger than $n$. *Ann. Stat.*, 35(6):2313–2351, 2007. ISSN 0090-5364. doi: 10.1214/009053606000001523. URL http://dx.doi.org/10.1214/009053606000001523.

Mengjie Chen, Zhao Ren, Hongyu Zhao, and Harrison H. Zhou. Asymptotically normal and efficient estimation of covariate-adjusted gaussian graphical model. *arXiv preprint arXiv:1309.5923*, 2013. URL http://arxiv.org/abs/1309.5923.

Jie Cheng, Elizaveta Levina, and Ji Zhu. High-dimensional mixed graphical models. *ArXiv e-prints, arXiv:1304.2810*, April 2013.

Gabor Csardi and Tamas Nepusz. The igraph software package for complex network research. *InterJournal*, Complex Systems:1695, 2006. URL http://igraph.org.

Alexandre d'Aspremont, Onureena Banerjee, and Laurent El Ghaoui. First-order methods for sparse covariance selection. *SIAM J. Matrix Anal. Appl.*, 30(1):56–66, 2008. ISSN 0895-4798. doi: 10.1137/060670985. URL http://dx.doi.org/10.1137/060670985.

Victor de la Pena and Evarist Giné. *Decoupling: from dependence to independence*. Springer, 1999.

Paul Embrechts, Filip Lindskog, and Alexander McNeil. Modelling dependence with copulas and applications to risk management. In S. T. Rachev, editor, *Handbook of heavy tailed distributions in finance*, pages 329–384. Elsevier, 2003. URL http://books.google.com/books?hl=en&lr=&id=sv8jGSVFra8C&oi=fnd&pg=PA329&dq=Modelling+dependence+with+copulas+and+applications+to+risk+management&ots=Yvg5prGUXz&sig=b3Qs3l2OC8JgSQKPxUnY6Lq2Ql8.

Jianqing Fan, Y. Feng, and Y. Wu. Network exploration via the adaptive lasso and scad penalties. *Ann. Appl. Stat.*, 3(2):521–541, 2009. ISSN 1932-6157. doi: 10.1214/08-AOAS215. URL http://dx.doi.org/10.1214/08-AOAS215.

Jianqing Fan, Fang Han, and Han Liu. Page: Robust pattern guided estimation of large covariance matrix. Technical report, Technical report, Princeton University, 2014.

Kai Tai Fang, Samuel Kotz, and Kai Wang Ng. *Symmetric multivariate and related distributions*, volume 36 of *Monographs on Statistics and Applied Probability*. Chapman and Hall, Ltd., London, 1990. ISBN 0-412-31430-4. doi: 10.1007/978-1-4899-2937-2. URL http://dx.doi.org/10.1007/978-1-4899-2937-2.

Max H. Farrell. Robust inference on average treatment effects with possibly more covariates than observations. *arXiv preprint arXiv:1309.4686*, September 2013.

Jerome H. Friedman, Trevor J. Hastie, and Robert J. Tibshirani. Sparse inverse covariance estimation with the graphical lasso. *Biostatistics*, 9(3):432–441, 2008.

Q. Gu, Y. Cao, Y. Ning, and Han Liu. Local and global inference for high dimensional gaussian copula graphical models. *ArXiv e-prints, arXiv:1502.02347*, February 2015.

J. Guo, Elizaveta Levina, G. Michailidis, and J. Zhu. Joint estimation of multiple graphical models. *Biometrika*, 98(1):1–15, 2011a.

Jian Guo, Elizaveta Levina, George Michailidis, and Ji Zhu. Asymptotic properties of the joint neighborhood selection method for estimating categorical markov networks. Technical report, University of Michigan, 2011b.

Fang Han and Han Liu. Optimal rates of convergence for latent generalized correlation matrix estimation in transelliptical distribution. *ArXiv e-prints, arXiv:1305.6916*, May 2013.

Holger Höfling and Robert J. Tibshirani. Estimation of sparse binary pairwise markov networks using pseudo-likelihoods. *J. Mach. Learn. Res.*, 10:883–906, 2009. URL http://dl.acm.org/citation.cfm?id=1577101.

Adel Javanmard and Andrea Montanari. Nearly optimal sample size in hypothesis testing for high-dimensional regression. *arXiv preprint arXiv:1311.0274*, November 2013.

Adel Javanmard and Andrea Montanari. Confidence intervals and hypothesis testing for high-dimensional regression. *J. Mach. Learn. Res.*, 15(Oct):2869–2909, 2014. URL http://jmlr.org/papers/v15/javanmard14a.html.

Claudia Klüppelberg, Gabriel Kuhn, and Liang Peng. Semi-parametric models for the multivariate tail dependence function–the asymptotically dependent case. *Scand. J. Stat.*, 35(4):701–718, 2008. URL http://onlinelibrary.wiley.com/doi/10.1111/j.1467-9469.2008.00602.x/full.

C. Lam and Jianqing Fan. Sparsistency and rates of convergence in large covariance matrix estimation. *Ann. Stat.*, 37:4254–4278, 2009.

S. L. Lauritzen. *Graphical Models*, volume 17 of *Oxford Statistical Science Series*. The Clarendon Press Oxford University Press, New York, 1996. ISBN 0-19-852219-3. Oxford Science Publications.

Jason D. Lee and Trevor J. Hastie. Learning mixed graphical models. *ArXiv e-prints, arXiv:1205.5012*, May 2012.

Jason D. Lee, Dennis L. Sun, Yuekai Sun, and Jonathan E. Taylor. Exact post-selection inference with the lasso. *ArXiv e-prints, arXiv:1311.6238*, November 2013.

Filip Lindskog, Alexander McNeil, and Uwe Schmock. Kendall's tau for elliptical distributions. *Credit Risk*, pages 149–156, 2003. ISSN 1431-1933. doi: 10.1007/978-3-642-59365-9_8. URL http://dx.doi.org/10.1007/978-3-642-59365-9_8.

Han Liu and Lie Wang. Tiger: A tuning-insensitive approach for optimally estimating gaussian graphical models. *ArXiv e-prints, arXiv:1209.2437*, September 2012.

Han Liu, John D. Lafferty, and Larry A. Wasserman. The nonparanormal: Semiparametric estimation of high dimensional undirected graphs. *J. Mach. Learn. Res.*, 10:2295–2328, 2009.

Han Liu, Fang Han, Ming Yuan, John D. Lafferty, and Larry A. Wasserman. High-dimensional semiparametric Gaussian copula graphical models. *Ann. Stat.*, 40(4):2293–2326, 2012a. ISSN 0090-5364. doi: 10.1214/12-AOS1037. URL http://dx.doi.org/10.1214/12-AOS1037.

Han Liu, Fang Han, and Cun-Hui Zhang. Transelliptical graphical models. In *Proc. of NIPS*, pages 809–817. 2012b. URL http://books.nips.cc/papers/files/nips25/NIPS2012_0380.pdf.

Weidong Liu. Gaussian graphical model estimation with false discovery rate control. *Ann. Stat.*, 41 (6):2948–2978, 2013. ISSN 0090-5364. doi: 10.1214/13-AOS1169. URL http://dx.doi.org/10.1214/13-AOS1169.

Richard Lockhart, Jonathan E. Taylor, Robert J. Tibshirani, and Robert J. Tibshirani. A significance test for the lasso. *Ann. Stat.*, 42(2):413–468, 2014. ISSN 0090-5364. doi: 10.1214/13-AOS1175. URL http://dx.doi.org/10.1214/13-AOS1175.

Po-Ling Loh and Martin J. Wainwright. Regularized m-estimators with nonconvexity: Statistical and algorithmic theory for local optima. *arXiv preprint arXiv:1305.2436*, 2013. URL http://arxiv.org/abs/1305.2436.

Pascal Massart. *Concentration inequalities and model selection*, volume 1896 of *Lecture Notes in Mathematics*. Springer, Berlin, 2007. ISBN 978-3-540-48497-4; 3-540-48497-3. Lectures from the 33rd Summer School on Probability Theory held in Saint-Flour, July 6–23, 2003, With a foreword by Jean Picard.

Nicolas Meinshausen and Peter Bühlmann. High dimensional graphs and variable selection with the lasso. *Ann. Stat.*, 34(3):1436–1462, 2006.

Ritwik Mitra and Cun-Hui Zhang. Multivariate analysis of nonparametric estimates of large correlation matrices. *ArXiv e-prints, arXiv:1403.6195*, March 2014.

Thomas Peel, Sandrine Anthoine, and Liva Ralaivola. Empirical bernstein inequalities for u-statistics. In J.D. Lafferty, C.K.I. Williams, J. Shawe-Taylor, R.S. Zemel, and A. Culotta, editors, *Adv. Neural Inf. Process. Syst. 23*, pages 1903–1911. Curran Associates, Inc., 2010. URL http://papers.nips.cc/paper/4081-empirical-bernstein-inequalities-for-u-statistics.pdf.

R Core Team. *R: A Language and Environment for Statistical Computing*. R Foundation for Statistical Computing, Vienna, Austria, 2012. URL http://www.R-project.org/. ISBN 3-900051-07-0.

P. Ravikumar, Martin J. Wainwright, G. Raskutti, and B. Yu. High-dimensional covariance estimation by minimizing $\ell_1$-penalized log-determinant divergence. *Electron. J. Stat.*, 5:935–980, 2011.

Pradeep Ravikumar, Martin J. Wainwright, and J. D. Lafferty. High-dimensional ising model selection using $\ell_1$-regularized logistic regression. *Ann. Stat.*, 38(3):1287–1319, 2010. ISSN 0090-5364. doi: 10.1214/09-AOS691. URL http://dx.doi.org/10.1214/09-AOS691.

Zhao Ren, Tingni Sun, Cun-Hui Zhang, and Harrison H. Zhou. Asymptotic normality and optimalities in estimation of large gaussian graphical model. *arXiv preprint arXiv:1309.6024*, 2013. URL http://arxiv.org/abs/1309.6024.

Adam J. Rothman, Peter J. Bickel, Elizaveta Levina, and J. Zhu. Sparse permutation invariant covariance estimation. *Electron. J. Stat.*, 2:494–515, 2008. ISSN 1935-7524. doi: 10.1214/08-EJS176. URL http://dx.doi.org/10.1214/08-EJS176.

Nathan Srebro and Adi Shraibman. Rank, trace-norm and max-norm. In *Learning theory*, volume 3559 of *Lecture Notes in Comput. Sci.*, pages 545–560. Springer, Berlin, 2005. doi: 10.1007/11503415_37. URL http://dx.doi.org/10.1007/11503415_37.

Tingni Sun and Cun-Hui Zhang. Sparse matrix inversion with scaled lasso. February 2012a.

Tingni Sun and Cun-Hui Zhang. Comment: "minimax estimation of large covariance matrices under $\ell_1$-norm". *Statist. Sinica*, 22:1354–1358, 2012b.

Jonathan E. Taylor, Richard Lockhart, Robert J. Tibshirani, and Robert J. Tibshirani. Post-selection adaptive inference for least angle regression and the lasso. *arXiv preprint arXiv:1401.3889*, January 2014.

Sara A. van de Geer, Peter Bühlmann, Ya'acov Ritov, and Ruben Dezeure. On asymptotically optimal confidence regions and tests for high-dimensional models. *Ann. Stat.*, 42(3):1166–1202, Jun 2014. ISSN 0090-5364. doi: 10.1214/14-aos1221. URL http://dx.doi.org/10.1214/14-AOS1221.

Roman Vershynin. Introduction to the non-asymptotic analysis of random matrices. In Y. C. Eldar and G. Kutyniok, editors, *Compressed Sensing: Theory and Applications*. Cambridge University Press, 2012.

Marten Wegkamp and Yue Zhao. Adaptive estimation of the copula correlation matrix for semiparametric elliptical copulas. *ArXiv e-prints, arXiv:1305.6526*, May 2013.

Lingzhou Xue and Hui Zou. Regularized rank-based estimation of high-dimensional nonparanormal graphical models. *Ann. Stat.*, 40(5):2541–2571, 2012. ISSN 0090-5364. doi: 10.1214/12-AOS1041. URL http://dx.doi.org/10.1214/12-AOS1041.

Lingzhou Xue, Hui Zou, and Tianxi Ca. Nonconcave penalized composite conditional likelihood estimation of sparse ising models. *Ann. Stat.*, 40(3):1403–1429, 2012. URL http://projecteuclid.org/euclid.aos/1344610588.

Eunho Yang, Genevera I. Allen, Zhandong Liu, and Pradeep Ravikumar. Graphical models via generalized linear models. In *Advances in Neural Information Processing Systems 25*, pages 1358–1366. Curran Associates, Inc., 2012. URL http://papers.nips.cc/paper/4617-graphical-models-via-generalized-linear-models.pdf.

Eunho Yang, Pradeep Ravikumar, Genevera I. Allen, and Zhandong Liu. On graphical models via univariate exponential family distributions. *ArXiv e-prints, arXiv:1301.4183*, January 2013.

Eunho Yang, Yulia Baker, Pradeep Ravikumar, Genevera I. Allen, and Zhandong Liu. Mixed graphical models via exponential families. In *Proc. 17th Int. Conf, Artif. Intel. Stat.*, pages 1042–1050, 2014.

M. Yuan. High dimensional inverse covariance matrix estimation via linear programming. *J. Mach. Learn. Res.*, 11:2261–2286, 2010.

M. Yuan and Y. Lin. Model selection and estimation in the gaussian graphical model. *Biometrika*, 94(1):19–35, 2007.

Cun-Hui Zhang and Stephanie S. Zhang. Confidence intervals for low dimensional parameters in high dimensional linear models. *J. R. Stat. Soc. B*, 76(1):217–242, Jul 2013. ISSN 1369-7412. doi: 10.1111/rssb.12026. URL http://dx.doi.org/10.1111/rssb.12026.

Tuo Zhao and Han Liu. Calibrated precision matrix estimation for high dimensional elliptical distributions. *IEEE Trans. Inf. Theory*, pages 1–1, 2014. ISSN 1557-9654. doi: 10.1109/tit.2014.2360980. URL http://dx.doi.org/10.1109/TIT.2014.2360980.

Tuo Zhao, Han Liu, Kathryn Roeder, John D. Lafferty, and Larry A. Wasserman. *huge: High-dimensional Undirected Graph Estimation*, 2014. URL http://CRAN.R-project.org/package=huge. R package version 1.2.6.

# A  Real data experiment

In this section, we evaluate the performance of the ROCKET method on a real data set, and compare with the Gaussian graphical model based approach of Ren et al. [2013] (using Pearson correlation) and the nonparanormal estimator proposed in Liu et al. [2009] (details for these methods are given in Section 5).

We use stock price closing data obtained via the R package `huge` [Zhao et al., 2014], which was gathered from publicly available data from Yahoo Finance.[6] The data consists of daily closing prices of 452 S&P 500 companies over 1258 days. We transform the data to consider the log-returns,

$$X_{ij} = \log\left(\frac{\text{Closing price of stock } j \text{ on day } i+1}{\text{Closing price of stock } j \text{ on day } i}\right).$$

While in practice there is dependence across time in this data set, we treat each row of $X$ as independent.

We perform two experiments on this data set. In Experiment 1, we test whether empirical results agree with the asymptotic normality predicted by the theory for the three methods—we do this by splittting the data into disjoint subsamples and comparing estimates across these subsamples. In Experiment 2, we use the full sample size and compare the estimates and confidence intervals produced by each of the three methods.

## A.1  Experiment 1: checking asymptotic normality

In this real data example, there is no available "ground truth" to compare to—that is, we do not know the true distribution of the data, and cannot compare our estimates to an exact true precision matrix $\Omega$. However, we can still check whether the estimators produced by these methods exhibit asymptotic normality (as claimed in the theory), by splitting the data into many subsamples and considering the empirical distribution of the estimators across these subsamples.

To construct our subsampled data, we randomly select $L = 25$ disjoint sets of size $n = 50$ from $\{1, \ldots, 1257\}$, denoted as $I_1, \ldots, I_L$. Due to this small sample size, we restrict our attention to companies in the categories `Materials` and `Consumer Staples`, which consist of 29 and 35 companies, respectively, for a total of $p_n = 64$ companies. For each $\ell = 1, \ldots, L$, define the $\ell$th data set $X^{(\ell)} = X_{I_\ell, S} \in \mathbb{R}^{n \times p}$, where $S \subset \{1, \ldots, 452\}$ identifies the $p = 64$ stocks of interest.

Next, for each pair $(a, b)$ of stocks, and for each subsample $\ell$, we compute $\breve{\Omega}_{ab}^{(\ell)}$ and $\breve{S}_{ab}^{(\ell)}$ using ROCKET. If the true distribution of the data follows the transelliptical model with precision matrix $\Omega$, then our main result, Theorem 3.5, implies that $\sqrt{n} \cdot (\breve{\Omega}_{ab}^{(\ell)} - \Omega_{ab})/\breve{S}_{ab}^{(\ell)}$ is approximately standard normal. Since $\breve{S}_{ab}^{(\ell)}$ concentrates near $S_{ab}$ (see Theorem 4.1), we should have

$$z_{ab}^{(\ell)} := \sqrt{n} \cdot \frac{\breve{\Omega}_{ab}^{(\ell)}}{\breve{S}_{ab}^{(\ell)}} \approx \sqrt{n} \cdot \frac{\Omega_{ab}}{S_{ab}} + N(0, 1).$$

In particular, this implies that the sample variance of the vector $(z_{ab}^{(1)}, \ldots, z_{ab}^{(L)})$ should have expectation approximately 1.

In Figure 4, we show a histogram of the sample variances $\mathsf{SampleVar}(z_{ab})$ across all $\binom{p_n}{2} = 2016$ pairs of variables. To compare to the Pearson and nonparanormal methods, we repeat this procedure for the estimators (and estimated variances) produced by the other two methods as well, which are also displayed in Figure 4. We see that ROCKET produces a mean sample variance $\approx 0.98$ (very near to 1), while the other two methods give mean sample variances of $\approx 1.28$ (Pearson) and $\approx 1.265$ (nonparanormal), substantially higher than the theoretical value of 1. This indicates that the normal approximation to the distribution of the estimator may be approximately valid for ROCKET, but does not have the correct scale (that is, the scale predicted by the theory) for the other two methods, on this data set.

---

[6]http://ichart.finance.yahoo.com

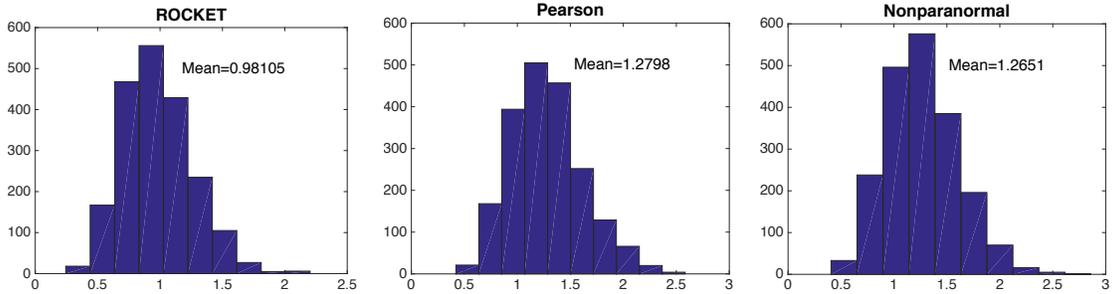

Figure 4: Sample variances of the rescaled estimator, $\sqrt{n} \cdot \breve{\Omega}_{ab}/\breve{S}_{ab}$, for each pair of variables $(a,b)$, using the subsampled stock data. The sample variances should be approximately 1 according to the theory (see Section A.1).

The vector $(z_{ab}^{(1)}, \ldots, z_{ab}^{(L)})$, in addition to having sample variance near 1, should also exhibit Gaussian-like tails according to the theory. To check this, we calculate the proportion of values in this vector lying near to the mean,

$$\frac{1}{L}\left|\left\{\ell : \left|z_{ab}^{(\ell)} - \overline{z}_{ab}\right| \leq 1.6449\sqrt{1 - \frac{1}{L}}\right\}\right|,$$

which should be approximately 90% according to the theory (using standard normal quantiles). The results are:

$$\begin{aligned} \text{ROCKET:} &\quad 90.55\% \text{ coverage} \\ \text{Pearson:} &\quad 85.01\% \text{ coverage} \\ \text{Nonparanormal:} &\quad 85.18\% \text{ coverage} \end{aligned}$$

We see that only the ROCKET method achieves the appropriate coverage.

### A.2 Experiment 2: estimating a graph

In the second experiment, we use the full sample size $n = 1257$ to estimate a sparse graph over the $p_n = 64$ stocks selected for Experiment 1, using each of the three methods. To do this, for each method we first produce a (approximate) p-value testing for the presence of an edge between each pair $(a,b)$ of variables. Recall that according to our main result, Theorem 3.5, if the pair of variables $(a,b)$ does not have an edge, then $\Omega_{ab} = 0$ and so $\sqrt{n} \cdot \breve{\Omega}_{ab}/\breve{S}_{ab}$ is approximately distributed as a standard normal variable. Then, using a two-sided z-test, we calculate a p-value $P_{ab} = 2 - 2\Phi\left(\left|\sqrt{n} \cdot \breve{\Omega}_{ab}/\breve{S}_{ab}\right|\right)$.

In Figure 5, we show the resulting graphs when edge $(a,b)$ is drawn whenever the p-value passes the threshold $P_{ab} < 0.00001$ or whenever $P_{ab} < 0.001$. The number of edges selected for each method is shown in the figures. Overall we see that ROCKET selects roughly the same number of edges as the Pearson method but less than the nonparanormal method, on this data set. Since the Pearson and nonparanormal methods do not exhibit approximately normal behavior across subsamples (Experiment 1), this should not be interpreted as a power comparison between the methods; the additional edges selected by the nonparanormal method, for instance, may not be as reliable since the p-value calculation is based on approximating the distribution of the estimator using a theoretical scaling that does not appear to hold for this method.

## B  Additional simulations

**Simulation 1 (chain graph).** Here we repeat Simulation 1 from the paper with a chain graph structure instead of a grid graph structure. The inverse covariance matrix $\Omega$ now has a chain structure with $\Omega^0_{j,j+1} = \Omega^0_{j+1,j} = 0.5$. We set $p_n = 1000$ and take sample size $n = 400$. Figures 6 and 7 show quantile-quantile (Q-Q) plots based on 1000 independent realizations of the test statistic error, $\sqrt{n} \cdot \frac{\breve{\Omega}_{ab} - \Omega_{ab}}{\breve{S}_{ab}}$, for the four methods together with the reference line showing quantiles of the standard normal distribution. The first row in the two figures illustrates actual performance of the methods, while the second row illustrates performance of an oracle procedure that does not need to solve a high-dimensional variable selection problem, but instead knows the sparsity pattern of $\Omega$. From these two figures, we observe that the quantiles of the test statistic error $\sqrt{n} \frac{\breve{\Omega}_{ab} - \Omega_{ab}}{\breve{S}_{ab}}$ based on ROCKET and the pseudo score estimators are closest to the quantiles of the standard normal random variable. We further quantify these results in Table 3, which reports empirical coverage of the confidence intervals based on $\sqrt{n} \frac{\breve{\Omega}_{ab} - \Omega_{ab}}{\breve{S}_{ab}}$. From the table, we can observe that the coverage of the confidence intervals based on ROCKET and the pseudo score are closest to nominal coverage of 95%.

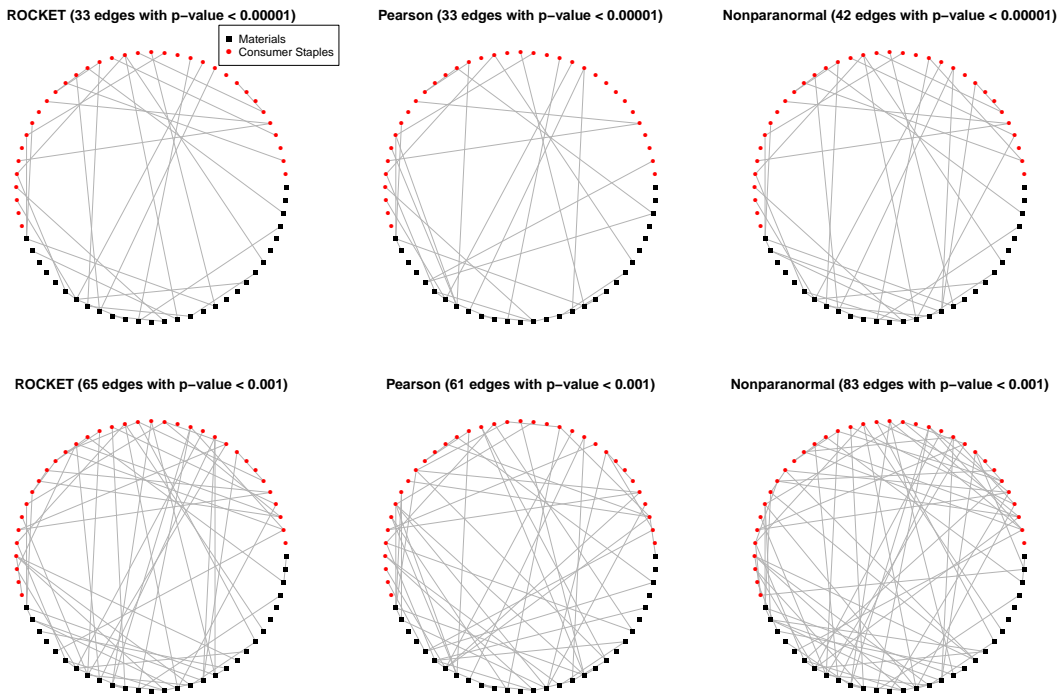

Figure 5: Estimated graph for the stock data, using the ROCKET, Pearson, and nonparanormal methods (see Section A.2). An edge is displayed for each pair of variables $(a, b)$ with p-value $P_{ab} <$ 0.00001 (top row) and $P_{ab} < 0.001$ (bottom row). Graphs were drawn using the `igraph` package [Csardi and Nepusz, 2006] in R [R Core Team, 2012].

|  | ROCKET | Pearson | NPN | Pseudo Score |
| --- | --- | --- | --- | --- |
| $\omega_{10,11} = 10.38$ | 94.5 (10.16) | 56.1 (4.90) | 31.8 (3.53) | 93.5 (11.12) |
| $\omega_{10,12} = 0$ | 96.3 (9.78) | 62.2 (4.62) | 74.0 (3.32) | 93.8 (10.11) |
| $\omega_{10,20} = 0$ | 95.9 (16.04) | 62.3 (6.04) | 81.9 (4.25) | 95.8 (17.69) |

Table 3: Simulation 1 (transelliptical data). Percent empirical coverage (average length) of 95% confidence intervals based on 1000 independent simulation runs. $\Omega$ corresponds to a chain graph structure.

**Simulation 4.** In this additional simulation, we evaluate robustness of the procedures to three data contamination mechanisms: random row contamination, deterministic row contamination, and element contamination mechanism. Let $r \in (0, 1)$ be the contamination level. As before we set $n = 400$ and $p_n = 900$. The precision matrix $\Omega$ is chosen as the grid structure used in earlier simulations. For both row contamination and element contamination, we consider two settings, where the (uncorrupted) samples are drawn either from $N(0, \Omega^{-1})$, or from $\mathsf{E}(0, \Omega^{-1}, \xi)$ with $\xi \sim t_5$.

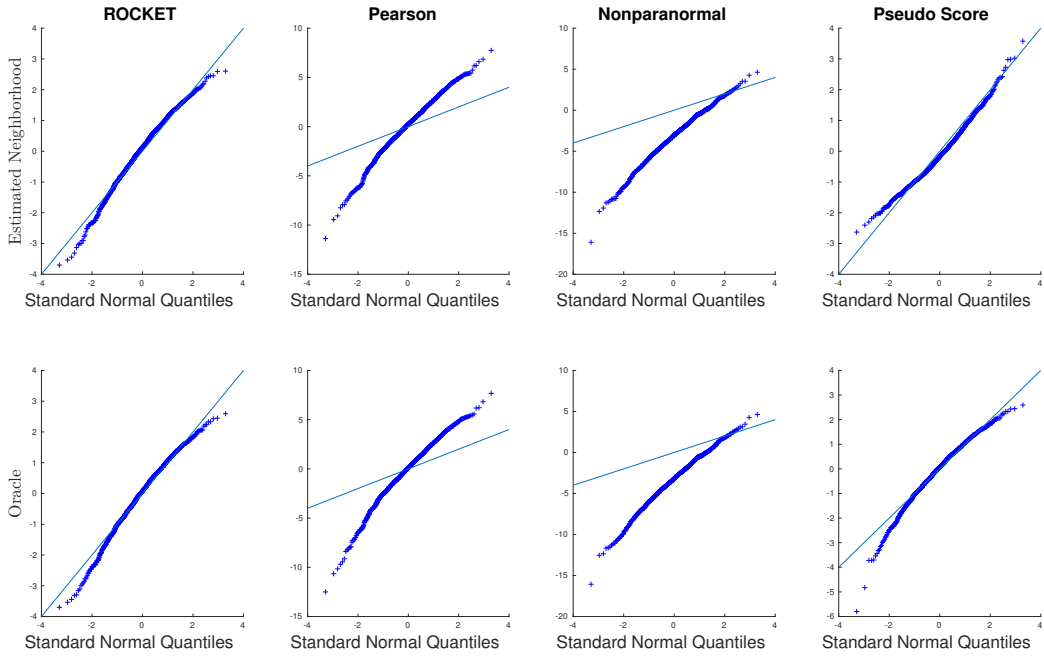

Figure 6: Simulation 1 (transelliptical data). Q-Q plot of $\sqrt{n} \cdot \frac{\breve{\Omega}_{ab} - \Omega_{ab}}{\breve{S}_{ab}}$ with $a = 10$ and $b = 11$ (true edge) when $\Omega$ corresponds to a chain graph structure.

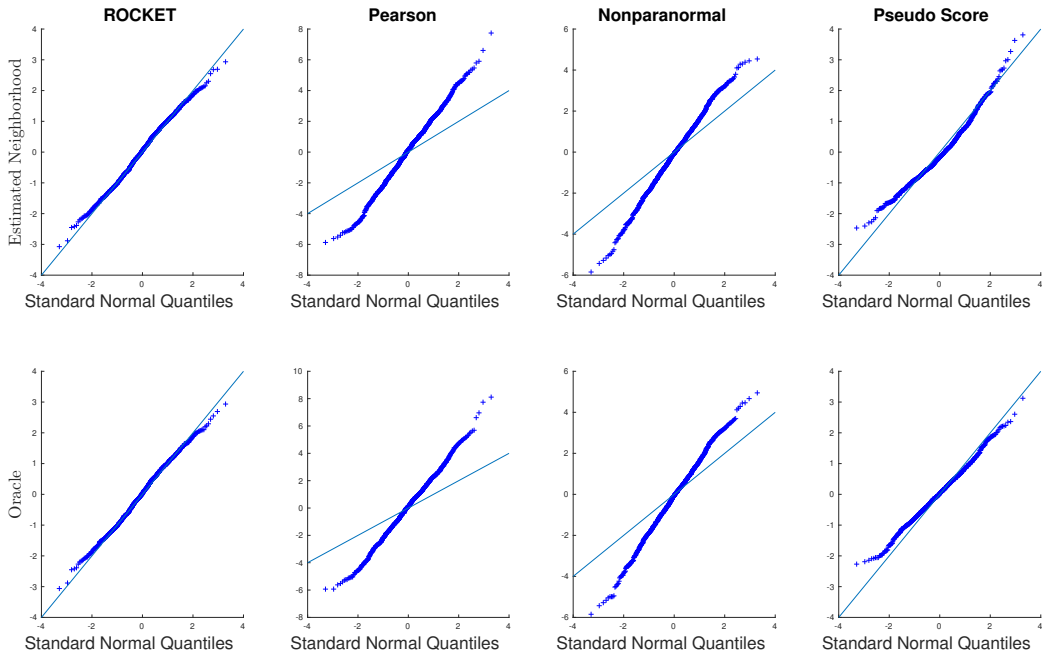

Figure 7: Simulation 1 (transelliptical data). Q-Q plot of $\sqrt{n} \cdot \frac{\breve{\Omega}_{ab} - \Omega_{ab}}{\breve{S}_{ab}}$ with $a = 10$ and $b = 12$ (non-edge close to a true edge) when $\Omega$ corresponds to a chain graph structure.

|  |  | $r$ | ROCKET | Pearson | NPN | Pseudo Score |
|---|---|---|---|---|---|---|
| Gaussian | $\omega_{(2,2),(2,3)} = 0.37$ | 0.01 | 94.9 (0.34) | 93.6 (0.32) | 93.4 (0.32) | 96.0 (0.33) |
|  |  | 0.05 | 94.1 (0.35) | 93.8 (0.32) | 94.4 (0.32) | 95.4 (0.32) |
|  |  | 0.1 | 95.1 (0.37) | 92.5 (0.32) | 92.2 (0.32) | 96.7 (0.37) |
|  | $\omega_{(2,2),(3,3)} = 0$ | 0.01 | 95.3 (0.35) | 95.3 (0.32) | 95.4 (0.31) | 95.1 (0.33) |
|  |  | 0.05 | 93.7 (0.36) | 92.6 (0.32) | 92.7 (0.31) | 95.0 (0.33) |
|  |  | 0.1 | 96.2 (0.38) | 92.9 (0.32) | 93.4 (0.31) | 94.5 (0.34) |
|  | $\omega_{(2,2),(10,10)} = 0$ | 0.01 | 95.8 (0.36) | 95.9 (0.33) | 96.1 (0.32) | 95.5 (0.36) |
|  |  | 0.05 | 94.7 (0.38) | 93.8 (0.33) | 94.0 (0.32) | 95.3 (0.38) |
|  |  | 0.1 | 94.7 (0.40) | 93.5 (0.33) | 93.0 (0.32) | 94.6 (0.43) |
| $t_5$ | $\omega_{(2,2),(2,3)} = 0.37$ | 0.01 | 94.4 (0.38) | 85.1 (0.35) | 89.0 (0.32) | 96.1 (0.39) |
|  |  | 0.05 | 93.7 (0.39) | 84.6 (0.36) | 89.7 (0.32) | 96.4 (0.39) |
|  |  | 0.1 | 94.9 (0.40) | 83.5 (0.36) | 86.4 (0.31) | 96.0 (0.44) |
|  | $\omega_{(2,2),(3,3)} = 0$ | 0.01 | 95.5 (0.39) | 88.8 (0.35) | 92.9 (0.31) | 94.6 (0.34) |
|  |  | 0.05 | 94.1 (0.41) | 87.4 (0.35) | 92.3 (0.31) | 95.4 (0.37) |
|  |  | 0.1 | 95.0 (0.42) | 85.8 (0.36) | 88.8 (0.31) | 96.1 (0.38) |
|  | $\omega_{(2,2),(10,10)} = 0$ | 0.01 | 95.2 (0.41) | 88.7 (0.36) | 90.5 (0.32) | 94.2 (0.41) |
|  |  | 0.05 | 95.4 (0.42) | 87.2 (0.36) | 91.7 (0.32) | 94.7 (0.42) |
|  |  | 0.1 | 94.6 (0.44) | 84.5 (0.37) | 89.6 (0.32) | 94.3 (0.45) |

Table 4: Simulation 4 (Gaussian and transelliptical data; row corruption mechanism). Percent empirical coverage (average length) of 95% confidence intervals based on 1000 independent simulation runs. $\Omega$ corresponds to a grid graph structure. $r$ denotes the fraction of corrupted rows.

The row contamination mechanism corrupts $\lfloor nr \rfloor$ rows of the data matrix. In the random row contamination mechanism, each corrupted row is filled with i.i.d. entries drawn from a $t_{1.5}$ distribution. These contaminated rows are heavy tailed and do not come from the same elliptical copula family as the uncorrupted data. In the deterministic row contamination mechanism, each corrupted row is equal to the vector $(+5, -5, +5, -5, \ldots) \in \mathbb{R}^p$, where the numbers $+5$ and $-5$ occur in an alternating way. The element contamination mechanism selects $\lfloor npr \rfloor$ elements and substitutes them with with a value drawn from either $N(3, 3)$ or $N(-3, 3)$ with equal probability.

The results of the simulation are summarized in Tables 4, 5, and 6. The deterministic row corruption mechanism is the most malicious and hurts the non-robust procedures the most. From the simulation, we observe that the procedures using the Kendall's tau correlation matrix tend to be more robust. Furthermore, coverage of non-zero elements tends to be more severely affected by corruption than the coverage of zero-elements.

|  |  | $r$ | ROCKET | Pearson | NPN | Pseudo Score |
|---|---|---|---|---|---|---|
| Gaussian | $\omega_{(2,2),(2,3)} = 0.37$ | 0.01 | 94.4 (0.35) | 90.2 (0.40) | 86.7 (0.32) | 93.6 (0.39) |
|  |  | 0.05 | 80.6 (0.51) | 40.1 (0.85) | 85.4 (0.45) | 75.7 (1.00) |
|  |  | 0.1 | 53.8 (0.86) | 23.6 (1.40) | 67.1 (0.63) | 51.6 (0.46) |
|  | $\omega_{(2,2),(3,3)} = 0$ | 0.01 | 94.6 (0.37) | 93.7 (0.39) | 84.7 (0.32) | 93.5 (0.29) |
|  |  | 0.05 | 93.6 (0.52) | 92.8 (0.83) | 93.6 (0.40) | 95.2 (0.47) |
|  |  | 0.1 | 93.1 (0.83) | 88.2 (1.30) | 92.6 (0.49) | 96.0 (1.40) |
|  | $\omega_{(2,2),(10,10)} = 0$ | 0.01 | 94.4 (0.37) | 94.8 (0.37) | 88.2 (0.33) | 93.1 (0.42) |
|  |  | 0.05 | 93.4 (0.49) | 83.4 (0.75) | 97.8 (0.41) | 92.9 (0.59) |
|  |  | 0.1 | 92.0 (0.64) | 79.4 (1.20) | 87.3 (0.56) | 91.3 (0.49) |
| $t_5$ | $\omega_{(2,2),(2,3)} = 0.37$ | 0.01 | 92.4 (0.39) | 78.2 (0.43) | 82.5 (0.32) | 91.9 (0.47) |
|  |  | 0.05 | 77.0 (0.61) | 52.8 (0.74) | 83.2 (0.47) | 72.1 (1.20) |
|  |  | 0.1 | 60.4 (1.00) | 36.5 (1.10) | 73.0 (0.65) | 62.8 (0.75) |
|  | $\omega_{(2,2),(3,3)} = 0$ | 0.01 | 94.8 (0.41) | 88.4 (0.43) | 84.2 (0.32) | 92.3 (0.35) |
|  |  | 0.05 | 94.1 (0.63) | 88.1 (0.73) | 91.7 (0.42) | 95.6 (0.63) |
|  |  | 0.1 | 94.0 (1.00) | 87.6 (1.10) | 90.1 (0.51) | 96.5 (2.30) |
|  | $\omega_{(2,2),(10,10)} = 0$ | 0.01 | 94.5 (0.42) | 90.2 (0.42) | 86.3 (0.33) | 91.9 (0.51) |
|  |  | 0.05 | 93.5 (0.59) | 75.0 (0.66) | 90.7 (0.42) | 94.6 (0.62) |
|  |  | 0.1 | 92.8 (0.76) | 68.7 (0.97) | 77.7 (0.58) | 93.1 (0.57) |

Table 5: Simulation 4 (Gaussian and transelliptical data; deterministic row corruption mechanism). Percent empirical coverage (average length) of 95% confidence intervals based on 1000 independent simulation runs. $\Omega$ corresponds to a grid graph structure. $r$ denotes the fraction of corrupted rows.

|  | $r$ | ROCKET | Pearson | NPN | Pseudo Score |
|---|---|---|---|---|---|
| Gaussian | | | | | |
| $\omega_{(2,2),(2,3)} = 0.37$ | 0.01 | 94.2 (0.32) | 72.6 (0.28) | 90.4 (0.30) | 96.5 (0.35) |
| | 0.05 | 93.6 (0.29) | 8.80 (0.23) | 52.0 (0.27) | 93.5 (0.40) |
| | 0.10 | 85.9 (0.27) | 0.70 (0.21) | 15.3 (0.24) | 82.6 (0.52) |
| $\omega_{(2,2),(3,3)} = 0$ | 0.01 | 95.6 (0.34) | 92.3 (0.28) | 94.6 (0.30) | 94.4 (0.32) |
| | 0.05 | 91.4 (0.30) | 76.8 (0.23) | 88.0 (0.26) | 92.7 (0.33) |
| | 0.1 | 89.7 (0.28) | 73.8 (0.20) | 81.7 (0.24) | 92.7 (0.34) |
| $\omega_{(2,2),(10,10)} = 0$ | 0.01 | 95.4 (0.35) | 95.4 (0.29) | 95.8 (0.31) | 95.6 (0.35) |
| | 0.05 | 95.3 (0.31) | 94.2 (0.23) | 95.3 (0.27) | 94.3 (0.32) |
| | 0.1 | 95.5 (0.28) | 94.9 (0.20) | 94.7 (0.24) | 95.7 (0.30) |
| $t_5$ | | | | | |
| $\omega_{(2,2),(2,3)} = 0.37$ | 0.01 | 94.4 (0.36) | 80.9 (0.32) | 85.2 (0.30) | 95.8 (0.39) |
| | 0.05 | 93.1 (0.32) | 21.3 (0.25) | 51.4 (0.27) | 92.5 (0.45) |
| | 0.1 | 94.8 (0.29) | 4.0 (0.22) | 20.9 (0.24) | 83.8 (0.54) |
| $\omega_{(2,2),(3,3)} = 0$ | 0.01 | 94.2 (0.38) | 87.3 (0.31) | 89.9 (0.30) | 95.0 (0.34) |
| | 0.05 | 93.1 (0.33) | 79.6 (0.25) | 86.3 (0.27) | 94.2 (0.36) |
| | 0.1 | 94.4 (0.29) | 75.7 (0.22) | 85.0 (0.24) | 93.7 (0.36) |
| $\omega_{(2,2),(10,10)} = 0$ | 0.01 | 94.6 (0.39) | 90.6 (0.32) | 92.2 (0.31) | 95.6 (0.39) |
| | 0.05 | 94.6 (0.34) | 92.9 (0.25) | 92.5 (0.27) | 94.9 (0.35) |
| | 0.1 | 96.3 (0.30) | 94.8 (0.22) | 94.7 (0.24) | 94.6 (0.30) |

Table 6: Simulation 4 (Gaussian and transelliptical data; element corruption mechanism). Percent empirical coverage (average length) of 95% confidence intervals based on 1000 independent simulation runs. $\Omega$ corresponds to a grid graph structure. $r$ denotes the fraction of corrupted elements.

# C Gaussian vectors and the sign-subgaussian property

In this section we prove Lemma 4.5, which shows that that a Gaussian vector $Z \sim N(\mu, \Sigma)$ satisfies the sign-subgaussianity property, that is, the centered sign vector $\text{sign}(Z) - \mathbb{E}[\text{sign}(Z)]$ is itself subgaussian. (In this paper, we will apply this lemma only with $\mu = 0$, in which case $\mathbb{E}[\text{sign}(Z)] = 0$ and so $\text{sign}(Z)$ is therefore subgaussian.)

**Lemma 4.5.** *Let $Z \sim N(\mu, \Sigma)$. Then $\text{sign}(Z) - \mathbb{E}[\text{sign}(Z)]$ is $\mathsf{C}(\Sigma)$-subgaussian.*

*Proof of Lemma 4.5.* Without loss of generality, rescale so that $\lambda_{\min}(\Sigma) = 1$ and then $\mathsf{C}(\Sigma) = \lambda_{\max}(\Sigma)$ (note that the mean $\mu$ is therefore rescaled as well). Write $\Sigma = AA^\top + \mathbf{I}_p$ for some matrix $A \in \mathbb{R}^{n \times n}$. Then we can write $Z = \mu + X + AY$, where $X, Y \overset{iid}{\sim} N(0, \mathbf{I}_p)$. Then, for any fixed vector $v \in \mathbb{R}^p$,

$$\mathbb{E}\left[e^{v^\top \text{sign}(Z)}\right] = \mathbb{E}\left[\mathbb{E}\left[e^{v^\top \text{sign}(Z)} \mid Y\right]\right] = \mathbb{E}\left[\mathbb{E}\left[e^{v^\top \text{sign}(\mu + X + AY)} \mid Y\right]\right]$$
$$= \mathbb{E}\left[\prod_i \mathbb{E}\left[e^{v_i \text{sign}(\mu_i + X_i + (AY)_i)} \mid Y\right]\right],$$

where the last step holds because, conditional on $Y$, each of the terms $\text{sign}(\mu_i + X_i + (AY)_i)$ depends on $X_i$ only, and therefore these terms are conditionally independent. Next, observe that

$$\mathbb{E}[\text{sign}(\mu_i + X_i + (AY)_i) \mid Y] = \mathbb{E}[\text{sign}(N(0,1) + \mu_i + (AY)_i) \mid Y]$$
$$= \Phi(\mu_i + (AY)_i) - \Phi(-\mu_i - (AY)_i) = \psi(\mu_i + (AY)_i),$$

where we define $\psi(z) = \Phi(z) - \Phi(-z)$ for $z \in \mathbb{R}$. Then, for each $i$,

$$\mathbb{E}\left[e^{v_i \text{sign}(\mu_i + X_i + (AY)_i)} \mid Y\right]$$
$$= \mathbb{E}\left[e^{v_i(\text{sign}(\mu_i + X_i + (AY)_i) - \psi(\mu_i + (AY)_i))} \mid Y\right] \cdot e^{v_i \psi(\mu_i + (AY)_i)}$$
$$\leq e^{v_i^2/2} \cdot e^{v_i \psi(\mu_i + (AY)_i)},$$

where the inequality is proved by applying Hoeffding's Lemma (see, for example, Massart [2007, Lemma 2.6]) to the bounded mean-zero random variable

$$v_i \left(\text{sign}(\mu_i + X_i + (AY)_i) - \psi(\mu_i + (AY)_i)\right)$$

(where, since we are conditioning on $Y$, only $X$ is treated as random). Combining the calculations so far, we have

$$\mathbb{E}\left[e^{v^\top \text{sign}(Z)}\right] \leq e^{\|v\|_2^2/2} \cdot \mathbb{E}\left[e^{v^\top \psi(\mu + AY)}\right],$$

where $\psi(\mu + AY)$ applies the function $\psi(\cdot)$ elementwise to the vector $\mu + AY$.

Next we show that $y \mapsto v^\top \psi(\mu + Ay)$ is Lipschitz over $y \in \mathbb{R}^n$. Note that $z \mapsto \psi(z)$ is 1-Lipschitz over $z \in \mathbb{R}$ since the density of the standard normal distribution is bounded uniformly as $\phi(z) \leq \frac{1}{\sqrt{2\pi}} \leq \frac{1}{2}$. For any $y, y' \in \mathbb{R}^n$, we have

$$\left|v^\top \psi(\mu + Ay) - v^\top \psi(\mu + Ay')\right| \leq \sum_i |v_i| \cdot \left|\psi((\mu + Ay)_i) - \psi((\mu + Ay')_i)\right|$$
$$\leq \sum_i |v_i| \cdot |(\mu + Ay)_i - (\mu + Ay')_i| \leq \|v\|_2 \cdot \|A(y - y')\|_2$$
$$\leq \|v\|_2 \cdot \sqrt{\lambda_{\max}(\Sigma) - 1} \cdot \|y - y'\|_2,$$

where the last step is true because

$$||A||_{\mathsf{op}} = \sqrt{||\Sigma - \mathbf{I}_p||_{\mathsf{op}}} = \sqrt{\mathsf{C}(\Sigma) - 1}\ .$$

Therefore, $y \mapsto v^\top \psi(\mu + Ay)$ is $\left(||v||_2 \cdot \sqrt{\mathsf{C}(\Sigma) - 1}\right)$-Lipschitz in $y$. We can now apply standard concentration results for Lipschitz functions of a Gaussian: by Massart [2007, Proposition 3.5],

$$\mathbb{E}\left[e^{v^\top(\psi(\mu + AY) - \mathbb{E}[\psi(\mu + AY)])}\right] \leqslant e^{||v||_2^2 (\mathsf{C}(\Sigma) - 1)/2}.$$

Therefore,
$$\mathbb{E}\left[e^{v^\top \operatorname{sign}(Z)}\right] \leqslant \mathbb{E}\left[e^{||v||_2^2/2 + v^\top \psi(\mu + AY)}\right] \leqslant e^{||v||_2^2 \cdot \mathsf{C}(\Sigma)/2} \cdot e^{v^\top \mathbb{E}[\psi(\mu + AY)]}\ .$$

Finally, for each $i$,

$$\mathbb{E}\left[\psi(\mu_i + (AY)_i)\right] = \mathbb{E}\left[\mathbb{E}\left[\operatorname{sign}(\mu_i + X_i + (AY)_i) \mid Y\right]\right]$$
$$= \mathbb{E}\left[\operatorname{sign}(\mu_i + X_i + (AY)_i\right] = \mathbb{E}\left[\operatorname{sign}(Z_i)\right],$$

or in other words, $\mathbb{E}\left[\psi(\mu + AY)\right] = \mathbb{E}\left[\operatorname{sign}(Z)\right]$. Combining everything, we have proved that $\mathbb{E}\left[e^{v^\top \operatorname{sign}(Z)}\right] \leqslant e^{||v||_2^2 \cdot \mathsf{C}(\Sigma)/2} \cdot e^{v^\top \mathbb{E}[\operatorname{sign}(Z)]}$, as desired. □

# D  Proof of main result

## D.1  Preliminaries

We first compute bounds on $||\gamma_c||_2$ and $||\gamma_c||_1$ for each $c = a, b$, which we will use many times in the proofs below. First, for $c = a, b$ note that

$$||\gamma_c||_2 = ||\Sigma_I^{-1} \Sigma_{Ic}||_2 \leqslant ||\Sigma_I^{-1}|| \cdot ||\Sigma_{Ic}||_2 \leqslant [\lambda_{\min}(\Sigma)]^{-1} \cdot \lambda_{\max}(\Sigma) \leqslant C_{\mathsf{cov}} \tag{D.1}$$

by Assumption 3.1. Next,

$$\begin{aligned}
||\gamma_c||_1 &= ||\Sigma_I^{-1} \Sigma_{Ic}||_1 \\
&= || - \Omega_{I,ab} \Theta_{ab,c}||_1 \quad \text{(by matrix blockwise inversion)} \\
&= \sum_{j \in I} |\Omega_{j,ab} \Theta_{ab,c}| \leqslant \sum_{j \in I} ||\Omega_{j,ab}||_1 ||\Theta_{ab,c}||_\infty \\
&\leqslant C_{\mathsf{cov}} \sum_{j \in I} ||\Omega_{j,ab}||_1
\end{aligned}$$

since
$$||\Theta||_\infty \leqslant \lambda_{\max}(\Theta) = (\lambda_{\min}(\Omega_{ab,ab}))^{-1} \leqslant (\lambda_{\min}(\Omega))^{-1} = \lambda_{\max}(\Sigma) \leqslant \mathsf{C}_{\mathsf{cov}}.$$

Therefore,
$$||\gamma_c||_1 \leqslant C_{\mathsf{cov}} (||\Omega_a||_1 + ||\Omega_b||_1) \leqslant 2 C_{\mathsf{cov}} C_{\mathsf{sparse}} \sqrt{k_n} \tag{D.2}$$

by applying Assumption 3.2.

## D.2 Proof of Theorem 4.1: asymptotic normality of the oracle estimator

**Theorem 4.1.** *Suppose that Assumptions 3.1, 3.2, and 3.4 hold. Then there exist constants $C_{\mathsf{normal}}, C_{\mathsf{variance}}$ depending on $C_{\mathsf{cov}}, C_{\mathsf{sparse}}, C_{\mathsf{kernel}}$ but not on $(n, p_n, k_n)$, such that*

$$\sup_{t \in \mathbb{R}} \left| \mathbb{P}\left\{ \sqrt{n} \cdot \frac{\widetilde{\Theta}_{ab} - \Theta_{ab}}{S_{ab} \cdot \det(\Theta)} \leq t \right\} - \Phi(t) \right| \leq C_{\mathsf{normal}} \cdot \frac{k_n \log(p_n)}{\sqrt{n}} + \frac{1}{2p_n} \,,$$

*where $S_{ab}$ is defined in the proof and satisfies $S_{ab} \cdot \det(\Theta) \geq C_{\mathsf{variance}} > 0$.*

*Proof of Theorem 4.1.* We first show that the error $\widetilde{\Theta}_{ab} - \Theta_{ab}$ can be approximated by a linear function of the Kendall's tau estimator $\widehat{T}$. Define vectors $u, v \in \mathbb{R}^{p_n}$ with entries

$$u_a = 1, u_b = 0, u_I = -\gamma_a \text{ and } v_a = 0, v_b = 1, v_I = -\gamma_b \,.$$

Then by definition, we have $\widetilde{\Theta}_{ab} = u^\top \widehat{\Sigma} v$ and $\Theta_{ab} = u^\top \Sigma v$, that is, the error is given by

$$\widetilde{\Theta}_{ab} - \Theta_{ab} = u^\top (\widehat{\Sigma} - \Sigma) v \,.$$

Next, since $\widehat{\Sigma} = \sin\left(\frac{\pi}{2}\widehat{T}\right)$ and $\Sigma = \sin\left(\frac{\pi}{2}T\right)$, we take a second-order Taylor expansion of $\sin(\cdot)$ to see that, for some $t \in [0,1]$,

$$\widetilde{\Theta}_{ab} - \Theta_{ab} = u^\top \Bigg[ \frac{\pi}{2} \cos\left(\frac{\pi}{2}T\right) \circ (\widehat{T} - T) -$$
$$\frac{1}{2} \cdot \left(\frac{\pi}{2}\right)^2 \cdot \sin\left(\frac{\pi}{2}(t \cdot T + (1-t) \cdot \widehat{T})\right) \circ (\widehat{T} - T) \circ (\widehat{T} - T) \Bigg] v \,. \quad (\text{D.3})$$

Next, we rewrite this linear term. We have

$$L := u^\top \left[ \cos\left(\frac{\pi}{2}T\right) \circ \widehat{T} \right] v = \frac{1}{\binom{n}{2}} \sum_{i < i'} \operatorname{sign}(X_i - X_{i'})^\top \left( uv^\top \circ \cos\left(\frac{\pi}{2}T\right) \right) \operatorname{sign}(X_i - X_{i'}) \,,$$

which is a U-statistic of order 2 with respect to the data $(X_1, \ldots, X_n)$. Note that

$$L - \mathbb{E}[L] = u^\top \left[\frac{\pi}{2}\cos\left(\frac{\pi}{2}T\right)\circ \widehat{T}\right]v - u^\top \left[\frac{\pi}{2}\cos\left(\frac{\pi}{2}T\right)\circ \mathbb{E}[\widehat{T}]\right]v = u^\top \left[\frac{\pi}{2}\cos\left(\frac{\pi}{2}T\right)\circ (\widehat{T} - T)\right]v \,.$$

Define the kernel $g(X, X') = \operatorname{sign}(X - X')^\top \left( uv^\top \circ \cos\left(\frac{\pi}{2}T\right) \right) \operatorname{sign}(X - X')$, and let $g_1(X) = \mathbb{E}[g(X, X') \mid X]$, where $X, X' \stackrel{iid}{\sim} \mathsf{TE}(\Sigma, \xi; f_1, \ldots, f_p)$. Let $\nu_{g_1}^2 = \operatorname{Var}(g_1(X))$ and $\eta_g^3 = \mathbb{E}\left[|g(X, X')|^3\right]$. By Callaert and Janssen [1978, Section 2], we have

$$\sup_{t \in \mathbb{R}} \left| \mathbb{P}\left\{ \frac{\sqrt{n}(L - \mathbb{E}[L])}{2\nu_{g_1}} \leq t \right\} - \Phi(t) \right| \leq C_{\mathsf{Ustat}} \cdot \frac{\eta_g^3}{\nu_{g_1}^3} \cdot \frac{1}{\sqrt{n}} \,, \quad (\text{D.4})$$

for a universal constant $C_{\mathsf{Ustat}}$. Next we bound the ratio $\frac{\eta_g^3}{\nu_{g_1}^3}$ in the following lemma, which is proved in Appendix F.8.

**Lemma D.1.** *Suppose that Assumptions 3.1, 3.2 and 3.4 hold. Let $g(X, X')$ and $g_1(X)$ be defined as in the proof of Theorem 4.1. Then*

$$\nu_{g_1}^2 = \operatorname{Var}(g_1(X)) \geq \frac{1}{\pi^2} C_{\mathsf{variance}}^2$$

*and*

$$\nu_{g_1}^3 \leq \eta_g^3 = \mathbb{E}\left[|g(X, X')|^3\right] \leq C_{\mathsf{moment}}$$

*where $C_{\mathsf{variance}}, C_{\mathsf{moment}}$ are constants depending only on $C_{\mathsf{cov}}, C_{\mathsf{kernel}}$ and not on $(n, p_n, k_n)$.*

In particular, this lemma implies that $S_{ab} := \pi\nu_{g_1} (\det(\Theta))^{-1} \geq C_{\mathsf{variance}} \cdot (\det(\Theta))^{-1}$.

Next, the linear term $L$ analysed here provides only an approximation to $\widetilde{\Theta}_{ab} - \Theta_{ab}$. Define

$$\Delta = \widetilde{\Theta}_{ab} - \Theta_{ab} - \frac{\pi}{2}(L - \mathbb{E}[L]) \ .$$

Then we have

$$\begin{aligned}
|\Delta| &= \left| u^\top \left[ \frac{1}{2} \cdot \left(\frac{\pi}{2}\right)^2 \cdot \sin\left(\frac{\pi}{2}(t \cdot T + (1-t) \cdot \widehat{T})\right) \circ (\widehat{T} - T) \circ (\widehat{T} - T) \right] v \right| \\
&\leq ||u||_1 ||v||_1 \left\| \frac{1}{2} \cdot \left(\frac{\pi}{2}\right)^2 \cdot \sin\left(\frac{\pi}{2}(t \cdot T + (1-t) \cdot \widehat{T})\right) \circ (\widehat{T} - T) \circ (\widehat{T} - T) \right\|_\infty \\
&\leq \frac{\pi^2}{8} ||u||_1 ||v||_1 ||\widehat{T} - T||_\infty^2 \\
&\leq \frac{\pi^2}{8} \cdot k_n \cdot (1 + 2C_{\mathsf{cov}} C_{\mathsf{sparse}})^2 \cdot ||\widehat{T} - T||_\infty^2 \ ,
\end{aligned} \tag{D.5}$$

where the last inequality holds by (D.2).

Finally, the next lemma is proved in de la Pena and Giné [1999].

**Lemma D.2** ([de la Pena and Giné, 1999, Theorem 4.1.8]). *For any $\delta > 0$, with probability at least $1 - \delta$,*

$$||\widehat{T} - T||_\infty \leq \sqrt{\frac{4 \log\left(2\binom{p_n}{2}/\delta\right)}{n}} \ .$$

Applying this lemma with $\delta = \frac{1}{2p_n}$, we have $||\widehat{T} - T||_\infty^2 \leq \frac{4\log(2p_n^3)}{n} \leq \frac{16 \log(p_n)}{n}$ with probability at least $1 - \frac{1}{2p_n}$.

To summarize the computations so far, we have $\widetilde{\Theta}_{ab} - \Theta_{ab} = \frac{\pi}{2}(L - \mathbb{E}[L]) + \Delta$, where (D.4) gives an asymptotic normality result for the linear term $(L - \mathbb{E}[L])$, while (D.5) gives a bound on $\Delta$. To prove therefore that $\widetilde{\Theta}_{ab} - \Theta_{ab}$ is asymptotically normal, we will use the following lemma (proved in Appendix F.1):

**Lemma D.3.** *Let $A, B, C$ be random variables such that*

$$\sup_{t \in \mathbb{R}} |\mathbb{P}\{A \leq t\} - \Phi(t)| \leq \epsilon_A \quad \text{and} \quad \mathbb{P}\{|B| \leq \delta_B, |C| \leq \delta_C\} \geq 1 - \epsilon_{BC} \ ,$$

*where $\epsilon_A, \epsilon_{BC}, \delta_B, \delta_C \in (0, 1)$. Then the variable $(A + B) \cdot (1 + C)^{-1}$ converges to a standard normal distribution with rate*

$$\sup_{t \in \mathbb{R}} \left| \mathbb{P}\left\{(A + B) \cdot (1 + C)^{-1} \leq t\right\} - \Phi(t) \right| \leq \delta_B + \frac{\delta_C}{1 - \delta_C} + \epsilon_A + \epsilon_{BC} \ .$$

We apply this lemma with $A = \frac{\pi}{2} \cdot \sqrt{n} \cdot \frac{L - \mathbb{E}[L]}{S_{ab} \cdot \det(\Theta)}$ and $B = \sqrt{n} \cdot \frac{\Delta}{S_{ab} \cdot \det(\Theta)}$ and $C = 0$. We have

$$\sup_{t \in \mathbb{R}} |\mathbb{P}\{A \leq t\} - \Phi(t)| \leq C_{\mathsf{Ustat}} \cdot \frac{C_{\mathsf{moment}}}{\left(\frac{1}{\pi^2} C_{\mathsf{variance}}^2\right)^{1.5}} \cdot \frac{1}{\sqrt{n}} =: \epsilon_A$$

by (D.4) and Lemma D.1. Furthermore,

$$\mathbb{P}\left\{ |B| \leq \sqrt{n} \cdot \frac{\frac{\pi^2}{8} \cdot k_n \cdot (1 + 2C_{\mathsf{cov}} C_{\mathsf{sparse}})^2 \cdot \frac{16 \log(p_n)}{n}}{C_{\mathsf{variance}}} =: \delta_B \right\}$$

$$\leq \mathbb{P}\left\{ ||\widehat{T} - T||_\infty^2 \leq \frac{16 \log(p_n)}{n} \right\} \geq 1 - \frac{1}{2p_n} =: 1 - \epsilon_{BC}$$

by (D.5) and Lemmas D.1 and D.2. Noting that $\sqrt{n} \cdot \frac{\widetilde{\Theta}_{ab} - \Theta_{ab}}{S_{ab}} = A + B$, and defining

$$C_{\text{normal}} = \frac{2\pi^2(1 + 2C_{\text{cov}}C_{\text{sparse}})^2}{C_{\text{variance}}} + C_{\text{Ustat}} \cdot \frac{C_{\text{moment}}}{\left(\frac{1}{\pi^2}C_{\text{variance}}^2\right)^{1.5}},$$

we have proved the desired result. □

## D.3 Proof of Theorem 4.2: gap between the estimator and the oracle estimator, and estimation of the variance

**Theorem 4.2.** *Suppose that Assumptions 3.1, 3.2, and 3.3 hold. Then there exists a constant $C_{\text{oracle}}$, depending on $C_{\text{cov}}, C_{\text{sparse}}, C_{\text{est}}$ but not on $(n, p_n, k_n)$, such that, if[7] $n \geq 15 k_n \log(p_n)$, then, with probability at least $1 - \frac{1}{2p_n} - \delta_n$,*

$$||\breve{\Theta} - \widetilde{\Theta}||_\infty \leq C_{\text{oracle}} \cdot \frac{k_n \log(p_n)}{n}$$

*and*

$$\left|\breve{S}_{ab} \cdot \det(\breve{\Theta}) - S_{ab} \cdot \det(\Theta)\right| \leq C_{\text{oracle}} \cdot \sqrt{\frac{k_n^2 \log(p_n)}{n}}.$$

The first part of Theorem 4.2, which bounds the distance between our estimator $\breve{\Theta}$ of $\Theta$ and the oracle estimator $\widetilde{\Theta}$, is established using bounds on $\widehat{\Sigma} - \Sigma$ in Section 4.3. Details are given in Appendix D.3.1. The second part of Theorem 4.2, which bounds the error in estimating variance, $\left|\breve{S}_{ab} \cdot \det(\breve{\Theta}) - S_{ab} \cdot \det(\Theta)\right|$, is treated in Appendix D.3.2.

### D.3.1 Bounds on $\breve{\Theta} - \widetilde{\Theta}$

We use our bounds on the covariance error, $\widehat{\Sigma} - \Sigma$, to derive a bound on the difference between our empirical estimator $\breve{\Theta}$ and the oracle estimator $\widetilde{\Theta}$ of $\Theta$. The bounds we give here are deterministic (given that our initial assumptions hold).

The following lemma is proved in Appendix F.7:

**Lemma D.4.** *If Assumptions 3.1, 3.2, and 3.3 hold, then with probability at least $1 - \delta_n$,*

$$||\breve{\Theta} - \widetilde{\Theta}||_\infty \leq C_{\text{submatrix}} \left( \frac{k_n \log(p_n)}{n} + ||\widehat{\Sigma} - \Sigma||_{\mathcal{S}_{k_n}} \cdot \sqrt{\frac{k_n \log(p_n)}{n}} \right),$$

*where $C_{\text{submatrix}}$ is a constant depending on $C_{\text{cov}}$, $C_{\text{est}}$, and $C_{\text{sparse}}$ but not on $(n, p_n, k_n)$.*

From this point on, we combine Corollary 4.8 and Lemma D.4 to obtain our probabilistic bound on $||\breve{\Theta} - \widetilde{\Theta}||_\infty$ (Theorem 4.2). Looking first at Corollary 4.8, and setting $\delta_1 = \delta_2 = \frac{1}{6p_n}$, we see that by the assumption $p_n \geq 2, k_n \geq 1$ and the assumption $n \geq 15 k_n \log(p_n)$ stated in Theorem 4.2, the conditions of Corollary 4.8 must hold. Then, with probability at least $1 - \delta_1 - \delta_2 = 1 - \frac{1}{3p_n}$,

$$||\widehat{\Sigma} - \Sigma||_{\mathcal{S}_{k_n}}$$
$$\leq \frac{\pi^2}{8} \cdot k_n \cdot \frac{4 \log\left(12 p_n \binom{p_n}{2}\right)}{n} + 2\pi \cdot 16(1 + \sqrt{5}) C_{\text{cov}} \cdot \sqrt{\frac{\log(12 p_n) + 2 k_n \log(12 p_n)}{n}}$$
$$\leq C_{\text{cov}} \cdot C' \cdot \sqrt{\frac{k_n \log(p_n)}{n}}, \quad \text{(D.6)}$$

---
[7]Note that the additional condition $n \geq 15 k_n \log(p_n)$ can be assumed to hold in our main result Theorem 3.5, since if this inequality does not hold, then the claim in Theorem 3.5 is trivial.

where we choose the universal constant $C' = 3\pi^2 + 2\pi \cdot 16(1 + \sqrt{5})\sqrt{15}$ which guarantees that the last inequality holds (using the assumptions $n \geqslant k_n \log(p_n)$, $p_n \geqslant 2$, and $k_n \geqslant 1$).

Now combining this result with Lemma D.4, we obtain

$$\|\check{\Theta} - \widetilde{\Theta}\|_\infty \leqslant \frac{k_n \log(p_n)}{n} \cdot C_{\mathsf{submatrix}} \left(1 + C_{\mathsf{cov}} \cdot C'\right).$$

Taking $C_{\mathsf{oracle}} \geqslant C_{\mathsf{submatrix}} (1 + C_{\mathsf{cov}} \cdot C')$, we have proved that the first bound in Theorem 4.2 holds with probability at least $1 - \frac{1}{3p_n} - \delta_n$.

### D.3.2 Variance estimate

For the second part of the theorem, that is, bounding the error in the variance estimate $\check{S}_{ab}$, we state this bound as a lemma and defer the proof to Appendix F.11, since we need to develop some additional technical results before treating this bound.

**Lemma D.5.** *Under the assumptions and definitions of Theorem 4.2, with probability at least* $1 - \frac{1}{6p_n}$, *if* $n \geqslant k_n^2 \log(p_n)$, *on the event that the bounds* (3.1) *in Assumption 3.3 hold,*

$$\left| \check{S}_{ab} \cdot \det(\check{\Theta}) - S_{ab} \cdot \det(\Theta) \right| \leqslant C_{\mathsf{oracle}} \cdot \sqrt{\frac{k_n^2 \log(p_n)}{n}} \ .$$

Combining this lemma with the work above, and using Assumption 3.3, we have proved that both bounds stated in Theorem 4.2 hold with probability at least $1 - \frac{1}{2p_n} - \delta_n$, as desired.

## D.4 Proof of Theorem 3.5: main result

We now prove our main result, Theorem 3.5.

*Proof of Theorem 3.5.* Recall that our goal is to prove that $\frac{\sqrt{n}(\check{\Omega}_{ab} - \Omega_{ab})}{\check{S}_{ab}}$ converges to the $N(0,1)$ distribution. Recalling that $\Theta = (\Omega_{ab,ab})^{-1}$ and using the formula for a $2 \times 2$ matrix inverse, we separate this random variable into several terms:

$$\frac{\sqrt{n}(\check{\Omega}_{ab} - \Omega_{ab})}{\check{S}_{ab}}$$

$$= \frac{\sqrt{n}\left(\frac{-\check{\Theta}_{ab}}{\det(\check{\Theta})} - \frac{-\Theta_{ab}}{\det(\Theta)}\right)}{\check{S}_{ab}} = \frac{\sqrt{n}\left(-\check{\Theta}_{ab} + \Theta_{ab} \cdot \frac{\det(\check{\Theta})}{\det(\Theta)}\right)}{\check{S}_{ab} \cdot \det(\check{\Theta})}$$

$$= \frac{\sqrt{n}\left(\Theta_{ab} - \widetilde{\Theta}_{ab} + \widetilde{\Theta}_{ab} - \check{\Theta}_{ab} - \Theta_{ab} \cdot \left(1 - \frac{\det(\check{\Theta})}{\det(\Theta)}\right)\right)}{\check{S}_{ab} \cdot \det(\check{\Theta})}$$

$$= \left[ -\frac{\sqrt{n}\left(\widetilde{\Theta}_{ab} - \Theta_{ab}\right)}{S_{ab} \cdot \det(\Theta)} + \frac{\sqrt{n}\left(\widetilde{\Theta}_{ab} - \check{\Theta}_{ab}\right)}{S_{ab} \cdot \det(\Theta)} + \frac{\sqrt{n} \cdot \Omega_{ab} \cdot \left(\det(\Theta) - \det(\check{\Theta})\right)}{S_{ab} \cdot \det(\Theta)} \right]$$

$$\times \left[ 1 + \frac{\check{S}_{ab} \cdot \det(\check{\Theta}) - S_{ab} \cdot \det(\Theta)}{S_{ab} \cdot \det(\Theta)} \right]^{-1} .$$

To show that $\frac{\sqrt{n}(\check{\Omega}_{ab} - \Omega_{ab})}{\check{S}_{ab}}$ converges to the standard normal distribution, we will can apply Lemma D.3 (stated in Appendix D.2). In order to apply this lemma and obtain the desired result, we assemble the following pieces:

First, the variable $A := -\frac{\sqrt{n}(\widetilde{\Theta}_{ab} - \Theta_{ab})}{S_{ab} \cdot \det(\Theta)}$ satisfies $\sup_{t \in \mathbb{R}} |\mathbb{P}\{A \leqslant t\} - \Phi(t)| \leqslant C_{\mathsf{normal}} \cdot \frac{k_n \log(p_n)}{\sqrt{n}} +$ $\frac{1}{2p_n} =: \epsilon_A$, as shown in Theorem 4.1.

Second, we define variables $B := \frac{\sqrt{n}(\widetilde{\Theta}_{ab} - \breve{\Theta}_{ab})}{S_{ab} \cdot \det(\Theta)} + \frac{\sqrt{n} \cdot \Omega_{ab} \cdot (\det(\Theta) - \det(\breve{\Theta}))}{S_{ab} \cdot \det(\Theta)}$ and $C := \frac{\breve{S}_{ab} \cdot \det(\breve{\Theta}) - S_{ab} \cdot \det(\Theta)}{S_{ab} \cdot \det(\Theta)}$, and set

$$\delta_B = \frac{k_n \log(p_n)}{\sqrt{n}} \cdot \left( \frac{C_{\text{oracle}} + 4C_{\text{cov}}^2 C_{\text{oracle}} + 2C_{\text{cov}} C_{\text{oracle}}^2}{C_{\text{variance}}} \right)$$

and

$$\delta_C = \frac{C_{\text{oracle}}}{C_{\text{variance}}} \cdot \sqrt{\frac{k_n^2 \log(p_n)}{n}} \ .$$

We now show that, by Theorem 4.2, with probability at least $1 - \frac{1}{2p_n} - \delta_n =: 1 - \epsilon_{BC}$ it holds that $|B| \leq \delta_B$ and $|C| \leq \delta_C$. For the variable $C$, this is a trivial consequence of the bound on $\left| \breve{S}_{ab} \cdot \det(\breve{\Theta}) - S_{ab} \cdot \det(\Theta) \right|$ in Theorem 4.2 combined with the lower bound $S_{ab} \cdot \det(\Theta) \geq C_{\text{variance}}$ from Theorem 4.1.

Now we turn to the bound on $B$. To prove this bound, observe that $||\breve{\Theta} - \widetilde{\Theta}||_\infty \leq C_{\text{oracle}} \cdot \frac{k_n \log(p_n)}{n}$ by Theorem 4.2 (with the stated probability). We also have

$$|B| = \left| \frac{\sqrt{n}\left(\widetilde{\Theta}_{ab} - \breve{\Theta}_{ab}\right)}{S_{ab} \cdot \det(\Theta)} \right| \leq \sqrt{n} \cdot \frac{1}{S_{ab} \cdot \det(\Theta)} \cdot ||\breve{\Theta} - \widetilde{\Theta}||_\infty \leq \frac{\sqrt{n}}{C_{\text{variance}}} \cdot ||\breve{\Theta} - \widetilde{\Theta}||_\infty \ ,$$

where the last step follows from Theorem 4.1. And,

$$\left| \det(\breve{\Theta}) - \det(\Theta) \right| = \left| \left( \breve{\Theta}_{aa} \breve{\Theta}_{bb} - \breve{\Theta}_{ab}^2 \right) - \left( \Theta_{aa} \Theta_{bb} - \Theta_{ab}^2 \right) \right|$$
$$\leq 4||\Theta||_\infty ||\breve{\Theta} - \Theta||_\infty + 2||\breve{\Theta} - \Theta||_\infty^2$$

and

$$|\Omega_{ab}| \leq \lambda_{\max}(\Omega) = (\lambda_{\min}(\Sigma))^{-1} \leq C_{\text{cov}} \ .$$

Therefore,

$$\left| \frac{\sqrt{n} \cdot \Omega_{ab} \cdot \left( \det(\Theta) - \det(\breve{\Theta}) \right)}{S_{ab} \cdot \det(\Theta)} \right| \leq \sqrt{n} \cdot \frac{|\Omega_{ab}|}{S_{ab} \cdot \det(\Theta)} \cdot \left( 4||\Theta||_\infty ||\breve{\Theta} - \Theta||_\infty + 2||\breve{\Theta} - \Theta||_\infty^2 \right)$$
$$\leq \sqrt{n} \cdot \frac{C_{\text{cov}}}{C_{\text{variance}}} \cdot \left( 4C_{\text{cov}} ||\breve{\Theta} - \Theta||_\infty + 2||\breve{\Theta} - \Theta||_\infty^2 \right) \ ,$$

where the last step follows from Theorem 4.1 along with the fact that

$$||\Theta||_\infty \leq \lambda_{\max}(\Theta) = (\lambda_{\min}(\Omega_{ab,ab}))^{-1} \leq (\lambda_{\min}(\Omega))^{-1} = \lambda_{\max}(\Sigma) \leq C_{\text{cov}} \ .$$

Combining everything, we have

$$|B| \leq \sqrt{n} \cdot \frac{1}{C_{\text{variance}}} \cdot ||\breve{\Theta} - \widetilde{\Theta}||_\infty + \sqrt{n} \cdot \frac{C_{\text{cov}}}{C_{\text{variance}}} \cdot \left( 4C_{\text{cov}} ||\breve{\Theta} - \Theta||_\infty + 2||\breve{\Theta} - \Theta||_\infty^2 \right)$$
$$\leq \frac{k_n \log(p_n)}{\sqrt{n}} \left[ \frac{C_{\text{oracle}} + 4C_{\text{cov}}^2 C_{\text{oracle}} + 2C_{\text{cov}} C_{\text{oracle}}^2 \cdot \frac{k_n \log(p_n)}{n}}{C_{\text{variance}}} \right] \ .$$

If $n < k_n \log(p_n)$, then the main result in Theorem 3.5 holds trivially. Assuming then that $n \geq k_n \log(p_n)$, we have proved the desired bound on $|B|$.

Given these convergence results, we apply Lemma D.3 to obtain the following result:

$$\sup_{t \in \mathbb{R}} \left| \mathbb{P}\left\{ \frac{\sqrt{n}(\breve{\Omega}_{ab} - \Omega_{ab})}{\breve{S}_{ab}} \leqslant t \right\} - \Phi(t) \right| \leqslant \delta_B + \frac{\delta_C}{1 - \delta_C} + \epsilon_A + \epsilon_{BC}$$

$$= \frac{k_n \log(p_n)}{\sqrt{n}} \cdot \left( \frac{C_{\text{oracle}} + 4C_{\text{cov}}^2 C_{\text{oracle}} + 2C_{\text{cov}} C_{\text{oracle}}^2}{C_{\text{variance}}} \right) +$$

$$\frac{\frac{C_{\text{oracle}}}{C_{\text{variance}}} \cdot \sqrt{\frac{k_n^2 \log(p_n)}{n}}}{1 - \frac{C_{\text{oracle}}}{C_{\text{variance}}} \cdot \sqrt{\frac{k_n^2 \log(p_n)}{n}}} + C_{\text{normal}} \cdot \frac{k_n \log(p_n)}{\sqrt{n}} + \frac{1}{2p_n} + \frac{1}{2p_n} + \delta_n \ .$$

If $\frac{C_{\text{oracle}}}{C_{\text{variance}}} \cdot \sqrt{\frac{k_n^2 \log(p_n)}{n}} > \frac{1}{2}$, then the result of Theorem 3.5 holds trivially, and so assuming that this is not the case, we have

$$\sup_{t \in \mathbb{R}} \left| \mathbb{P}\left\{ \frac{\sqrt{n}(\breve{\Omega}_{ab} - \Omega_{ab})}{\breve{S}_{ab}} \leqslant t \right\} - \Phi(t) \right| \leqslant C_{\text{converge}} \cdot \sqrt{\frac{k_n^2 \log^2(p_n)}{n}} + \frac{1}{p_n} + \delta_n \ ,$$

where

$$C_{\text{converge}} := \frac{C_{\text{oracle}} + 4C_{\text{cov}}^2 C_{\text{oracle}} + 2C_{\text{cov}} C_{\text{oracle}}^2}{C_{\text{variance}}} + \frac{2C_{\text{oracle}}}{C_{\text{variance}}} + C_{\text{normal}} \ .$$

$\square$

## E  Accuracy of the initial Lasso estimator

**Corollary 3.8.** *Suppose that Assumption 3.1 holds. Assume additionally that the columns $\Omega_a, \Omega_b$ of the true inverse covariance $\Omega = \Sigma^{-1}$ are $k_n$-sparse. Then there exist constants $C_{\text{sample}}, C_{\text{Lasso}}$, depending on $C_{\text{cov}}$ but not on $(n, k_n, p_n)$, such that if $n \geqslant C_{\text{sample}} k_n \log(p_n)$ then, with probability at least $1 - \frac{1}{2p_n}$, any local minimizer $\breve{\gamma}_a$ of the objective function*

$$\frac{1}{2} \gamma^\top \widehat{\Sigma}_I \gamma - \gamma^\top \widehat{\Sigma}_{Ia} + \lambda ||\gamma||_1$$

*over the set $\{\gamma \in \mathbb{R}^I : ||\gamma||_1 \leqslant C_{\text{cov}} \sqrt{2k_n}\}$ satisfies*

$$||\breve{\gamma}_a - \gamma_a||_2 \leqslant 3\sqrt{2} C_{\text{cov}} \lambda \sqrt{k_n} \text{ and } ||\breve{\gamma}_a - \gamma_a||_1 \leqslant 24 C_{\text{cov}} \lambda \sqrt{k_n} \ .$$

*where we choose $\lambda = C_{\text{Lasso}} \cdot \sqrt{\frac{\log(p_n)}{n}}$. The same result holds for estimating $\gamma_b$.*

*Proof of Corollary 3.8.* Define

$$A = \widehat{\Sigma}_I, \ z = \widehat{\Sigma}_{Ia}, \ x^\star = \gamma_a, \ p = p_n - 2, \ k = k_n \ .$$

Now we apply Theorem 3.7 to this sparse recovery problem. In order to do so, we need to check that the conditions (3.3), (3.4), and (3.5) hold, and that $\gamma_a$ is feasible under the condition $||\gamma||_1 \leqslant R$. Once these conditions are satisfied, the result of Theorem 3.7 can be applied to this setting.

**Feasibility of $\gamma_a$** Define $R = C_{\text{cov}} \sqrt{2k_n}$. As proved in (D.1), $||\gamma_a||_2 \leqslant C_{\text{cov}}$, and furthermore $||\gamma_a||_0 \leqslant ||\Omega_a||_0 + ||\Omega_b||_0 \leqslant 2k_n$ (this is true because $\gamma_a = -\Omega_{I,ab}\Theta_{ab,a}$ by (D.2)). Therefore, $||\gamma_a||_1 \leqslant C_{\text{cov}} \sqrt{2k_n} = R$.

**Condition (3.3) (restricted strong convexity)** Now we need to check that the restricted strong convexity conditions (3.3) hold for our matrix $A = \hat{\Sigma}_I$. We will show that Corollary 4.8 implies that there exists a constant $C_{\mathsf{RSC}}$ depending only on $C_{\mathsf{cov}}$, such that if $n \geqslant 16 \log(p_n)$, then with probability at least $1 - \frac{1}{8p_n}$, for all $v \in \mathbb{R}^{p_n}$,

$$\left| v^\top \left( \hat{\Sigma}_I - \Sigma_I \right) v \right| \leqslant \frac{1}{2C_{\mathsf{cov}}} \left( ||v||_2^2 + ||v||_1^2 \cdot \frac{C_{\mathsf{RSC}} \log(p_n)}{n} \right) .$$

To see this why this holds, set $k = \frac{n}{C_{\mathsf{RSC}} \log(p_n)}$, and apply Lemma 4.9 to obtain

$$\left| v^\top \left( \hat{\Sigma}_I - \Sigma_I \right) v \right| \leqslant (||v||_2 + ||v||_1/\sqrt{k})^2 ||\hat{\Sigma} - \Sigma||_{\mathcal{S}_k} \leqslant 2(||v||_2^2 + ||v||_1^2/k) ||\hat{\Sigma} - \Sigma||_{\mathcal{S}_k} .$$

Then, applying Corollary 4.8 with this value of $k$ and with $\delta_1 = \delta_2 = \frac{1}{16p_n}$, we see that with probability at least $1 - \frac{1}{8p_n}$,

$$||\hat{\Sigma} - \Sigma||_{\mathcal{S}_k} \leqslant \frac{1}{4C_{\mathsf{cov}}},$$

as long as we set the constant $C_{\mathsf{RSC}}$ large enough.

Then, if this event holds, for all $v \in \mathbb{R}^I$ we have

$$v^\top \hat{\Sigma}_I v \geqslant v^\top \Sigma_I v - \left| v^\top \left( \hat{\Sigma}_I - \Sigma_I \right) v \right| \geqslant C_{\mathsf{cov}}^{-1} \cdot ||v||_2^2 - \left| v^\top \left( \hat{\Sigma}_I - \Sigma_I \right) v \right|$$

$$\geqslant \frac{1}{2C_{\mathsf{cov}}} \cdot ||v||_2^2 - \frac{C_{\mathsf{RSC}}}{2C_{\mathsf{cov}}} \cdot \frac{\log(p_n)}{n} \cdot ||v||_1^2 .$$

Therefore, with probability at least $1 - \frac{1}{8p_n}$, the restricted strong convexity condition (3.3) holds with

$$\alpha_1 = \frac{1}{2C_{\mathsf{cov}}} \text{ and } \tau_1 = \frac{C_{\mathsf{RSC}}}{2C_{\mathsf{cov}}} .$$

**Condition (3.5) (penalty parameter)** Below, we will prove that, with probability at least $1 - \frac{3}{8p_n}$,

$$||Ax^\star - z||_\infty = ||\hat{\Sigma}_I \gamma_a - \hat{\Sigma}_{Ia}||_\infty \leqslant \frac{\pi}{2} C_{\mathsf{feasible}} \sqrt{\frac{\log(p_n)}{n}} + \sqrt{\frac{\log(p_n)}{n}} \cdot \left[ \frac{1.5\sqrt{3}\pi^2 \sqrt{1 + C_{\mathsf{cov}}^2}}{\sqrt{C_{\mathsf{sample}}}} \right] , \quad \text{(E.1)}$$

for a constant $C_{\mathsf{feasible}}$ depending only on $C_{\mathsf{cov}}$, as long as we set

$$C_{\mathsf{sample}} \geqslant \left[ 16(1 + \sqrt{5}) C_{\mathsf{cov}} \sqrt{1 + C_{\mathsf{cov}}^2} \right]^2 .$$

Given that this is true, we now require that condition (3.5) holds, that is,

$$\max \left\{ 4||Ax^\star - z||_\infty, 4\alpha_1 \sqrt{\frac{\log(p)}{n}} \right\} \leqslant \lambda \leqslant \frac{\alpha_1}{6R} .$$

Define

$$C_{\mathsf{Lasso}} = \max \left\{ 4 \left[ \frac{\pi}{2} C_{\mathsf{feasible}} + \frac{1.5\sqrt{3}\pi^2 \sqrt{1 + C_{\mathsf{cov}}^2}}{\sqrt{C_{\mathsf{sample}}}} \right], \frac{2}{C_{\mathsf{cov}}} \right\} ,$$

Plugging in the bound (E.1), we see that the lower bound on $\lambda$ is satisfied for $\lambda = C_{\mathsf{Lasso}} \sqrt{\frac{\log(p_n)}{n}}$. To check the upper bound, we only need

$$\lambda = C_{\mathsf{Lasso}} \sqrt{\frac{\log(p_n)}{n}} \leqslant \frac{\alpha_1}{6R} = \frac{1}{2C_{\mathsf{cov}}^2 \sqrt{2k_n}} .$$

Assuming that

$$n \geqslant 8 C_{\mathsf{Lasso}}^2 C_{\mathsf{cov}}^4 \cdot k_n \log(p_n) , \quad \text{(E.2)}$$

then this follows directly. Therefore, (3.5) is satisfied with probability at least $1 - \frac{3}{8p_n}$.

**Condition (3.4) (sample size)** To satisfy (3.4), by plugging in the definitions of $R$, $\alpha_1$, and $\tau_1$ above, we see that it is sufficient to require

$$n \geq 64 C_{\mathsf{cov}}^2 C_{\mathsf{RSC}} \max\{1, 2C_{\mathsf{RSC}}\} \cdot k_n \log(p_n) \ . \tag{E.3}$$

**Conclusion** Combining all of our work above, we see that the conditions (3.3), (3.4), and (3.5), and the feasibility of $\gamma_a$, are all satisfied with probability at least $1 - \frac{1}{2p_n}$, as long as

$$n \geq C_{\mathsf{sample}} k_n \log(p_n)$$

for

$$C_{\mathsf{sample}} := \max\left\{16, \left[16(1+\sqrt{5})C_{\mathsf{cov}}\sqrt{1+C_{\mathsf{cov}}^2}\right]^2, 8C_{\mathsf{Lasso}}^2 C_{\mathsf{cov}}^4, 64 C_{\mathsf{cov}}^2 C_{\mathsf{RSC}} \max\{1, 2C_{\mathsf{RSC}}\}\right\} \ .$$

Therefore, applying Theorem 3.7, if these high probability events hold, then then for any $\breve{\gamma}_a$ that is a local minimizer of

$$L(x) = \frac{1}{2}\gamma^\top \widehat{\Sigma}_I \gamma - \gamma^\top \widehat{\Sigma}_{Ia} + \lambda ||\gamma||_1$$

over the set $\{\gamma \in \mathbb{R}^I : ||\gamma||_1 \leq 2C_{\mathsf{cov}} C_{\mathsf{sparse}} \sqrt{k_n}\}$, it holds that

$$||\breve{\gamma}_a - \gamma_a||_2 \leq \frac{1.5\lambda \cdot \sqrt{2k_n}}{\alpha_1} = 3\sqrt{2}C_{\mathsf{cov}}\lambda\sqrt{k_n} \text{ and } ||\breve{\gamma}_a - \gamma_a||_1 \leq \frac{6\lambda \cdot 2k_n}{\alpha_1} = 24 C_{\mathsf{cov}} \lambda \sqrt{k_n} \ .$$

By the same arguments, the same results hold for estimating $\gamma_b$.

**Proving (E.1)** Now we consider the term $||Ax^\star - z||_\infty = ||\widehat{\Sigma}_I \gamma_a - \widehat{\Sigma}_{Ia}||_\infty$. Since $\gamma_a = \Sigma_I^{-1} \Sigma_{Ia}$, we have

$$||\widehat{\Sigma}_I \gamma_a - \widehat{\Sigma}_{Ia}||_\infty = ||(\widehat{\Sigma}_I - \Sigma_I)\gamma_a - (\widehat{\Sigma}_{Ia} - \Sigma_{Ia})||_\infty = ||(\widehat{\Sigma} - \Sigma)u||_\infty \ ,$$

where $u \in \mathbb{R}^{p_n}$ is the fixed vector with

$$u_a = 1, u_b = 0, u_I = -\gamma_a \ .$$

By the Taylor expansion of $\widehat{\Sigma} - \Sigma$ (calculated as in (D.3)), we have

$$||(\widehat{\Sigma} - \Sigma)u||_\infty = \max_j \left|\mathbf{e}_j^\top (\widehat{\Sigma} - \Sigma) u\right|$$

$$\leq \frac{\pi}{2} \max_j \left|\mathbf{e}_j^\top \left(\cos\left(\frac{\pi}{2}T\right) \circ (\widehat{T} - T)\right) u\right| + \frac{\pi^2}{8}\left|\mathbf{e}_j^\top \left(\sin\left(\frac{\pi}{2}\overline{T}\right) \circ (\widehat{T} - T) \circ (\widehat{T} - T)\right) u\right|$$

$$\leq \frac{\pi}{2} \max_j \left|\mathbf{e}_j^\top \left(\cos\left(\frac{\pi}{2}T\right) \circ (\widehat{T} - T)\right) u\right| + \frac{\pi^2}{8}||u||_1 ||\widehat{T} - T||_\infty^2 \ . \tag{E.4}$$

Next we bound each term in this final expression (E.4) separately. Beginning with the second term, by (D.1), we know that $||u||_1 \leq \sqrt{||u||_0}||u||_2 \leq \sqrt{1+2k_n} \cdot \sqrt{1+C_{\mathsf{cov}}^2} \leq \sqrt{k_n} \cdot \sqrt{3(1+C_{\mathsf{cov}}^2)}$, where to bound $||u||_0$ we use the calculation $||\gamma_a||_0 \leq 2k_n$ from before. Furthermore, by Lemma D.2, with probability at least $1 - \frac{1}{8p_n}$,

$$||\widehat{T} - T||_\infty \leq \sqrt{\frac{12\log(8p_n)}{n}} \leq \sqrt{\frac{48\log(p_n)}{n}} \ ,$$

using $p_n \geq 2$. Therefore, the second term in (E.4) is bounded as

$$\frac{\pi^2}{8}||u||_1 ||\widehat{T} - T||_\infty^2 \leq \frac{\pi^2}{8} \sqrt{k_n} \cdot \sqrt{3(1+C_{\mathsf{cov}}^2)} \frac{48\log(p_n)}{n} \leq \sqrt{\frac{\log(p_n)}{n}} \cdot \left[\frac{6\sqrt{3}\pi^2 \sqrt{1+C_{\mathsf{cov}}^2}}{\sqrt{C_{\mathsf{sample}}}}\right] \ , \tag{E.5}$$

where we use the assumption $n \geq C_{\mathsf{sample}} \log(p_n)$.

Next we turn to the first term in (E.4). In order to bound this term, we begin by stating two lemmas (proved in Appendix F.9):

**Lemma E.1.** *There exist vectors $a_1, a_2, \ldots$ and $b_1, b_2, \ldots$ with $||a_r||_\infty, ||b_r||_\infty \leq 1$ for all $r \geq 1$, and a sequence $t_1, t_2, \cdots \geq 0$ with $\sum_r t_r = 4$, such that $\cos\left(\frac{\pi}{2}T\right) = \sum_{r \geq 1} t_r a_r b_r^\top$.*

**Lemma E.2.** *For fixed $u, v$ with $||u||_2, ||v||_2 \leq 1$, for any $|t| \leq \frac{n}{4(1+\sqrt{5})C_{\text{cov}}}$,*

$$\mathbb{E}\left[\exp\left(t \cdot u^\top (\widehat{T} - T)v\right)\right] \leq \exp\left(\frac{[4(1+\sqrt{5})]^2 t^2 \cdot C_{\text{cov}}^2}{n}\right).$$

By Lemma E.1, we can write

$$\cos\left(\frac{\pi}{2}T\right) = \sum_{r \geq 1} t_r \cdot a_r b_r^\top,$$

where $t_r \geq 0$, $\sum_r t_r = 4$, and $||a_r||_\infty, ||b_r||_\infty \leq 1$. Then

$$\mathbf{e}_j^\top \left(\cos\left(\frac{\pi}{2}T\right) \circ (\widehat{T} - T)\right) u = \left\langle \cos\left(\frac{\pi}{2}T\right) \circ \mathbf{e}_j u^\top, \widehat{T} - T \right\rangle = \sum_{r \geq 1} t_r \cdot (a_r \circ \mathbf{e}_j)^\top (\widehat{T} - T)(b_r \circ u).$$

Note that

$$||a_r \circ \mathbf{e}_j||_2 \leq ||a_r||_\infty \cdot ||\mathbf{e}_j||_2 \leq 1$$

and, by (D.1),

$$||b_r \circ u||_2 \leq ||b_r||_\infty \cdot ||u||_2 \leq \sqrt{1 + C_{\text{cov}}^2}.$$

Then for any $|t| \leq \frac{n}{16(1+\sqrt{5})C_{\text{cov}}\sqrt{1+C_{\text{cov}}^2}}$,

$$\mathbb{E}\left[\exp\left\{t \cdot \mathbf{e}_j^\top \left(\cos\left(\frac{\pi}{2}T\right) \circ (\widehat{T} - T)\right) u\right\}\right]$$

$$= \mathbb{E}\left[\exp\left\{\sum_{r \geq 1} t_r \left[t \cdot (a_r \circ \mathbf{e}_j)^\top (\widehat{T} - T)(b_r \circ u)\right]\right\}\right]$$

$$\leq \sum_r \frac{t_r}{4} \mathbb{E}\left[\exp\left\{4\left[t \cdot (a_r \circ \mathbf{e}_j)^\top (\widehat{T} - T)(b_r \circ u)\right]\right\}\right] \quad \text{(by Jensen's inequality)}$$

$$\leq \sum_r \frac{t_r}{4} \exp\left(\frac{[4(1+\sqrt{5})]^2 16t^2 \cdot C_{\text{cov}}^2 (1 + C_{\text{cov}}^2)}{n}\right) \quad \text{(by Lemma E.2)}$$

$$= \exp\left(\frac{[4(1+\sqrt{5})]^2 16t^2 \cdot C_{\text{cov}}^2 (1 + C_{\text{cov}}^2)}{n}\right).$$

Observe that we can set $t = \pm\sqrt{n \log(p_n)}$, which satisfies $|t| \leq \frac{n}{16(1+\sqrt{5})C_{\text{cov}}\sqrt{1+C_{\text{cov}}^2}}$ as long as we set $C_{\text{sample}} \geq \left[16(1+\sqrt{5})C_{\text{cov}}\sqrt{1+C_{\text{cov}}^2}\right]^2$, due to the assumption $n \geq C_{\text{sample}} \log(p_n)$. Then, we see that for any $C > 0$,

$$\mathbb{P}\left\{\left|\mathbf{e}_j^\top \left(\cos\left(\frac{\pi}{2}T\right) \circ (\widehat{T} - T)\right) u\right| > C\sqrt{\frac{\log(p_n)}{n}}\right\}$$

$$\leq \mathbb{E}\left[e^{\sqrt{n\log(p_n)} \cdot \mathbf{e}_j^\top \left(\cos\left(\frac{\pi}{2}T\right)\circ(\widehat{T}-T)\right)u - \sqrt{n\log(p_n)} \cdot C\sqrt{\frac{\log(p_n)}{n}}}\right]$$

$$+ \mathbb{E}\left[e^{-\sqrt{n\log(p_n)} \cdot \mathbf{e}_j^\top \left(\cos\left(\frac{\pi}{2}T\right)\circ(\widehat{T}-T)\right)u - \sqrt{n\log(p_n)} \cdot C\sqrt{\frac{\log(p_n)}{n}}}\right]$$

$$\leq 2\exp\left(\frac{[4(1+\sqrt{5})]^2 16(\sqrt{n\log(p_n)})^2 \cdot C_{\text{cov}}^2 (1+C_{\text{cov}}^2)}{n}\right) \cdot \exp\left\{-\sqrt{n\log(p_n)} \cdot C\sqrt{\frac{\log(p_n)}{n}}\right\}$$

$$\leq 2p_n^{-\left(C - [4(1+\sqrt{5})]^2 \cdot 16 C_{\text{cov}}^2 (1+C_{\text{cov}}^2)\right)} = 2p_n^{-5} \leq \frac{1}{4p_n^2},$$

where we set $C = C_{\text{feasible}} := 5 + \left[4(1+\sqrt{5})\right]^2 \cdot 16 C_{\text{cov}}^2 (1 + C_{\text{cov}}^2)$. Therefore,

$$\mathbb{P}\left\{\max_j \left|\mathbf{e}_j^\top \left(\cos\left(\frac{\pi}{2}T\right) \circ (\widehat{T} - T)\right) u\right| > C_{\text{feasible}} \sqrt{\frac{\log(p_n)}{n}}\right\} \leqslant \frac{1}{4p_n} \, . \tag{E.6}$$

Combining (E.5) and (E.6), and returning to (E.4), we have

$$\|\widehat{\Sigma}_I \gamma_a - \widehat{\Sigma}_{Ia}\|_\infty \leqslant \frac{\pi}{2} C_{\text{feasible}} \sqrt{\frac{\log(p_n)}{n}} + \sqrt{\frac{\log(p_n)}{n}} \cdot \left[\frac{1.5\sqrt{3}\pi^2 \sqrt{1+C_{\text{cov}}^2}}{\sqrt{C_{\text{sample}}}}\right] ,$$

with probability at least $1 - \frac{3}{8p_n}$. This proves (E.1). $\square$

# F Proofs of lemmas

## F.1 Proof of the normal convergence lemma

**Lemma D.3.** *Let $A, B, C$ be random variables such that*

$$\sup_{t \in \mathbb{R}} |\mathbb{P}\{A \leqslant t\} - \Phi(t)| \leqslant \epsilon_A \quad \text{and} \quad \mathbb{P}\{|B| \leqslant \delta_B, |C| \leqslant \delta_C\} \geqslant 1 - \epsilon_{BC} ,$$

*where $\epsilon_A, \epsilon_{BC}, \delta_B, \delta_C \in [0,1)$. Then the variable $(A+B) \cdot (1+C)^{-1}$ converges to a standard normal distribution with rate*

$$\sup_{t \in \mathbb{R}} \left|\mathbb{P}\left\{(A+B) \cdot (1+C)^{-1} \leqslant t\right\} - \Phi(t)\right| \leqslant \delta_B + \frac{\delta_C}{1-\delta_C} + \epsilon_A + \epsilon_{BC} .$$

*Proof of Lemma D.3.* First, define truncated versions of $B$ and $C$:

$$\widetilde{B} = \text{sign}(B) \cdot \min\{|B|, \delta_B\}, \ \widetilde{C} = \text{sign}(C) \cdot \min\{|C|, \delta_C\} .$$

Then, for any $t \in \mathbb{R}$,

$$\left|\mathbb{P}\left\{(A+B) \cdot (1+C)^{-1} \leqslant t\right\} - \mathbb{P}\left\{(A+\widetilde{B}) \cdot (1+\widetilde{C})^{-1} \leqslant t\right\}\right|$$
$$\leqslant \mathbb{P}\left\{B \neq \widetilde{B} \text{ or } C \neq \widetilde{C}\right\} \leqslant \epsilon_{BC} .$$

Note that $|\widetilde{B}| \leqslant \delta_B$ and $|\widetilde{C}| \leqslant \delta_C$ with probability 1.

Next, fix any $t \geqslant 0$ and suppose that $A \leqslant t(1-\delta_C) - \delta_B$. Then

$$(A + \widetilde{B}) \cdot (1+\widetilde{C})^{-1} \leqslant ((t(1-\delta_C) - \delta_B) + \delta_B) \cdot (1-\delta_C)^{-1} = t ,$$

and so

$$\mathbb{P}\left\{(A+\widetilde{B}) \cdot (1+\widetilde{C})^{-1} \leqslant t\right\} \geqslant \mathbb{P}\{A \leqslant t(1-\delta_C) - \delta_B\}$$
$$\geqslant \Phi(t(1-\delta_C) - \delta_B) - \epsilon_A$$
$$= \Phi(t) - \mathbb{P}\{t(1-\delta_C) - \delta_B < N(0,1) < t(1-\delta_C)\} - \mathbb{P}\{t(1-\delta_C) < N(0,1) < t\} - \epsilon_A$$
$$\geqslant \Phi(t) - \delta_B - \mathbb{P}\{t(1-\delta_C) < N(0,1) < t\} - \epsilon_A \quad \text{(since the density of } N(0,1) \text{ is } \leqslant \frac{1}{\sqrt{2\pi}} \leqslant 1\text{)}$$
$$\geqslant \Phi(t) - \delta_B - \frac{\delta_C}{1-\delta_C} - \epsilon_A. \quad \text{(by applying Lemma F.1 stated below)}$$

To prove the reverse bound, suppose that $(A + \widetilde{B}) \cdot (1 + \widetilde{C})^{-1} \leq t$. Then

$$A = (A + \widetilde{B}) \cdot (1 + \widetilde{C})^{-1} \cdot (1 + \widetilde{C}) - \widetilde{B} \leq t(1 + \delta_C) + \delta_B \ .$$

Therefore,

$$\mathbb{P}\left\{(A + \widetilde{B}) \cdot (1 + \widetilde{C})^{-1} \leq t\right\} \leq \mathbb{P}\{A \leq t(1 + \delta_C) + \delta_B\}$$
$$\leq \Phi\left(t(1 + \delta_C) + \delta_B\right) + \epsilon_A$$
$$= \Phi(t) + \mathbb{P}\{t(1+\delta_C) < N(0,1) < t(1+\delta_C) + \delta_B\} + \mathbb{P}\{t < N(0,1) < t(1+\delta_C)\} + \epsilon_A \ .$$
$$\leq \Phi(t) + \delta_B + \mathbb{P}\{t < N(0,1) < t(1+\delta_C)\} + \epsilon_A \quad \left(\text{since the density of } N(0,1) \text{ is} \leq \frac{1}{\sqrt{2\pi}} \leq 1\right)$$
$$\leq \Phi(t) + \delta_B + \delta_C + \epsilon_A. \quad \text{(by applying Lemma F.1 stated below)}$$

Therefore, for all $t \geq 0$,

$$\left|\mathbb{P}\left\{(A + \widetilde{B}) \cdot (1 + \widetilde{C})^{-1} \leq t\right\} - \Phi(t)\right| \leq \delta_B + \frac{\delta_C}{1 - \delta_C} + \epsilon_A \ .$$

By identical arguments, we can prove the same for $t \leq 0$. $\square$

**Lemma F.1.** *For any $0 \leq a \leq b$,*

$$\mathbb{P}\{a < N(0,1) < b\} \leq \left(\frac{b}{a} - 1\right) \cdot \frac{1}{\sqrt{2\pi e}} \leq \left(\frac{b}{a} - 1\right) \ .$$

*Proof.*

$$\mathbb{P}\{a < N(0,1) < b\} = \int_{t=a}^{b} \frac{1}{\sqrt{2\pi}} e^{-t^2/2} \, \mathsf{d}t$$
$$\leq (b-a) \cdot \frac{1}{\sqrt{2\pi}} e^{-a^2/2}$$
$$= \left(\frac{b}{a} - 1\right) \cdot \frac{1}{\sqrt{2\pi}} \cdot a \cdot e^{-a^2/2}$$
$$\leq \left(\frac{b}{a} - 1\right) \cdot \frac{1}{\sqrt{2\pi e}} \ ,$$

where the last step holds because $\sup_{t>0}\{t \cdot e^{-t^2/2}\} = \frac{1}{\sqrt{e}}$. $\square$

### F.2 Sign vector of a transelliptical distribution

**Lemma 4.4.** *Let $X, X' \overset{iid}{\sim} \mathsf{TE}(\Sigma, \xi; f_1, \ldots, f_p)$. Suppose that $\Sigma$ is positive definite, and that $\xi > 0$ with probability 1. Then $\text{sign}(X - X')$ is equal in distribution to $\text{sign}(Z)$, where $Z \sim N(0, \Sigma)$.*

*Proof of Lemma 4.4.* First, since the $f_j$'s are strictly monotone, we see that $\text{sign}(X - X')$ has the same distribution regardless of the choice of the $f_j$'s. Therefore it suffices to consider the case that the $f_j$'s are each the identity function, and so $X, X' \overset{iid}{\sim} \mathsf{E}(\mathbf{0}, \Sigma, \xi)$, that is, a zero-mean elliptical distribution. In this case, by Lindskog et al. [2003, Lemma 1], $X - X' \sim \mathsf{E}(\mathbf{0}, \Sigma, \zeta)$ where the distribution of the random variable $\zeta \geq 0$ obeys $\varphi_\zeta(t) = \varphi_\xi(t)^2$, where $\varphi_\zeta$ and $\varphi_\xi$ are the characteristic functions of $\zeta$ and $\xi$, respectively. Note that for two independent copies $\xi_1, \xi_2 \overset{iid}{\sim} \xi$, we have $\varphi_{\xi_1+\xi_2} = \varphi_{\xi_1} \cdot \varphi_{\xi_2} = \varphi_\xi^2 = \varphi_\zeta$, and therefore, $\zeta \overset{\mathcal{D}}{=} \xi_1 + \xi_2$. Since $\xi > 0$ with probability 1, this proves that $\zeta > 0$ with probability 1.

Next take $Z \sim N(0, \Sigma)$. Then $\frac{\Sigma^{-1/2}Z}{||\Sigma^{-1/2}Z||_2}$ is uniformly distributed on the unit sphere, and so

$$\zeta \cdot \Sigma^{1/2} \cdot \frac{\Sigma^{-1/2}Z}{||\Sigma^{-1/2}Z||_2} \sim \mathsf{E}(\mathbf{0}, \Sigma, \zeta) ,$$

which is the distribution of $X - X'$. Using the fact that $\zeta > 0$ with probability 1, we see that $\text{sign}(X - X')$ is equal in distribution to

$$\text{sign}\left(\zeta \cdot \Sigma^{1/2} \cdot \frac{\Sigma^{-1/2}Z}{||\Sigma^{-1/2}Z||_2}\right) = \text{sign}(\Sigma^{1/2} \cdot \Sigma^{-1/2}Z) = \text{sign}(Z) ,$$

as desired. $\square$

### F.3 Proof of Lemma 4.7

**Lemma 4.7.** *The following bound holds deterministically: for any $k \geq 1$,*

$$||\widehat{\Sigma} - \Sigma||_{\mathcal{S}_k} \leq \frac{\pi^2}{8} \cdot k ||\widehat{T} - T||_\infty^2 + 2\pi ||\widehat{T} - T||_{\mathcal{S}_k} .$$

*Proof of Lemma 4.7.* By Taylor's theorem,

$$\widehat{\Sigma} = \Sigma + \frac{\pi}{2} \cos\left(\frac{\pi}{2}T\right) \circ \left(\widehat{T} - T\right) - \frac{\pi^2}{8} \sin\left(\frac{\pi}{2}\bar{T}\right) \circ \left(\widehat{T} - T\right) \circ \left(\widehat{T} - T\right)$$

where $\bar{T}$ has entries $\bar{\tau}_{ab} = (1 - t_{ab})\tau_{ab} + t_{ab}\widehat{\tau}_{ab}$, with $t_{ab} \in [0, 1]$ for each $a, b$. Taking any $u, v \in \mathcal{S}_k$, then,

$$\left|u^\top(\widehat{\Sigma} - \Sigma)v\right| \leq \frac{\pi}{2}\left|u^\top\left[\cos\left(\frac{\pi}{2}T\right) \circ \left(\widehat{T} - T\right)\right]v\right| + \frac{\pi^2}{8}\left|u^\top\left[\sin\left(\frac{\pi}{2}\bar{T}\right) \circ \left(\widehat{T} - T\right) \circ \left(\widehat{T} - T\right)\right]v\right| .$$

First, to bound the $\sin(\cdot)$ matrix term, note that

$$\left|u^\top\left[\sin\left(\frac{\pi}{2}\bar{T}\right) \circ \left(\widehat{T} - T\right) \circ \left(\widehat{T} - T\right)\right]v\right|$$
$$\leq ||u||_1 ||v||_1 \left\|\sin\left(\frac{\pi}{2}\bar{T}\right) \circ \left(\widehat{T} - T\right) \circ \left(\widehat{T} - T\right)\right\|_\infty \leq ||u||_1 ||v||_1 ||\widehat{T} - T||_\infty^2 ,$$

where the last step holds since the $\sin(\cdot)$ function lies in $[-1, 1]$. Furthermore, $||u||_1, ||v||_1 \leq \sqrt{k}$ for all $u, v \in \mathcal{S}_k$ by definition.

Next, we bound the $\cos(\cdot)$ matrix term. By Lemma E.1, we can express $\cos\left(\frac{\pi}{2}T\right)$ as a convex combination,

$$\cos\left(\frac{\pi}{2}T\right) = \sum_{r \geq 1} t_r \cdot a_r b_r^\top ,$$

where $a_r, b_r \in \mathbb{R}^p$ satisfy $||a_r||_\infty, ||b_r||_\infty \leq 1$ for all $r$, and $t_r \geq 0$ satisfy $\sum_r t_r = 4$. Furthermore, for $u, v \in \mathcal{S}_k$ and for each $r$, note that $u \circ a_r, v \circ b_r \in \mathcal{S}_k$ due to the bound on $||a_r||_\infty, ||b_r||_\infty$. Then

$$\left|u^\top\left[\cos\left(\frac{\pi}{2}T\right) \circ \left(\widehat{T} - T\right)\right]v\right| \leq \sum_{r \geq 1} t_r \left|u^\top\left[a_r b_r^\top \circ \left(\widehat{T} - T\right)\right]v\right|$$
$$= \sum_{r \geq 1} t_r \left|(u \circ a_r)^\top \left(\widehat{T} - T\right)(v \circ b_r)\right| \leq 4 \sup_{u', v' \in \mathcal{S}_k} \left|u'^\top(\widehat{T} - T)v'\right| = 4||\widehat{T} - T||_{\mathcal{S}_k} .$$

Using the definition $||\widehat{\Sigma} - \Sigma||_{\mathcal{S}_k} = \max_{u, v \in \mathcal{S}_k} u^\top(\widehat{\Sigma} - \Sigma)v$, this proves the lemma. $\square$

## F.4 Proof of Lemma 4.6

**Lemma 4.6.** *Suppose that $k \geq 1$ and $\delta \in (0,1)$ satisfy $\log(2/\delta) + 2k\log(12p) \leq n$. Then with probability at least $1 - \delta$ it holds that*

$$\|\widehat{T} - T\|_{\mathcal{S}_k} \leq 32(1 + \sqrt{5})\mathsf{C}(\Sigma) \cdot \sqrt{\frac{\log(2/\delta) + 2k\log(12p)}{n}} \;.$$

*Proof of Lemma 4.6.* This lemma is a straightforward combination of Lemma E.2 (stated in Appendix E) together with the following result:

**Lemma F.2** (Adapted from Lemma 5.1 and Theorem 5.2 of Baraniuk et al. [2008]). *Let $A$ be a random matrix satisfying*

$$\exp\{t \cdot u^\top A u\} \leq \exp\left\{\frac{c_1 t^2}{n}\right\} \text{ for all } |t| \leq c_0 n \text{ and all unit vectors } u \in \mathbb{R}^p \quad \text{(F.1)}$$

*for some constants $c_0, c_1$. Then for any $k \geq 1$ and any $\delta \in (0,1)$ satisfying*

$$\log(2/\delta) + k\log(12p) \leq n c_0^2 c_1 \;,$$

*with probability at least $1 - \delta$ it holds that*

$$|u^\top A u| \leq \sqrt{\frac{16 c_1}{n}\left(\log(2/\delta) + k\log(12p)\right)} \text{ for all } k\text{-sparse unit vectors } u \in \mathbb{R}^p. \quad \text{(F.2)}$$

Combined, Lemma F.2 (applied with $2k$ in place of $k$, with $c_0 = \frac{1}{4(1+\sqrt{5})C_{\text{cov}}}$, and $c_1 = c_0^{-1}$) and Lemma E.2 immediately yield the bound

$$\sup_{u \in \mathcal{S}_{2k}} \left|u^\top(\widehat{T} - T)u\right| \leq 16(1 + \sqrt{5})\mathsf{C}(\Sigma) \cdot \sqrt{\frac{\log(2/\delta) + 2k\log(12p)}{n}} \;,$$

with probability at least $1 - \delta$, as long as $\log(2/\delta) + 2k\log(12p) \leq n c_0^2 c_1 = n$. Next take any $u, v \in \mathcal{S}_k$. Then $\frac{u+v}{2}, \frac{u-v}{2} \in \mathcal{S}_{2k}$.

$$u^\top(\widehat{T} - T)v = \left(\frac{u+v}{2}\right)^\top (\widehat{T} - T)\left(\frac{u+v}{2}\right) - \left(\frac{u-v}{2}\right)^\top (\widehat{T} - T)\left(\frac{u+v}{2}\right)$$

$$\leq 16(1 + \sqrt{5})\mathsf{C}(\Sigma) \cdot \sqrt{\frac{\log(2/\delta) + 2k\log(12p)}{n}} \;.$$

This proves Lemma 4.6, as desired. $\square$

We next turn to the proof of Lemma F.2.

*Proof of Lemma F.2.* (Adapted from Lemma 5.1 and Theorem 5.2 of Baraniuk et al. [2008].) First fix any $S \subset [p]$ with $|S| = k$. Let $\epsilon = \sqrt{\frac{16 c_1}{n}\left(\log(2/\delta) + k\log(12p)\right)}$. Following the same arguments as in Baraniuk et al. [2008, Lemma 5.1], we can take a set $\mathcal{U} \subset \mathbb{R}^S$ of unit vectors, with $|\mathcal{U}| \leq 12^k$, such that

$$\sup_{\text{unit } u \in \mathbb{R}^S} \left|u^\top A u\right| \leq 2 \sup_{\widetilde{u} \in \mathcal{U}} \left|\widetilde{u}^\top A \widetilde{u}\right| \;.$$

Furthermore, for any fixed $\widetilde{u} \in \mathcal{U}$, for any $0 < t \leq c_0 n$,

$$\mathbb{P}\left\{\widetilde{u}^\top A \widetilde{u} > \epsilon/2\right\} \leq \mathbb{E}\left[t \cdot \widetilde{u}^\top A \widetilde{u} - t \cdot \epsilon/2\right]$$

$$\leq \exp\left(\frac{c_1 t^2}{n} - t \cdot \epsilon/2\right)$$

$$= \exp\left(-\frac{n\epsilon^2}{16 c_1}\right) \;,$$

where for the last step we set $t = \frac{n\epsilon}{4c_1} \leq c_0 n$. Similarly,

$$\mathbb{P}\left\{\widetilde{u}^\top A \widetilde{u} < -\epsilon/2\right\} \leq \exp\left(-\frac{n\epsilon^2}{16c_1}\right).$$

Therefore,

$$\mathbb{P}\left\{\sup_{\widetilde{u}\in\mathcal{U}} |\widetilde{u}^\top A \widetilde{u}| > \epsilon/2\right\} \leq 2 \cdot 12^k \cdot \exp\left(-\frac{n\epsilon^2}{16c_1}\right),$$

and so

$$\mathbb{P}\left\{\sup_{\text{unit } u\in\mathbb{R}^S} |u^\top A u| > \epsilon\right\} \leq 2 \cdot 12^k \cdot \exp\left(-\frac{n\epsilon^2}{16c_1}\right).$$

Finally, taking all $\binom{p}{k} \leq p^k$ choices for $S$, we see that

$$\mathbb{P}\left\{\sup_{u\in\mathcal{S}_k} \{u^\top A u\} \leq \epsilon\right\} \geq 1 - 2(12p)^k \cdot \exp\left(-\frac{n\epsilon^2}{16c_1}\right) = 1 - \delta.$$

$\square$

### F.5 Proof of Corollary 4.8

**Corollary 4.8.** *Take any $\delta_1, \delta_2 \in (0,1)$ and any $k \geq 1$ such that $\log(2/\delta_2) + 2k\log(12p) \leq n$. Then, with probability at least $1 - \delta_1 - \delta_2$, the following bound on $\widehat{\Sigma} - \Sigma$ holds:*

$$||\widehat{\Sigma} - \Sigma||_{\mathcal{S}_k} \leq$$

$$\frac{\pi^2}{8} \cdot k \cdot \frac{4\log\left(2\binom{p}{2}/\delta_1\right)}{n} + 2\pi \cdot 32(1+\sqrt{5})\mathsf{C}(\Sigma) \cdot \sqrt{\frac{\log(2/\delta_2) + 2k\log(12p)}{n}}.$$

*Proof of Corollary 4.8.* This proof is a straightforward combination of Lemmas D.2, 4.6, and 4.7. We have

$$||\widehat{\Sigma} - \Sigma||_{\mathcal{S}_k}$$
$$\leq \frac{\pi^2}{8} \cdot k ||\widehat{T} - T||_\infty^2 + 2\pi ||\widehat{T} - T||_{\mathcal{S}_k} \quad \text{by Lemma 4.7}$$
$$\leq \frac{\pi^2}{8} \cdot k \frac{4\log\left(2\binom{p_n}{2}/\delta_1\right)}{n} + 2\pi \cdot 32(1+\sqrt{5})\mathsf{C}(\Sigma) \cdot \sqrt{\frac{\log(2/\delta_2) + 2k\log(12p)}{n}},$$

where the last step holds by applying Lemma D.2 with $\delta = \delta_1$ and Lemma 4.6 with $\delta = \delta_2$. $\square$

### F.6 Proof of Lemma 4.9

**Lemma 4.9** (Based on Proposition 5 of Sun and Zhang [2012a]). *For any fixed matrix $M \in \mathbb{R}^{p\times p}$ and vectors $u, v \in \mathbb{R}^p$, and any $k \geq 1$,*

$$|u^\top M v| \leq \left(||u||_2 + ||u||_1/\sqrt{k}\right) \cdot \left(||v||_2 + ||v||_1/\sqrt{k}\right) \cdot \sup_{u',v'\in\mathcal{S}_k} |u'^\top M v'|.$$

*Proof of Lemma 4.9.* We first introduce two vector norms used in Proposition 5 of Sun and Zhang [2012a],

$$||w||_{(2,k)} = \sup_{|S|\leq k} ||w_S||_2$$

and its dual norm $||w||_{(2,k)}^*$. Note that

$$||w||_{(2,k)} = \sup_{z\in\mathcal{S}_k} |z^\top w|,$$

from the norm definition. We then have

$$|u^\top M v| \leq ||u||^*_{(2,k)} ||Mv||_{(2,k)}$$

by definition of the dual norm, and

$$||Mv||_{2,k} = \sup_{w \in \mathcal{S}_k} |w^\top M v| \leq \sup_{w \in \mathcal{S}_k} ||v||^*_{(2,k)} ||M^\top w||_{(2,k)} = \sup_{w \in \mathcal{S}_k} ||v||^*_{(2,k)} \left( \sup_{z \in \mathcal{S}_k} |z^\top M^\top w| \right),$$

and so combining everything,

$$|u^\top M v| \leq ||u||^*_{(2,k)} ||v||^*_{(2,k)} \cdot \sup_{w,z \in \mathcal{S}_k} |w^\top M z| = ||u||^*_{(2,k)} ||v||^*_{(2,k)} \cdot ||M||_{\mathcal{S}_k}.$$

Finally, Proposition 5(ii) of Sun and Zhang [2012a] (applied with $q = 2$ and $m = k$, in their notation), proves that, for any $w$,

$$||w||^*_{2,k} \leq ||w||_{(2,4k)} + ||w||_1/\sqrt{k} \leq ||w||_2 + ||w||_1/\sqrt{k},$$

where the last inequality follows trivially from the norm definition. This proves the lemma. $\square$

## F.7 Proof of Lemma D.4

**Lemma D.4.** *If Assumptions 3.1, 3.2, and 3.3 hold, then with probability at least $1 - \delta_n$,*

$$||\breve{\Theta} - \widetilde{\Theta}||_\infty \leq C_{\mathsf{submatrix}} \left( \frac{k_n \log(p_n)}{n} + ||\widehat{\Sigma} - \Sigma||_{\mathcal{S}_{k_n}} \cdot \sqrt{\frac{k_n \log(p_n)}{n}} \right),$$

*where $C_{\mathsf{submatrix}}$ is a constant depending on $C_{\mathsf{cov}}$, $C_{\mathsf{est}}$, and $C_{\mathsf{sparse}}$ but not on $(n, p_n, k_n)$.*

*Proof of Lemma D.4.* Choose any $c, d \in \{a, b\}$; we will bound the $(c, d)$th entry of the error, that is, $\left|\breve{\Theta}_{cd} - \widetilde{\Theta}_{cd}\right|$. Write $\Delta_c = \breve{\gamma}_c - \gamma_c$ for each $c = a, b$. First, with probability at least $1 - \delta_n$, the bounds on $\Delta_c, \Delta_d$ in Assumption 3.3 all hold; assume that this event occurs from this point on.

We have

$$\left|\breve{\Theta}_{cd} - \widetilde{\Theta}_{cd}\right|$$
$$= \left|\left(\widehat{\Sigma}_{cd} - \breve{\gamma}_c^\top \widehat{\Sigma}_{Id} - \widehat{\Sigma}_{Ic}^\top \breve{\gamma}_d + \breve{\gamma}_c^\top \widehat{\Sigma}_I \breve{\gamma}_d\right) - \left(\widehat{\Sigma}_{cd} - \gamma_c^\top \widehat{\Sigma}_{Id} - \widehat{\Sigma}_{Ic}^\top \gamma_d + \gamma_c^\top \widehat{\Sigma}_I \gamma_d\right)\right|$$
$$= \left|\breve{\gamma}_c^\top \widehat{\Sigma}_I \breve{\gamma}_d - \gamma_c^\top \widehat{\Sigma}_I \gamma_d - \Delta_c^\top \widehat{\Sigma}_{Id} - \widehat{\Sigma}_{Ic}^\top \Delta_d\right|$$
$$= \left|\Delta_c^\top \widehat{\Sigma}_I \gamma_d + \gamma_c^\top \widehat{\Sigma}_I \Delta_d + \Delta_c^\top \widehat{\Sigma}_I \Delta_d - \Delta_c^\top \widehat{\Sigma}_{Id} - \widehat{\Sigma}_{Ic}^\top \Delta_d\right|$$
$$\leq \left|\Delta_c^\top \Sigma_I \gamma_d + \gamma_c^\top \Sigma_I \Delta_d + \Delta_c^\top \Sigma_I \Delta_d - \Delta_c^\top \Sigma_{Id} - \Sigma_{Ic}^\top \Delta_d\right| \quad \text{(F.3)}$$
$$+ \left|\Delta_c^\top (\widehat{\Sigma}_I - \Sigma_I)\gamma_d + \gamma_c^\top (\widehat{\Sigma}_I - \Sigma_I)\Delta_d + \Delta_c^\top (\widehat{\Sigma}_I - \Sigma_I)\Delta_d - \Delta_c^\top (\widehat{\Sigma}_{Id} - \Sigma_{Id}) - (\widehat{\Sigma}_{Ic} - \Sigma_{Ic})^\top \Delta_d\right|.$$

Now we bound each of these terms. To bound the first term on the right-hand side of (F.3), we have

$$\left|\Delta_c^\top \Sigma_I \gamma_d + \gamma_c^\top \Sigma_I \Delta_d + \Delta_c^\top \Sigma_I \Delta_d - \Delta_c^\top \Sigma_{Id} - \Sigma_{Ic}^\top \Delta_d\right|$$
$$= \left|\Delta_c^\top \Sigma_I \Sigma_I^{-1} \Sigma_{Id} + \Sigma_{Ic}^\top \Sigma_I^{-1} \Sigma_I \Delta_d + \Delta_c^\top \Sigma_I \Delta_d - \Delta_c^\top \Sigma_{Id} - \Sigma_{Ic}^\top \Delta_d\right|$$
$$= \left|\Delta_c^\top \Sigma_I \Delta_d\right|$$
$$\leq ||\Delta_c||_2 \cdot ||\Delta_d||_2 \cdot ||\Sigma||_{\mathsf{op}}$$
$$\leq ||\Delta_c||_2 \cdot ||\Delta_d||_2 \cdot \mathsf{C}(\Sigma)$$
$$\leq C_{\mathsf{cov}} (C_{\mathsf{est}})^2 \frac{k_n \log(p_n)}{n},$$

where the next-to-last step holds because

$$||\Sigma_I|| \leq ||\Sigma|| = \lambda_{\min}(\Sigma) \cdot \mathsf{C}(\Sigma) \,,$$

and we must have $\lambda_{\min}(\Sigma) \leq 1$ because $\mathrm{diag}(\Sigma) = \mathbf{1}$, while the last step holds by Assumptions 3.1 and 3.3.

Finally, to bound the second term on the right-hand side of (F.3), we have

$$\left|\Delta_c^\top(\widehat{\Sigma}_I - \Sigma_I)\gamma_d + \gamma_c^\top(\widehat{\Sigma}_I - \Sigma_I)\Delta_d + \Delta_c^\top(\widehat{\Sigma}_I - \Sigma_I)\Delta_d - \Delta_c^\top(\widehat{\Sigma}_{Id} - \Sigma_{Id}) - (\widehat{\Sigma}_{Ic} - \Sigma_{Ic})^\top\Delta_d\right|$$

$$\leq \left|\Delta_c^\top(\widehat{\Sigma}_I - \Sigma_I)\gamma_d\right| + \left|\gamma_c^\top(\widehat{\Sigma}_I - \Sigma_I)\Delta_d\right| + \left|\Delta_c^\top(\widehat{\Sigma}_I - \Sigma_I)\Delta_d\right| + \left|\Delta_c^\top(\widehat{\Sigma}_{Id} - \Sigma_{Id})\right| + \left|(\widehat{\Sigma}_{Ic} - \Sigma_{Ic})^\top\Delta_d\right|$$

$$= \left|\Delta_c^\top(\widehat{\Sigma}_I - \Sigma_I)\gamma_d\right| + \left|\gamma_c^\top(\widehat{\Sigma}_I - \Sigma_I)\Delta_d\right| + \left|\Delta_c^\top(\widehat{\Sigma}_I - \Sigma_I)\Delta_d\right| + \left|\Delta_c^\top(\widehat{\Sigma}_I - \Sigma_I)\mathbf{e}_d\right| + \left|\mathbf{e}_c^\top(\widehat{\Sigma}_I - \Sigma_I)^\top\Delta_d\right|,$$

where $\mathbf{e}_c$ and $\mathbf{e}_d$ are the basis vectors in $\mathbb{R}^I$ corresponding to nodes $c$ and $d$. Now, for each vector $w \in \{\gamma_c, \gamma_d, \mathbf{e}_c, \mathbf{e}_d\}$ and each vector $z \in \{\Delta_c, \Delta_d\}$, we have

$$||w||_2 + ||w||_1/\sqrt{k} \leq C_{\mathsf{cov}} + 2C_{\mathsf{cov}}C_{\mathsf{sparse}} + 2$$

by (D.1) and (D.2), and

$$||z||_2 + ||z||_1/\sqrt{k} \leq 2C_{\mathsf{est}}\sqrt{\frac{k_n \log(p_n)}{n}}$$

by Assumption 3.3, and so

$$\left|w^\top(\widehat{\Sigma}_I - \Sigma_I)z\right| \leq ||\widehat{\Sigma}_I - \Sigma_I||_{\mathcal{S}_{k_n}} \cdot (C_{\mathsf{cov}} + 2C_{\mathsf{cov}}C_{\mathsf{sparse}} + 2) \cdot 2C_{\mathsf{est}}\sqrt{\frac{k_n \log(p_n)}{n}}$$

by applying Lemma 4.9. Noting that $||\widehat{\Sigma}_I - \Sigma_I||_{\mathcal{S}_{k_n}} \leq ||\widehat{\Sigma} - \Sigma||_{\mathcal{S}_{k_n}}$ trivially, the desired result of the lemma follows trivially from these bounds by setting $C_{\mathsf{submatrix}}$ appropriately. $\square$

### F.8 Proof of Lemma D.1

**Lemma D.1.** *Suppose that Assumptions 3.1, 3.2 and 3.4 hold. Let $g(X, X')$ and $g_1(X)$ be defined as in the proof of Theorem 4.1. Then*

$$\nu_{g_1}^2 = \mathsf{Var}(g_1(X)) \geq \frac{1}{\pi^2} C_{\mathsf{variance}}^2$$

*and*

$$\nu_{g_1}^3 \leq \eta_g^3 = \mathbb{E}\left[|g(X, X')|^3\right] \leq C_{\mathsf{moment}}$$

*where $C_{\mathsf{variance}}, C_{\mathsf{moment}}$ are constants depending only on $C_{\mathsf{cov}}, C_{\mathsf{kernel}}$ and not on $(n, p_n, k_n)$.*

*Proof of Lemma D.1.* First, we have

$$\begin{aligned}
g_1(X) &= \mathbb{E}\left[\mathrm{sign}(X - X')^\top \left(uv^\top \circ \cos\left(\frac{\pi}{2}T\right)\right)\mathrm{sign}(X - X') \mid X\right] \\
&= \mathbb{E}\left[\left(\mathrm{sign}(X - X') \otimes \mathrm{sign}(X - X')\right)^\top \mathsf{vec}\left(uv^\top \circ \cos\left(\frac{\pi}{2}T\right)\right) \mid X\right] \\
&= \mathbb{E}\left[\left(\mathrm{sign}(X - X') \otimes \mathrm{sign}(X - X')\right) \mid X\right]^\top \mathsf{vec}\left(uv^\top \circ \cos\left(\frac{\pi}{2}T\right)\right) \\
&= h_1(X)^\top \mathsf{vec}\left(uv^\top \circ \cos\left(\frac{\pi}{2}T\right)\right) \,,
\end{aligned}$$

where $h_1(X)$ is defined in Assumption 3.4, and has variance $\Sigma_{h_1}$. Therefore,

$$
\begin{aligned}
\nu_{g_1}^2 &= \mathsf{Var}(g_1(X)) \\
&= \mathsf{vec}\left(uv^\top \circ \cos\left(\frac{\pi}{2}T\right)\right)^\top \cdot \Sigma_{h_1} \cdot \mathsf{vec}\left(uv^\top \circ \cos\left(\frac{\pi}{2}T\right)\right) \\
&\geqslant C_{\mathsf{kernel}} \cdot \mathsf{vec}\left(uv^\top \circ \cos\left(\frac{\pi}{2}T\right)\right)^\top \cdot \Sigma_h \cdot \mathsf{vec}\left(uv^\top \circ \cos\left(\frac{\pi}{2}T\right)\right) \quad \text{(by Assumption 3.4)} \\
&= C_{\mathsf{kernel}} \cdot \mathsf{Var}\left(\mathsf{vec}\left(uv^\top \circ \cos\left(\frac{\pi}{2}T\right)\right)^\top h(X, X')\right) \\
&= C_{\mathsf{kernel}} \cdot \mathsf{Var}\left(\mathrm{sign}(X-X')^\top \left(uv^\top \circ \cos\left(\frac{\pi}{2}T\right)\right) \mathrm{sign}(X-X')\right) \\
&= C_{\mathsf{kernel}} \cdot \mathsf{Var}\left(\mathrm{sign}(Z)^\top \left(uv^\top \circ \cos\left(\frac{\pi}{2}T\right)\right) \mathrm{sign}(Z)\right) \\
&\geqslant C_{\mathsf{kernel}} \cdot C_{\mathsf{signs}} \cdot \left(uv^\top \circ \cos\left(\frac{\pi}{2}T\right)\right)_{ab}^2,
\end{aligned}
\tag{F.4}
$$

where for the next-to-last step, we take $Z \sim N(0, \Sigma)$ and apply Lemma 4.4 to see that $\mathrm{sign}(X - X')$ has the same distribution as $\mathrm{sign}(Z)$, and for the last step, we apply the following lemma (proved in Appendix F.10).

**Lemma F.3.** *Take any positive definite $\Sigma \in \mathbb{R}^{p \times p}$, any distinct $a, b \in \{1, \ldots, p\}$, and any matrix $M \in \mathbb{R}^{p \times p}$ with $M_{ja} = 0$ for all $j$. Let $Z \sim N(0, \Sigma)$. Then there exists a constant $C_{\mathsf{signs}}$ depending on $\mathsf{C}(\Sigma)$ only, such that*
$$\mathsf{Var}\left(\mathrm{sign}(Z)^\top M \, \mathrm{sign}(Z)\right) \geqslant C_{\mathsf{signs}} \cdot M_{ab}^2.$$

Applying this lemma with $M = uv^\top \circ \cos\left(\frac{\pi}{2}T\right)$ yields the lower bound (F.4) (since $v_a = 0$ so that $M_{ja} = 0$ for all $j$).

Finally, we have
$$\left(uv^\top \circ \cos\left(\frac{\pi}{2}T\right)\right)_{ab}^2 = u_a^2 v_b^2 \cos\left(\frac{\pi}{2}T_{ab}\right)^2 \geqslant (C_{\mathsf{cov}})^{-2},$$

where the last step holds because $u_a = v_b = 1$ and
$$\cos\left(\frac{\pi}{2}T_{ab}\right) = \sqrt{1 - \sin\left(\frac{\pi}{2}T_{ab}\right)^2} = \sqrt{1 - \Sigma_{ab}^2} \geqslant \lambda_{\min}(\Sigma_{ab,ab}) \geqslant (C_{\mathsf{cov}})^{-1}.$$

To summarize, we have
$$\nu_{g_1}^2 \geqslant \frac{C_{\mathsf{kernel}} C_{\mathsf{signs}}}{C_{\mathsf{cov}}^2} =: \frac{1}{\pi^2} C_{\mathsf{variance}}^2.$$

Next, we give an upper bound on $\nu_{g_1}^2$:

$$
\begin{aligned}
\nu_{g_1}^2 &= \mathsf{Var}(g_1(X)) \\
&= \mathsf{vec}\left(uv^\top \circ \cos\left(\frac{\pi}{2}T\right)\right)^\top \cdot \Sigma_{h_1} \cdot \mathsf{vec}\left(uv^\top \circ \cos\left(\frac{\pi}{2}T\right)\right) \\
&\leqslant \mathsf{vec}\left(uv^\top \circ \cos\left(\frac{\pi}{2}T\right)\right)^\top \cdot \Sigma_h \cdot \mathsf{vec}\left(uv^\top \circ \cos\left(\frac{\pi}{2}T\right)\right) \\
&= \mathsf{Var}\left(\mathrm{sign}(X-X')^\top \left(uv^\top \circ \cos\left(\frac{\pi}{2}T\right)\right) \mathrm{sign}(X-X')\right) \quad \text{(as for the lower bound earlier)} \\
&= \mathsf{Var}\left(g(X, X')\right) \leqslant \mathbb{E}\left[|g(X, X')|^2\right] \leqslant \mathbb{E}\left[|g(X, X')|^3\right]^{2/3} = \eta_g^2.
\end{aligned}
$$

Finally, we compute an upper bound on $\eta_g^3 = \mathbb{E}\left[|g(X, X')|^3\right]$. By Lemma E.1, there exists a decomposition
$$\cos\left(\frac{\pi}{2}T\right) = \sum_{r \geqslant 1} t_r a_r b_r^\top$$

where $t_r \geq 0$, $\sum_r t_r \leq 4$, and $||a_r||_\infty, ||b_r||_\infty \leq 1$. Note that, by (D.1),
$$||u||_2 = \sqrt{1 + ||\gamma_a||_2^2} \leq \sqrt{1 + C_{\text{cov}}^2}$$
and similarly $||v||_2 \leq \sqrt{1 + C_{\text{cov}}^2}$. Then for each $r$,
$$||u \circ a_r||_2 \vee ||v \circ b_r||_2 \leq \sqrt{1 + C_{\text{cov}}^2} \ .$$

Then we have

$$\mathbb{E}\left[|g(X, X')|^3\right] = \mathbb{E}\left[\left|\text{sign}(X-X')^\top \left(uv^\top \circ \cos\left(\frac{\pi}{2}T\right)\right)\text{sign}(X-X')\right|^3\right]$$

$$= \mathbb{E}\left[\left|\sum_{r \geq 1} t_r \cdot \text{sign}(X-X')^\top \left(uv^\top \circ a_r b_r^\top\right)\text{sign}(X-X')\right|^3\right]$$

$$\leq \sum_{r \geq 1} \frac{t_r}{4} \cdot \mathbb{E}\left[\left|4\,\text{sign}(X-X')^\top \left(uv^\top \circ a_r b_r^\top\right)\text{sign}(X-X')\right|^3\right] \quad \text{(by Jensen's inequality)}$$

$$\leq 64 \cdot \max_r \mathbb{E}\left[\left|\text{sign}(X-X')^\top \left(uv^\top \circ a_r b_r^\top\right)\text{sign}(X-X')\right|^3\right]$$

$$= 64 \cdot \max_r \mathbb{E}\left[\left|\text{sign}(X-X')^\top (u \circ a_r)\right|^3 \cdot \left|\text{sign}(X-X')^\top (v \circ b_r)\right|^3\right]$$

$$= 64 \cdot \max_r \sqrt{\mathbb{E}\left[\left|\text{sign}(X-X')^\top (u \circ a_r)\right|^6\right]} \cdot \sqrt{\mathbb{E}\left[\left|\text{sign}(X-X')^\top (v \circ b_r)\right|^6\right]}$$

$$\leq 64 ||u \circ a_r||_2^3 \cdot ||v \circ b_r||_2^3 \cdot \max_{||w||_2 = 1} \mathbb{E}\left[\left|\text{sign}(X-X')^\top w\right|^6\right]$$

$$\leq 64(1 + C_{\text{cov}}^2)^3 \cdot \max_{||w||_2 = 1} \mathbb{E}\left[\left|\text{sign}(X-X')^\top w\right|^6\right]$$

$$\leq 64(1 + C_{\text{cov}}^2)^3 \cdot C_{\text{cov}}^3 \cdot 6! \cdot 2\sqrt{e} =: C_{\text{moment}} \ ,$$

where the last inequality holds because $\text{sign}(X - X')$ is $C_{\text{cov}}$-subgaussian by Lemmas 4.4 and 4.5. $\square$

### F.9 Proofs of lemmas for the initial estimators

**Lemma E.1.** *There exist vectors $a_1, a_2, \ldots$ and $b_1, b_2, \ldots$ with $||a_r||_\infty, ||b_r||_\infty \leq 1$ for all $r \geq 1$, and a sequence $t_1, t_2, \cdots \geq 0$ with $\sum_r t_r = 4$, such that $\cos\left(\frac{\pi}{2}T\right) = \sum_{r \geq 1} t_r a_r b_r^\top$.*

*Proof of Lemma E.1.* We will use the matrix max norm, defined for a matrix $M \in \mathbb{R}^{d_1 \times d_2}$ as

$$||M||_{\max} =$$
$$\min\left\{\max_{1 \leq i \leq d_1} ||A_{(i)}||_2 \cdot \max_{1 \leq j \leq d_2} ||B_{(j)}||_2 \ : \ r \geq 1, A \in \mathbb{R}^{d_1 \times r}, B \in \mathbb{R}^{d_2 \times r} \text{ s.t. } M = A \cdot B^\top\right\} ,$$

where $A_{(i)}$ and $B_{(j)}$ denote the $i$th row of $A$ and the $j$th row of $B$, respectively. The matrix max norm satisfies several key properties that we will use here [Srebro and Shraibman, 2005]: first,

$$W \succeq 0 \ \Rightarrow \ ||W||_{\max} \leq \max_i W_{ii} \ ; \quad (\text{F.5})$$

second,

$$||W||_{\max} \leq 1 \ \Rightarrow \ \frac{W}{2} \in \text{ConvexHull}\left\{ab^\top : ||a||_\infty, ||b||_\infty \leq 1\right\} \ ; \quad (\text{F.6})$$

and finally,

$$||W \circ (uv^\top)||_* \leq ||W||_{\max} \text{ for all unit vectors } u, v \text{ and all matrices } W \ , \quad (\text{F.7})$$

where recall that $||\cdot||_*$ is the matrix nuclear norm (the sum of the singular values).

For our matrix $\cos\left(\frac{\pi}{2}T\right)$, Wegkamp and Zhao [2013] show that

$$\cos\left(\frac{\pi}{2}T\right) = \sum_{r\geq 0} \binom{1/2}{r}(-1)^r \Sigma \circ_{2r} \Sigma, \text{ and } \Sigma \circ_{2r} \Sigma \succeq 0 \text{ for all } r,$$

where $\Sigma \circ_{2r} \Sigma$ is the matrix with entries given by elementwise powers of $\Sigma$, that is, $(\Sigma \circ_{2r} \Sigma)_{jk} = (\Sigma_{jk})^{2r}$. Then for each $r \geq 0$, applying (F.5),

$$||\Sigma \circ_{2r} \Sigma||_{\max} \leq \max_i (\Sigma \circ_{2r} \Sigma)_{ii} = \max_i (\Sigma_{ii})^{2r} = 1,$$

since $\Sigma$ is a correlation matrix. Then

$$\left\|\cos\left(\frac{\pi}{2}T\right)\right\|_{\max} = \left\|\sum_{r\geq 0}\binom{1/2}{r}(-1)^r \Sigma \circ_{2r} \Sigma\right\|_{\max}$$

$$\leq \sum_{r\geq 0}\left|\binom{1/2}{r}\right| \cdot ||\Sigma \circ_{2r} \Sigma||_{\max} \leq \sum_{r\geq 0}\left|\binom{1/2}{r}\right| = 2,$$

where the last identity comes from Wegkamp and Zhao [2013]. Finally, by (F.6), we have

$$\frac{\cos\left(\frac{\pi}{2}T\right)}{4} \in \mathsf{ConvexHull}\left\{ab^\top : ||a||_\infty, ||b||_\infty \leq 1\right\}$$

and so $\cos\left(\frac{\pi}{2}T\right)$ can be expressed as a convex combination as stated in the lemma. $\square$

**Lemma E.2.** *For fixed $u, v$ with $||u||_2, ||v||_2 \leq 1$, for any $|t| \leq \frac{n}{4(1+\sqrt{5})C_{\mathsf{cov}}}$,*

$$\mathbb{E}\left[\exp\left(t \cdot u^\top (\widehat{T} - T)v\right)\right] \leq \exp\left(\frac{[4(1+\sqrt{5})]^2 t^2 \cdot C_{\mathsf{cov}}^2}{n}\right).$$

*Proof of Lemma E.2.* We start with a simple observation that

$$u^T\left(\widehat{T} - T\right)v = \frac{1}{4}(u+v)^\top(\widehat{T}-T)(u+v) - \frac{1}{4}(u-v)^\top(\widehat{T}-T)(u-v),$$

which gives us (via Cauchy-Schwartz)

$$\mathbb{E}\left[\exp\left(t \cdot u^\top(\widehat{T}-T)v\right)\right]$$
$$= \mathbb{E}\left[\exp\left(t \cdot \frac{1}{4}(u+v)^\top(\widehat{T}-T)(u+v) - t \cdot \frac{1}{4}(u-v)^\top(\widehat{T}-T)(u-v)\right)\right]$$
$$\leq \sqrt{\mathbb{E}\left[\exp\left(t \cdot \frac{1}{2}(u+v)^\top(\widehat{T}-T)(u+v)\right)\right]} \cdot \sqrt{\mathbb{E}\left[\exp\left(-t \cdot \frac{1}{2}(u-v)^\top(\widehat{T}-T)(u-v)\right)\right]}.$$

Note that $||\frac{1}{2}(u+v)||_2 \vee ||\frac{1}{2}(u-v)||_2 \leq 1$. Therefore, it will be sufficient to show that for any $|t| \leq \frac{n}{4(1+\sqrt{5})C_{\mathsf{cov}}}$ and any unit vector $w$,

$$\mathbb{E}\left[\exp\left(2t \cdot w^\top(\widehat{T}-T)w\right)\right] \leq \exp\left(\frac{[4(1+\sqrt{5})]^2 t^2 \cdot C_{\mathsf{cov}}^2}{n}\right). \qquad \text{(F.8)}$$

We will prove (F.8) using the Chernoff bounding technique. To that end, denote $S_n$ the group of permutations of $[n]$, and for any $i$, let $X_{(i)}$ denote the $i$-th row of $X$. For a fixed $w \in \mathbb{R}^p$, for each $i \in [n/2]$ and $\sigma \in S_n$, define

$$Z_{\sigma,i} = w^\top \left(\mathsf{sign}\left((X_{(\sigma(i))} - X_{(\sigma(i+n/2))})(X_{(\sigma(i))} - X_{(\sigma(i+n/2))})^\top\right) - T\right)w.$$

Observe that
$$w^T \left(\widehat{T} - T\right) w = \frac{1}{n!} \sum_{\sigma \in S_n} \frac{2}{n} \sum_{i \in [n/2]} Z_{\sigma,i}, \tag{F.9}$$

and that for any fixed $\sigma \in S_n$, the $Z_{\sigma,i}$'s are i.i.d. for $i = 1, \ldots, n/2$, and are identically distributed as

$$\widetilde{Z} = w^\top \left(\text{sign}\left((X_{(1)} - X_{(1+n/2)})(X_{(1)} - X_{(1+n/2)})^\top\right) - T\right) w .$$

Using Lemma 4.4 and Lemma 4.5, for any fixed unit vector $w \in \mathbb{R}^p$, $w^T \text{sign}\left(X_{(i)} - X_{(i+n/2)}\right)$ is a $C_{\text{cov}}$-subgaussian random variable, and

$$\widetilde{Z} = \left(w^T \text{sign}\left(X_{(1)} - X_{(1+n/2)}\right)\right)^2 - \mathbb{E}\left[\left(w^T \text{sign}\left(X_{(1)} - X_{(1+n/2)}\right)\right)^2\right] .$$

Applying Lemma F.4 (stated below), for any $|t| \leq \frac{1}{2(1+\sqrt{5})C_{\text{cov}}}$,

$$\mathbb{E}\left[\exp\left(t\widetilde{Z}\right)\right] \leq \exp\left(\frac{32t^2 C_{\text{cov}}^2}{1 - 4tC_{\text{cov}}}\right) \leq \exp\left(\frac{32t^2 C_{\text{cov}}^2}{1 - \frac{2}{1+\sqrt{5}}}\right) = \exp\left(8(1+\sqrt{5})^2 t^2 C_{\text{cov}}^2\right) .$$

Then, referring back to (F.9), for $0 < t \leq \frac{n}{4(1+\sqrt{5})C_{\text{cov}}}$,

$$\mathbb{E}\left[\exp\left(tw^T\left(\widehat{T} - T\right)w\right)\right] = \mathbb{E}\left[\exp\left(\frac{t}{n!}\sum_{\sigma \in S_n} \frac{2}{n}\sum_{i \in [n/2]} Z_{\sigma,i}\right)\right]$$
$$\leq \frac{1}{n!} \sum_{\sigma \in S_n} \mathbb{E}\left[\exp\left(\frac{2t}{n}\sum_{i \in [n/2]} Z_{\sigma,i}\right)\right] \quad \text{(by Jensen's inequality)}$$
$$= \frac{1}{n!} \sum_{\sigma \in S_n} \left(\mathbb{E}\left[\exp\left(\frac{2t}{n}\widetilde{Z}\right)\right]\right)^{n/2},$$

since for any fixed $\sigma$, the $Z_{\sigma,i}$'s are i.i.d., and are each equal to $\widetilde{Z}$ in distribution. Continuing from this step, we have

$$\mathbb{E}\left[\exp\left(tw^T\left(\widehat{T} - T\right)w\right)\right] = \left(\mathbb{E}\left[\exp\left(\frac{2t}{n}\widetilde{Z}\right)\right]\right)^{n/2}$$
$$\leq \left(\exp\left(8(1+\sqrt{5})^2 (2t/n)^2 C_{\text{cov}}^2\right)\right)^{n/2}$$
$$= \exp\left(\frac{\left[4(1+\sqrt{5})\right]^2 t^2 \cdot C_{\text{cov}}^2}{n}\right) .$$

$\square$

**Lemma F.4.** *Suppose $Z$ is $C$-subgaussian, that is, $\mathbb{E}\left[\exp(tZ)\right] \leq e^{Ct^2/2}$ for all $t \in \mathbb{R}$. Then*

$$\mathbb{E}\left[\exp\left\{t(Z^2 - \mathbb{E}[Z^2])\right\}\right] \leq \exp\left(\frac{32t^2 C^2}{1 - 4|t|C}\right)$$

*for all $|t| < \frac{1}{4C}$.*

We remark that it is well known that the square of a subgaussian random variable satisfies subgaussian tails near to its mean [see, for example, Lemmas 5.5, 5.14, 5.15 in Vershynin, 2012], but here we obtain small explicit constants.

*Proof.* The first part of this proof follows the arguments in Vershynin [2012, Lemma 5.5]. First, we bound $\mathbb{E}[Z^{2k}]$ for all integers $k \geq 1$. We have

$$\mathbb{E}[Z^{2k}] = \frac{C^k}{(2k)^k}\mathbb{E}\left[\left(\sqrt{\frac{2k}{C}}\cdot Z\right)^{2k}\right] \leq \frac{C^k}{(2k)^k}\mathbb{E}\left[(2k)!\cdot \exp\left\{\sqrt{\frac{2k}{C}}\cdot Z\right\}\right]$$

$$\leq \frac{C^k}{(2k)^k}(2k)!\cdot \exp\left\{\left(\sqrt{\frac{2k}{C}}\right)^2 \cdot C/2\right\} = \frac{C^k(2k)!e^k}{(2k)^k}.$$

Then, for any $t > 0$,

$$\mathbb{E}\left[e^{tZ^2}\right] = 1 + t\mathbb{E}[Z^2] + \sum_{k\geq 2}\mathbb{E}\left[\frac{(tZ^2)^k}{k!}\right] \leq 1 + t\mathbb{E}[Z^2] + \sum_{k\geq 2}\frac{t^k C^k e^k}{(2k)^k}\cdot \frac{(2k)!}{k!}$$

$$\leq 1 + t\mathbb{E}[Z^2] + \sum_{k\geq 2}\frac{t^k C^k e^k}{(2k)^k}\cdot \frac{e\cdot (2k)^{2k+1/2}\cdot e^{-2k}}{\sqrt{2\pi}\cdot k^{k+1/2}\cdot e^{-k}} \quad \text{(by Stirling's approximation for } (2k)! \text{ and } k!\text{)}$$

$$= 1 + t\mathbb{E}[Z^2] + \frac{e}{\sqrt{\pi}}\sum_{k\geq 2}(2tC)^k$$

$$= 1 + t\mathbb{E}[Z^2] + \frac{e}{\sqrt{\pi}}\cdot \frac{4t^2 C^2}{1-2tC}$$

$$\leq 1 + t\mathbb{E}[Z^2] + \frac{8t^2 C^2}{1-2tC}, \tag{F.10}$$

as long as $2tC < 1$ (we also use the fact $\frac{e}{\sqrt{\pi}} \leq 2$ to simplify). Next, trivially, for any $k \geq 2$, $\left|\mathbb{E}\left[(Z^2 - \mathbb{E}[Z^2])^k\right]\right| \leq 2^k \mathbb{E}[Z^{2k}]$. Then we have, for $|t| < \frac{1}{4C}$,

$$\mathbb{E}\left[\exp\left(t(Z^2 - \mathbb{E}[Z^2])\right)\right]$$

$$= 1 + \sum_{k\geq 2}\frac{t^k}{k!}\mathbb{E}\left[(Z^2 - \mathbb{E}[Z^2])^k\right] \leq 1 + \sum_{k\geq 2}\frac{|t|^k}{k!}\left|\mathbb{E}\left[(Z^2 - \mathbb{E}[Z^2])^k\right]\right|$$

$$\leq 1 + \sum_{k\geq 2}\frac{2^k|t|^k}{k!}\mathbb{E}[Z^{2k}] = \left(\sum_{k\geq 0}\frac{2^k|t|^k}{k!}\mathbb{E}[Z^{2k}]\right) - 2|t|\mathbb{E}[Z^2]$$

$$= \mathbb{E}\left[\exp\left(2|t|Z^2\right)\right] - 2|t|\mathbb{E}[Z^2]$$

$$\leq 1 + 2|t|\mathbb{E}[Z^2] + \frac{32t^2 C^2}{1-4|t|C} - 2|t|\mathbb{E}[Z^2] \quad \text{(by (F.10) with } 2|t| \text{ in place of } t\text{)}$$

$$\leq \exp\left\{\frac{32t^2 C^2}{1-4|t|C}\right\}.$$

□

### F.10 Lower bounds on variance for signs of a Gaussian

**Lemma F.3.** *Take any positive definite $\Sigma \in \mathbb{R}^{p\times p}$, any distinct $a, b \in \{1, \ldots, p\}$, and any matrix $M \in \mathbb{R}^{p\times p}$ with $M_{ja} = 0$ for all $j$. Let $Z \sim N(0, \Sigma)$. Then there exists a constant $C_{\mathsf{signs}}$ depending on $\mathsf{C}(\Sigma)$ only, such that*

$$\mathsf{Var}\left(\mathrm{sign}(Z)^\top M\,\mathrm{sign}(Z)\right) \geq C_{\mathsf{signs}}\cdot M_{ab}^2.$$

*Proof of Lemma F.3.* Let $(-a)$ denote the set $[p]\backslash\{a\}$. By the law of total variance,

$$\mathsf{Var}\left(\mathrm{sign}(Z)^\top M\,\mathrm{sign}(Z)\right) \geq \mathbb{E}\left[\mathsf{Var}\left(\mathrm{sign}(Z)^\top M\,\mathrm{sign}(Z) \mid Z_{(-a)}\right)\right].$$

Let $M_{j,(-a)} \in \mathbb{R}^{p-1}$ denote the $j$th row of $M$ with its $a$th entry removed, written as a column vector. Then, recalling that $M_{ja} = 0$ for all $j$, we have

$$\mathsf{Var}\left(\mathrm{sign}(Z)^\top M \mathrm{sign}(Z) \mid Z_{(-a)}\right)$$
$$= \mathsf{Var}\left(\mathrm{sign}(Z_a) \cdot M_{a,(-a)}^\top \mathrm{sign}(Z_{(-a)}) + \sum_{j \neq a} \mathrm{sign}(Z_j) \cdot M_{j,(-a)}^\top \mathrm{sign}(Z_{(-a)}) \mid Z_{(-a)}\right)$$
$$= \mathsf{Var}\left(\mathrm{sign}(Z_a) \cdot M_{a,(-a)}^\top \mathrm{sign}(Z_{(-a)}) \mid Z_{(-a)}\right)$$
$$= \mathsf{Var}\left(\mathrm{sign}(Z_a) \mid Z_{(-a)}\right) \cdot \left(M_{a,(-a)}^\top \mathrm{sign}(Z_{(-a)})\right)^2$$
$$= \mathsf{Var}\left(\mathrm{sign}(Z_{(-a)}^\top \beta_a + N(0,\nu_a^2))\right) \cdot \left(M_{a,(-a)}^\top \mathrm{sign}(Z_{(-a)})\right)^2,$$

where the last step holds since the distribution of $Z_a$ conditional on $Z_{(-a)}$ is given by $Z_{(-a)}^\top \beta_a + N(0, \nu_a^2)$ where $\beta_a = \Sigma_{(-a)}^{-1} \Sigma_{(-a),a}$ and $\nu_a^2 = \Sigma_{aa} - \Sigma_{(-a),a}^\top \Sigma_{(-a)}^{-1} \Sigma_{(-a),a}$. Continuing from this step, we have

$$\mathsf{Var}\left(\mathrm{sign}(Z)^\top M \mathrm{sign}(Z) \mid Z_{(-a)}\right)$$
$$= \left(1 - \mathbb{E}\left[\mathrm{sign}(Z_{(-a)}^\top \beta_a + N(0,\nu_a^2))\right]^2\right) \cdot \left(M_{a,(-a)}^\top \mathrm{sign}(Z_{(-a)})\right)^2$$
$$= \left(1 - \psi\left(\frac{Z_{(-a)}^\top \beta_a}{\nu_a}\right)^2\right) \cdot \left(M_{a,(-a)}^\top \mathrm{sign}(Z_{(-a)})\right)^2,$$

where $\psi(z) = \Phi(z) - \Phi(-z) \in (-1, 1)$.

Now we will give a lower bound on the expectation of this quantity. First consider the term $\psi\left(\frac{Z_{(-a)}^\top \beta_a}{\nu_a}\right)$. Note that

$$\frac{Z_{(-a)}^\top \beta_a}{\nu_a} \sim N\left(0, \frac{\beta_a^\top \Sigma_{(-a)} \beta_a}{\nu_a^2}\right) = N\left(0, \frac{\Sigma_{(-a),a}^\top \Sigma_{(-a)}^{-1} \Sigma_{(-a),a}}{\Sigma_{aa} - \Sigma_{(-a),a}^\top \Sigma_{(-a)}^{-1} \Sigma_{(-a),a}}\right)$$

and this variance is bounded by $\mathsf{C}(\Sigma)$. Then, for any $c \in (0, 1)$,

$$\mathbb{P}\left\{\left|\psi\left(\frac{Z_{(-a)}^\top \beta_a}{\nu_a}\right)\right| \leq \psi\left(\sqrt{\mathsf{C}(\Sigma)} \cdot \Phi^{-1}(1 - c/2)\right)\right\} = \mathbb{P}\left\{\left|\frac{Z_{(-a)}^\top \beta_a}{\nu_a}\right| \leq \sqrt{\mathsf{C}(\Sigma)} \cdot \Phi^{-1}(1 - c/2)\right\}$$
$$\geq \mathbb{P}\left\{|N(0, \mathsf{C}(\Sigma))| \leq \sqrt{\mathsf{C}(\Sigma)} \cdot \Phi^{-1}(1 - c/2)\right\} = 1 - c. \quad \text{(F.11)}$$

Next, note that $M_{a,(-a)}^\top \mathrm{sign}(Z_{(-a)})$ is $\left(\|M_{a,(-a)}\|_2^2 \cdot \mathsf{C}(\Sigma)\right)$-subgaussian by Lemma 4.5, and

$$\mathbb{E}\left[\left(M_{a,(-a)}^\top \mathrm{sign}(Z_{(-a)})\right)^2\right] \geq \|M_{a,(-a)}\|_2^2 \cdot \lambda_{\min}(T),$$

where $T = \mathbb{E}\left[\mathrm{sign}(Z)\mathrm{sign}(Z)^\top\right]$ (recall that $\Sigma = \sin\left(\frac{\pi}{2} T\right)$). Furthermore, by Wegkamp and Zhao [2013, Section 4.3], we have

$$T = \frac{2}{\pi} \sum_{k \geq 1} g(k) \Sigma \circ_k \Sigma,$$

where $g(k) \geq 0$ are nonnegative scalars, $g(1) = 1$, and $\Sigma \circ_k \Sigma$ is the $k$-fold Hadamard product, that is, $(\Sigma \circ_k \Sigma)_{ij} = (\Sigma_{ij})^k$. Wegkamp and Zhao [2013, Section 4.3] show also that $\Sigma \circ_k \Sigma \succeq 0$ for all $k$. Therefore,

$$T = \frac{2}{\pi} \Sigma + \frac{2}{\pi} \sum_{k \geq 2} g(k) \Sigma \circ_k \Sigma \succeq \frac{2}{\pi} \Sigma,$$

and so $\lambda_{\min}(T) \geq \frac{2}{\pi}\lambda_{\min}(\Sigma) \geq \frac{2}{\pi}(\mathsf{C}(\Sigma))^{-1}$. Applying Lemma F.5 (stated below),

$$\mathbb{P}\left\{\left(M_{a,(-a)}^\top \operatorname{sign}(Z_{(-a)})\right)^2 \geq \|M_{a,(-a)}\|_2^2 \cdot \lambda_{\min}(T)/2\right\} \geq \frac{1}{16e^{2\mathsf{C}(\Sigma)/\lambda_{\min}(T)}},$$

and so,

$$\mathbb{P}\left\{\left(M_{a,(-a)}^\top \operatorname{sign}(Z_{(-a)})\right)^2 \geq \|M_{a,(-a)}\|_2^2 \cdot \frac{1}{\pi\mathsf{C}(\Sigma)}\right\} \geq \frac{1}{16e^{\pi\mathsf{C}(\Sigma)^2}}.$$

Now set $c = \frac{1}{32e^{\pi\mathsf{C}(\Sigma)^2}}$ in (F.11). Then, we see that with probability at least $\frac{1}{32e^{\pi\mathsf{C}(\Sigma)^2}}$,

$$\left(1 - \psi\left(\frac{Z_{(-a)}^\top \beta_a}{\nu_a}\right)^2\right) \cdot \left(M_{a,(-a)}^\top \operatorname{sign}(Z_{(-a)})\right)^2 \geq$$

$$\left(1 - \psi\left(\sqrt{\mathsf{C}(\Sigma)} \cdot \Phi^{-1}\left(1 - \frac{1}{64e^{\pi\mathsf{C}(\Sigma)^2}}\right)\right)^2\right) \cdot \|M_{a,(-a)}\|_2^2 \cdot \frac{1}{\pi\mathsf{C}(\Sigma)}.$$

Therefore, combining everything,

$$\mathsf{Var}\left(\operatorname{sign}(Z)^\top M \operatorname{sign}(Z)\right) \geq \frac{1}{32e^{\pi\mathsf{C}(\Sigma)^2}} \cdot \left(1 - \psi\left(\sqrt{\mathsf{C}(\Sigma)} \cdot \Phi^{-1}\left(1 - \frac{1}{64e^{\pi\mathsf{C}(\Sigma)^2}}\right)\right)^2\right) \cdot \frac{\|M_{a,(-a)}\|_2^2}{\pi\mathsf{C}(\Sigma)}.$$

Noting that $\|M_{a,(-a)}\|_2^2 \geq M_{ab}^2$, this proves the desired result, where we define

$$C_{\mathsf{signs}} = \frac{1}{32e^{\pi\mathsf{C}(\Sigma)^2}} \cdot \left(1 - \psi\left(\sqrt{\mathsf{C}(\Sigma)} \cdot \Phi^{-1}\left(1 - \frac{1}{64e^{\pi\mathsf{C}(\Sigma)^2}}\right)\right)^2\right) \cdot \frac{1}{\pi\mathsf{C}(\Sigma)}.$$

$\square$

**Lemma F.5.** *Suppose that $W \in \mathbb{R}$ is a random variable with $\mathbb{E}[W] = 0$, $\mathbb{E}[W^2] \geq C_0$, and $\mathbb{E}[e^{tW}] \leq e^{C_1 t^2/2}$ for all $t \in \mathbb{R}$. Then*

$$\mathbb{P}\left\{W^2 \geq C_0/2\right\} \geq \frac{1}{16e^{2C_1/C_0}}.$$

*Proof of Lemma F.5.*

$$\begin{aligned}
C_0/2 &\leq \mathbb{E}[W^2] - C_0/2 \\
&= \mathbb{E}[W^2 \cdot \mathbb{1}\{W^2 \geq C_0/2\}] + \mathbb{E}[W^2 \cdot \mathbb{1}\{W^2 < C_0/2\}] - C_0/2 \\
&\leq \mathbb{E}[W^2 \cdot \mathbb{1}\{W^2 \geq C_0/2\}] \\
&\leq C_0 \mathbb{E}[(e^{W/\sqrt{C_0}} + e^{-W/\sqrt{C_0}}) \cdot \mathbb{1}\{W^2 \geq C_0/2\}] \quad (\text{since } t^2 \leq e^t + e^{-t} \text{ for all } t \in \mathbb{R}) \\
&= C_0 \mathbb{E}[e^{W/\sqrt{C_0}} \cdot \mathbb{1}\{W^2 \geq C_0/2\}] + C_0 \mathbb{E}[e^{-W/\sqrt{C_0}} \cdot \mathbb{1}\{W^2 \geq C_0/2\}] \\
&\leq C_0 \sqrt{\mathbb{E}[(e^{W/\sqrt{C_0}})^2] \cdot \mathbb{E}[\mathbb{1}\{W^2 \geq C_0/2\}^2]} + C_0 \sqrt{\mathbb{E}[(e^{-W/\sqrt{C_0}})^2] \cdot \mathbb{E}[\mathbb{1}\{W^2 \geq C_0/2\}^2]} \\
&= C_0 \sqrt{\mathbb{E}[e^{2W/\sqrt{C_0}}] \cdot \mathbb{P}\{W^2 \geq C_0/2\}} + C_0 \sqrt{\mathbb{E}[e^{-2W/\sqrt{C_0}}] \cdot \mathbb{P}\{W^2 \geq C_0/2\}} \\
&\leq C_0 \sqrt{e^{C_1/C_0 \cdot 2^2/2} \cdot \mathbb{P}\{W^2 \geq C_0/2\}} + C_0 \sqrt{e^{C_1/C_0 \cdot 2^2/2} \cdot \mathbb{P}\{W^2 \geq C_0/2\}},
\end{aligned}$$

and rearranging terms we have proved the lemma. $\square$

## F.11 Bounding the error in estimating the variance (Lemma D.5)

**Lemma D.5.** *Under the assumptions and definitions of Theorem 4.2, with probability at least $1 - \frac{1}{6p_n}$, if $n \geq k_n^2 \log(p_n)$, on the event that the bounds (3.1) in Assumption 3.3 hold,*

$$\left| \breve{S}_{ab} \cdot \det(\breve{\Theta}) - S_{ab} \cdot \det(\Theta) \right| \leq C_{\text{oracle}} \cdot \sqrt{\frac{k_n^2 \log(p_n)}{n}} \ .$$

*Proof of Lemma D.5.* Recall from the proof of Theorem 4.1 that we have defined

$$g(X, X') = \text{sign}(X - X')^\top \left( uv^\top \circ \cos\left(\frac{\pi}{2} T\right) \right) \text{sign}(X - X') \ ,$$

and $g_1(X) = \mathbb{E}\left[g(X, X') \mid X\right]$, where

$$u_a = 1, u_b = 0, u_I = -\gamma_a \text{ and } v_a = 0, v_b = 1, v_I = -\gamma_b \ .$$

Recall from the proof of Lemma D.1, given in Appendix F.8, that we have

$$\nu_{g_1}^2 = \text{Var}(g_1(X)) = \text{vec}\left( uv^\top \circ \cos\left(\frac{\pi}{2} T\right) \right)^\top \cdot \Sigma_{h_1} \cdot \text{vec}\left( uv^\top \circ \cos\left(\frac{\pi}{2} T\right) \right) \ ,$$

where $\Sigma_{h_1} = \text{Var}(h_1(X))$ for

$$h_1(X) = \mathbb{E}\left[ \text{sign}(X - X') \otimes \text{sign}(X - X') \mid X \right] \in \mathbb{R}^{p_n^2} \ .$$

Note that $\|\Sigma_{h_1}\|_\infty = 1$.

To estimate this variance, define vectors $\breve{u}$ and $\breve{v}$ with entries

$$\breve{u}_a = 1, \breve{u}_b = 0, \breve{u}_I = -\breve{\gamma}_a \text{ and } \breve{v}_a = 0, \breve{v}_b = 1, \breve{v}_I = -\breve{\gamma}_b \ ,$$

and define

$$\widehat{\Sigma}_{h_1} = \frac{1}{n} \sum_i \left( \widehat{h}_1(X_i) - \frac{1}{n} \sum_{i'} \widehat{h}_1(X_{i'}) \right) \left( \widehat{h}_1(X_i) - \frac{1}{n} \sum_{i'} \widehat{h}_1(X_{i'}) \right)^\top \ ,$$

where abusing notation, we write

$$\widehat{h}_1(X_i) = \frac{1}{n-1} \sum_{i' \neq i} h(X_i, X_i') = \frac{1}{n-1} \sum_{i' \neq i} \text{sign}(X_i - X_{i'}) \otimes \text{sign}(X_i - X_{i'}) \ .$$

We then define

$$\breve{\nu}_{g_1}^2 = \text{vec}\left( \breve{u}\breve{v}^\top \circ \cos\left(\frac{\pi}{2} \widehat{T}\right) \right)^\top \cdot \widehat{\Sigma}_{h_1} \cdot \text{vec}\left( \breve{u}\breve{v}^\top \circ \cos\left(\frac{\pi}{2} \widehat{T}\right) \right) \ .$$

Writing

$$x = \text{vec}\left( uv^\top \circ \cos\left(\frac{\pi}{2} T\right) \right) \text{ and } \breve{x} = \text{vec}\left( \breve{u}\breve{v}^\top \circ \cos\left(\frac{\pi}{2} \widehat{T}\right) \right) \ ,$$

we then have

$$\nu_{g_1}^2 = x^\top \Sigma_{h_1} x \text{ and } \breve{\nu}_{g_1}^2 = \breve{x}^\top \widehat{\Sigma}_{h_1} \breve{x} \ .$$

Define also

$$\overline{x} = \text{vec}\left( \breve{u}\breve{v}^\top \circ \cos\left(\frac{\pi}{2} T\right) \right) \ .$$

The following lemma, proved in Appendix F.12, carries out some elementary calculations on the norms of these vectors $x, \overline{x}, \breve{x}$.

**Lemma F.6.** *Define $x, \overline{x}, \breve{x}$ as in the proof of Lemma D.5, and assume $n \geqslant k_n^2 \log(p_n)$. If the bounds (3.1) in Assumption 3.3 hold then for constants $C_0, C_1, C_2, C_3$ that depend only on $C_{\mathsf{cov}}, C_{\mathsf{sparse}}, C_{\mathsf{est}}$,*

$$||x||_1 \leqslant C_0 k_n,$$

*and with probability at least $1 - \frac{1}{36 p_n}$, the following bounds all hold as well:*

$$||\breve{x} - \overline{x}||_1 \leqslant C_1 \sqrt{\frac{k_n^2 \log(p_n)}{n}},$$

$$||\breve{x} - x||_1 \leqslant C_2 \sqrt{\frac{k_n^3 \log(p_n)}{n}},$$

$$||\mathsf{mat}(\overline{x} - x)||_{\ell_1/\ell_2} \leqslant C_3 \sqrt{\frac{k_n^2 \log(p_n)}{n}},$$

*where $\mathsf{mat}(\cdot)$ reshapes a vector in $\mathbb{R}^{p_n^2}$ into a $p_n \times p_n$ matrix, and where we define the matrix $\ell_1/\ell_2$ norm as $M_{\ell_1/\ell_2} := \sum_j ||M_j||_2$, where $M_j$ is the $j$th column of $M$.*

We now continue bounding error in estimating $\nu_{g_1}$. We have:

$$\left| \breve{\nu}_{g_1}^2 - \nu_{g_1}^2 \right| = \left| \breve{x}^\top \widehat{\Sigma}_{h_1} \breve{x} - x^\top \Sigma_{h_1} x \right| \leqslant \left| x^\top (\widehat{\Sigma}_{h_1} - \Sigma_{h_1}) x \right| + \left| \breve{x}^\top \widehat{\Sigma}_{h_1} \breve{x} - x^\top \widehat{\Sigma}_{h_1} x \right|. \tag{F.12}$$

We bound each term separately. For the first term in (F.12), we apply the following lemma (proved in Appendix F.12 below):

**Lemma F.7.** *Under the same assumptions and notation as Lemmas D.1 and D.5, for a universal constant $C_{\mathsf{studentized}}$,*

$$\mathbb{P}\left\{ \left| x^\top (\widehat{\Sigma}_{h_1} - \Sigma_{h_1}) x \right| \leqslant C_{\mathsf{studentized}} \sqrt{\frac{k_n^2 \log(p_n)}{n}} \right\} \geqslant 1 - \frac{1}{36 p_n}.$$

For the second term in (F.12), since $\widehat{\Sigma}_{h_1} \geq 0$ and so $y \mapsto \sqrt{y^\top \widehat{\Sigma}_{h_1} y}$ is a norm and must satisfy the triangle inequality,

$$\left| \breve{x}^\top \widehat{\Sigma}_{h_1} \breve{x} - x^\top \widehat{\Sigma}_{h_1} x \right| = \left| \sqrt{\breve{x}^\top \widehat{\Sigma}_{h_1} \breve{x}} - \sqrt{x^\top \widehat{\Sigma}_{h_1} x} \right| \cdot \left| \sqrt{\breve{x}^\top \widehat{\Sigma}_{h_1} \breve{x}} + \sqrt{x^\top \widehat{\Sigma}_{h_1} x} \right|$$

$$\leqslant \left| \sqrt{\breve{x}^\top \widehat{\Sigma}_{h_1} \breve{x}} - \sqrt{x^\top \widehat{\Sigma}_{h_1} x} \right|^2 + \left| \sqrt{\breve{x}^\top \widehat{\Sigma}_{h_1} \breve{x}} - \sqrt{x^\top \widehat{\Sigma}_{h_1} x} \right| \cdot 2 \sqrt{x^\top \widehat{\Sigma}_{h_1} x}$$

$$\leqslant \left| \sqrt{\breve{x}^\top \widehat{\Sigma}_{h_1} \breve{x}} - \sqrt{x^\top \widehat{\Sigma}_{h_1} x} \right|^2 + \left| \sqrt{\breve{x}^\top \widehat{\Sigma}_{h_1} \breve{x}} - \sqrt{x^\top \widehat{\Sigma}_{h_1} x} \right| \cdot 2 \sqrt{x^\top \Sigma_{h_1} x + \left| x^\top (\widehat{\Sigma}_{h_1} - \Sigma_{h_1}) x \right|}. \tag{F.13}$$

To bound the difference term $\left| \sqrt{\breve{x}^\top \widehat{\Sigma}_{h_1} \breve{x}} - \sqrt{x^\top \widehat{\Sigma}_{h_1} x} \right|$ which appears twice in the expression above,

applying the triangle inequality several times, we have

$$\left|\sqrt{\breve{x}^\top \widehat{\Sigma}_{h_1} \breve{x}} - \sqrt{x^\top \widehat{\Sigma}_{h_1} x}\right| \leq \sqrt{(\breve{x}-x)^\top \widehat{\Sigma}_{h_1}(\breve{x}-x)} \tag{F.14}$$

$$\leq \sqrt{(\breve{x}-x)^\top \Sigma_{h_1}(\breve{x}-x)} + \sqrt{\left|(\breve{x}-x)^\top (\widehat{\Sigma}_{h_1}-\Sigma_{h_1})(\breve{x}-x)\right|} \tag{F.15}$$

$$\leq \sqrt{(\breve{x}-x)^\top \Sigma_{h_1}(\breve{x}-x)} + \sqrt{||\widehat{\Sigma}_{h_1}-\Sigma_{h_1}||_\infty \cdot ||\breve{x}-x||_1^2} \tag{F.16}$$

$$\leq \sqrt{(\breve{x}-\overline{x})^\top \Sigma_{h_1}(\breve{x}-\overline{x})} + \sqrt{(\overline{x}-x)^\top \Sigma_{h_1}(\overline{x}-x)} + \sqrt{||\widehat{\Sigma}_{h_1}-\Sigma_{h_1}||_\infty \cdot ||\breve{x}-x||_1^2} \tag{F.17}$$

$$\leq \sqrt{||\Sigma_{h_1}||_\infty ||\breve{x}-\overline{x}||_1^2} + \sqrt{(\overline{x}-x)^\top \Sigma_{h_1}(\overline{x}-x)} + \sqrt{||\widehat{\Sigma}_{h_1}-\Sigma_{h_1}||_\infty \cdot ||\breve{x}-x||_1^2} \tag{F.18}$$

$$\leq C_1 \sqrt{\frac{k_n^2 \log(p_n)}{n}} + \sqrt{(\overline{x}-x)^\top \Sigma_{h_1}(\overline{x}-x)} + \sqrt{||\widehat{\Sigma}_{h_1}-\Sigma_{h_1}||_\infty} \cdot C_2 \sqrt{\frac{k_n^3 \log(p_n)}{n}}. \tag{F.19}$$

Next, we state two lemmas, which are proved in Appendix F.12.

**Lemma F.8.** *With probability at least $1 - \frac{1}{9p_n}$,*

$$||\widehat{\Sigma}_{h_1} - \Sigma_{h_1}||_\infty \leq 100\sqrt{\frac{\log(p_n)}{n}}.$$

**Lemma F.9.** *Let $\Sigma_{h_1}$ be defined as in Assumption 3.4. For every $z \in \mathbb{R}^{p_n^2}$,*

$$z^\top \Sigma_{h_1} z \leq \lambda_{\max}(\Sigma) \cdot ||\mathsf{mat}(z)||_{\ell_1/\ell_2}^2,$$

*where $||\mathsf{mat}(z)||_{\ell_1/\ell_2}$ is defined as in the statement of Lemma F.6.*

From this point on, we assume that the bounds derived in Lemmas F.6 and F.8 all hold (which the lemmas have shown to be true with probability at least $1 - \frac{1}{6p_n}$, on the event that the bounds (3.1) of Assumption 3.3 hold.) By Lemmas F.6 and F.9,

$$(\overline{x}-x)^\top \Sigma_{h_1}(\overline{x}-x) \leq C_{\mathsf{cov}} \cdot \left(C_3 \sqrt{\frac{k_n^2 \log(p_n)}{n}}\right)^2.$$

Applying this bound, along with the high probability events of Lemmas F.7 and F.8, we return to (F.14) and obtain

$$\left|\sqrt{\breve{x}^\top \widehat{\Sigma}_{h_1} \breve{x}} - \sqrt{x^\top \widehat{\Sigma}_{h_1} x}\right| \leq$$

$$C_1 \sqrt{\frac{k_n^2 \log(p_n)}{n}} + \sqrt{C_{\mathsf{cov}} \cdot \left(C_3 \sqrt{\frac{k_n^2 \log(p_n)}{n}}\right)^2} + \sqrt{100\sqrt{\frac{\log(p_n)}{n}} \cdot C_2 \sqrt{\frac{k_n^3 \log(p_n)}{n}}}$$

$$= \sqrt{\frac{k_n^2 \log(p_n)}{n}} \cdot \left(C_1 + C_3\sqrt{C_{\mathsf{cov}}} + 10C_2 \sqrt[4]{\frac{k_n^2 \log(p_n)}{n}}\right) \leq \sqrt{\frac{k_n^2 \log(p_n)}{n}} \cdot C_4,$$

where for the last step we define $C_4 = C_1 + C_3\sqrt{C_{\mathsf{cov}}} + 10C_2$ and use the assumption $n \geq k_n^2 \log(p_n)$.

Next, returning to (F.13),

$$\left|\breve{x}^\top\widehat{\Sigma}_{h_1}\breve{x} - x^\top\widehat{\Sigma}_{h_1}x\right|$$
$$\leqslant \left|\sqrt{\breve{x}^\top\widehat{\Sigma}_{h_1}\breve{x}} - \sqrt{x^\top\widehat{\Sigma}_{h_1}x}\right|^2 + \left|\sqrt{\breve{x}^\top\widehat{\Sigma}_{h_1}\breve{x}} - \sqrt{x^\top\widehat{\Sigma}_{h_1}x}\right| \cdot 2\sqrt{x^\top\Sigma_{h_1}x + \left|x^\top(\widehat{\Sigma}_{h_1} - \Sigma_{h_1})x\right|}$$
$$\leqslant C_4^2 \cdot \frac{k_n^2\log(p_n)}{n} + C_4\sqrt{\frac{k_n^2\log(p_n)}{n}} \cdot 2\sqrt{C_{\mathsf{moment}}^{2/3} + C_{\mathsf{studentized}}\sqrt{\frac{k_n^2\log(p_n)}{n}}},$$

where the last step applies the high probability event of Lemma F.7, and uses the fact that $x^\top\Sigma_{h_1}x = \nu_{g_1}^2 \leqslant C_{\mathsf{moment}}^{2/3}$ by Lemma D.1. Defining $C_5 = C_4^2 + C_4 \cdot 2\sqrt{C_{\mathsf{moment}}^{2/3} + C_{\mathsf{studentized}}}$, and using the assumption $n \geqslant k_n^2\log(p_n)$, we have

$$\left|\breve{x}^\top\widehat{\Sigma}_{h_1}\breve{x} - x^\top\widehat{\Sigma}_{h_1}x\right| \leqslant C_5\sqrt{\frac{k_n^2\log(p_n)}{n}}.$$

Finally, returning to (F.12) and applying Lemma F.8, we see that

$$\left|\breve{\nu}_{g_1}^2 - \nu_{g_1}^2\right| \leqslant C_{\mathsf{studentized}} \cdot \sqrt{\frac{k_n^2\log(p_n)}{n}} + C_5\sqrt{\frac{k_n^2\log(p_n)}{n}}.$$

Next, we have

$$|\breve{\nu}_{g_1} - \nu_{g_1}| = \frac{|\breve{\nu}_{g_1}^2 - \nu_{g_1}^2|}{\breve{\nu}_{g_1} + \nu_{g_1}} \leqslant \frac{|\breve{\nu}_{g_1}^2 - \nu_{g_1}^2|}{\nu_{g_1}} \leqslant \frac{(C_{\mathsf{studentized}} + C_5)\sqrt{\frac{k_n^2\log(p_n)}{n}}}{\frac{1}{\pi}C_{\mathsf{variance}}},$$

where for the denominator we apply Lemma D.1. Finally, since we know that $S_{ab} = \pi\nu_{g_1} \cdot (\det(\Theta))^{-1}$ and $\breve{S}_{ab} = \pi\breve{\nu}_{g_1} \cdot (\det(\breve{\Theta}))^{-1}$, and then we have

$$\left|\breve{S}_{ab} \cdot \det(\breve{\Theta}) - S_{ab} \cdot \det(\Theta)\right| = \pi \cdot |\breve{\nu}_{g_1} - \nu_{g_1}| \leqslant \pi \cdot \frac{(C_{\mathsf{studentized}} + C_5)\sqrt{\frac{k_n^2\log(p_n)}{n}}}{\frac{1}{\pi}C_{\mathsf{variance}}}.$$

Defining

$$C_{\mathsf{oracle}} \geqslant \pi \cdot \frac{C_{\mathsf{studentized}} + C_5}{\frac{1}{\pi}C_{\mathsf{variance}}},$$

we see that

$$\left|\breve{S}_{ab} \cdot \det(\breve{\Theta}) - S_{ab} \cdot \det(\Theta)\right| \leqslant C_{\mathsf{oracle}} \cdot \sqrt{\frac{k_n^2\log(p_n)}{n}}.$$

□

### F.12 Calculations for the variance estimate (Lemma D.5)

**Lemma F.6.** *Define $x, \overline{x}, \breve{x}$ as in the proof of Lemma D.5, and assume $n \geqslant k_n^2\log(p_n)$. If the bounds (3.1) in Assumption 3.3 hold then for constants $C_0, C_1, C_2, C_3$ that depend only on $C_{\mathsf{cov}}, C_{\mathsf{sparse}}, C_{\mathsf{est}}$,*

$$||x||_1 \leqslant C_0 k_n,$$

and with probability at least $1 - \frac{1}{36p_n}$, the following bounds all hold as well:

$$||\breve{x} - \overline{x}||_1 \leqslant C_1 \sqrt{\frac{k_n^2 \log(p_n)}{n}},$$

$$||\breve{x} - x||_1 \leqslant C_2 \sqrt{\frac{k_n^3 \log(p_n)}{n}},$$

$$||\mathsf{mat}(\overline{x} - x)||_{\ell_1/\ell_2} \leqslant C_3 \sqrt{\frac{k_n^2 \log(p_n)}{n}},$$

where $\mathsf{mat}(\cdot)$ reshapes a vector in $\mathbb{R}^{p_n^2}$ into a $p_n \times p_n$ matrix, and where we define the matrix $\ell_1/\ell_2$ norm as $M_{\ell_1/\ell_2} := \sum_j ||M_j||_2$, where $M_j$ is the $j$th column of $M$.

*Proof of Lemma F.6.* We calculate

$$||x||_1 = \left\| uv^\top \circ \cos\left(\frac{\pi}{2}T\right) \right\|_1 \leqslant ||uv^\top||_1 \cdot \left\| \cos\left(\frac{\pi}{2}T\right) \right\|_\infty \leqslant ||u||_1 ||v||_1$$
$$\leqslant k_n(1 + 2C_{\mathsf{cov}}C_{\mathsf{sparse}})^2 =: C_0 k_n ,$$

where for the last inequality we apply (D.2).

Next,

$$||\breve{x} - \overline{x}||_1 = \left\| ||\breve{u}\breve{v}^\top \circ \left(\cos\left(\frac{\pi}{2}\widehat{T}\right) - \cos\left(\frac{\pi}{2}T\right)\right) \right\|_1$$
$$\leqslant (||u||_1 + ||\breve{u} - u||_1) \cdot (||v||_1 + ||\breve{v} - v||_1) \cdot \left\| \cos\left(\frac{\pi}{2}\widehat{T}\right) - \cos\left(\frac{\pi}{2}T\right) \right\|_\infty.$$

Applying (D.2) and Assumption 3.3, and the fact that $\cos(\cdot)$ is 1-Lipschitz, if the bounds in Assumption 3.3 hold, we then have

$$||\breve{x} - \overline{x}||_1 \leqslant \left( \sqrt{k_n}(1 + 2C_{\mathsf{cov}}C_{\mathsf{sparse}}) + C_{\mathsf{est}} \sqrt{\frac{k_n^2 \log(p_n)}{n}} \right)^2 \cdot \frac{\pi}{2} ||\widehat{T} - T||_\infty.$$

Furthermore, applying Lemma D.2, with probability at least $1 - \frac{1}{36p_n}$,

$$||\breve{x} - \overline{x}||_1 \leqslant \left( \sqrt{k_n}(1 + 2C_{\mathsf{cov}}C_{\mathsf{sparse}}) + C_{\mathsf{est}} \sqrt{\frac{k_n^2 \log(p_n)}{n}} \right)^2 \cdot \frac{\pi}{2} \sqrt{\frac{4 \log(36p_n^3)}{n}}.$$

Finally, since $4 \log(36p_n^3) \leqslant 36 \log(p_n) \leqslant 4n$, where the last step holds by assumption in Theorem 4.2,

$$||\breve{x} - \overline{x}||_1 \leqslant \sqrt{\frac{k_n^2 \log(p_n)}{n}} \cdot ((1 + 2C_{\mathsf{cov}}C_{\mathsf{sparse}}) + C_{\mathsf{est}})^2 \cdot \pi$$
$$= C_1 \sqrt{\frac{k_n^2 \log(p_n)}{n}} \text{ for } C_1 := ((1 + 2C_{\mathsf{cov}}C_{\mathsf{sparse}}) + C_{\mathsf{est}})^2 \cdot \pi .$$

Next,

$$||\breve{x} - x||_1 \leq ||\breve{x} - \overline{x}||_1 + ||\overline{x} - x||_1$$

$$\leq C_1 \sqrt{\frac{k_n^2 \log(p_n)}{n}} + ||\overline{x} - x||_1$$

$$= C_1 \sqrt{\frac{k_n^2 \log(p_n)}{n}} + \left\|(\breve{u}\breve{v}^\top - uv^\top) \circ \cos\left(\frac{\pi}{2}T\right)\right\|_1$$

$$\leq C_1 \sqrt{\frac{k_n^2 \log(p_n)}{n}} + \left\|\breve{u}(\breve{v} - v)^\top \circ \cos\left(\frac{\pi}{2}T\right)\right\|_1 + \left\|(\breve{u} - u)v^\top \circ \cos\left(\frac{\pi}{2}T\right)\right\|_1$$

$$\leq C_1 \sqrt{\frac{k_n^2 \log(p_n)}{n}} + ||\breve{u}||_1 ||\breve{v} - v||_1 \left\|\cos\left(\frac{\pi}{2}T\right)\right\|_\infty + ||\breve{u} - u||_1 ||v||_1 \left\|\cos\left(\frac{\pi}{2}T\right)\right\|_\infty$$

$$\leq C_1 \sqrt{\frac{k_n^2 \log(p_n)}{n}} + ||\breve{u}||_1 ||\breve{v} - v||_1 + ||\breve{u} - u||_1 ||v||_1$$

$$\leq C_1 \sqrt{\frac{k_n^2 \log(p_n)}{n}} + \left(\sqrt{k_n}(1 + 2C_{\mathsf{cov}}C_{\mathsf{sparse}}) + C_{\mathsf{est}}\sqrt{\frac{k_n^2 \log(p_n)}{n}}\right) \cdot C_{\mathsf{est}}\sqrt{\frac{k_n^2 \log(p_n)}{n}}$$

$$+ \sqrt{k_n}(1 + 2C_{\mathsf{cov}}C_{\mathsf{sparse}}) \cdot C_{\mathsf{est}}\sqrt{\frac{k_n^2 \log(p_n)}{n}}$$

$$\leq C_2 \sqrt{\frac{k_n^3 \log(p_n)}{n}},$$

where the next-to-last step applies (D.2) (assuming the bounds (3.1) in Assumption 3.3 hold), and where for the last step, $C_2 = C_1 + 2(1 + 2C_{\mathsf{cov}}C_{\mathsf{sparse}}) \cdot C_{\mathsf{est}} + C_{\mathsf{est}}^2$ and we use the assumption $n \geq k_n^2 \log(p)$.

Finally, noting that $\overline{x} - x = \mathsf{vec}\left((\breve{u}\breve{v}^\top - uv^\top) \circ \cos\left(\frac{\pi}{2}T\right)\right)$, we calculate the $\ell_1/\ell_2$ norm of this matrix:

$$||\mathsf{mat}(\overline{x} - x)||_{\ell_1/\ell_2} = \left\|(\breve{u}\breve{v}^\top - uv^\top) \circ \cos\left(\frac{\pi}{2}T\right)\right\|_{\ell_1/\ell_2}$$

$$= \sum_j \left\|\left[(\breve{u}\breve{v}^\top - uv^\top) \circ \cos\left(\frac{\pi}{2}T\right)\right]_j\right\|_2$$

$$\leq \sum_j \left\|[\breve{u}\breve{v}^\top - uv^\top]_j\right\|_2 \cdot \left\|\cos\left(\frac{\pi}{2}T\right)\right\|_\infty$$

$$\leq \sum_j \left\|[\breve{u}\breve{v}^\top - uv^\top]_j\right\|_2$$

$$\leq \sum_j ||\breve{u} \cdot (\breve{v}_j - v_j)||_2 + ||(\breve{u} - u) \cdot v_j||_2$$

$$= \sum_j ||\breve{u}||_2 \cdot |\breve{v}_j - v_j| + ||\breve{u} - u||_2 \cdot |v_j|$$

$$= ||\breve{u}||_2 ||\breve{v} - v||_1 + ||\breve{u} - u||_2 ||v||_1.$$

Next, applying (D.1), if the bounds (3.1) in Assumption 3.3 hold, we then have

$$||\mathsf{mat}(\overline{x} - x)||_{\ell_1/\ell_2}$$

$$\leq \left(\sqrt{1 + C_{\mathsf{cov}}^2} + C_{\mathsf{est}}\sqrt{\frac{k_n \log(p_n)}{n}}\right) \cdot C_{\mathsf{est}} \cdot \sqrt{\frac{k_n^2 \log(p_n)}{n}} + C_{\mathsf{est}} \cdot \sqrt{\frac{k_n \log(p_n)}{n}} \cdot \sqrt{k_n}\sqrt{1 + C_{\mathsf{cov}}}$$

$$\leq C_3 \sqrt{\frac{k_n^2 \log(p_n)}{n}},$$

where we define $C_3 = 2\sqrt{1 + C_{\mathsf{cov}}^2} \cdot C_{\mathsf{est}} + C_{\mathsf{est}}^2$ and use the assumption that $n \geq k_n^2 \log(p_n)$. □

**Lemma F.7.** *Under the same assumptions and notation as Lemmas D.1 and D.5, for a universal constant $C_{\text{studentized}}$,*

$$\mathbb{P}\left\{\left|x^\top(\widehat{\Sigma}_{h_1} - \Sigma_{h_1})x\right| \leq C_{\text{studentized}}\sqrt{\frac{k_n^2 \log(p_n)}{n}}\right\} \geq 1 - \frac{1}{36p_n}.$$

*Proof of Lemma F.7.* By definition, we have $\Sigma_{h_1} = \text{Var}(h_1(X))$ for

$$h_1(X) = \mathbb{E}\left[\text{sign}(X - X') \otimes \text{sign}(X - X') \mid X\right] \in \mathbb{R}^{p_n^2}.$$

Therefore, since $x$ is fixed,

$$x^\top \Sigma_{h_1} x = x^\top \text{Var}(h_1(X))x = \text{Var}(x^\top h_1(X)) = \text{Var}(g_1(X)) = \nu_{g_1}^2,$$

where we recall that $g_1(X) = \mathbb{E}\left[g(X, X') \mid X\right]$ where we define the kernel

$$g(X, X') = \text{sign}(X - X')^\top \left(uv^\top \circ \cos\left(\frac{\pi}{2}T\right)\right) \text{sign}(X - X') = x^\top h(X, X').$$

Define

$$\gamma(X, X', X'') = \frac{g(X, X')g(X, X'') + g(X', X)g(X', X'') + g(X'', X)g(X'', X')}{3}.$$

Note that $\gamma(X, X', X'')$ is a U-statistic of order 3, and that

$$\sup_{X, X'} |g(X, X')| \leq \|x\|_1 \sup_{X, X'} \|h(X, X')\|_\infty = \|x\|_1 \leq C_0 k_n,$$

where the last step applies Lemma F.6. So,

$$\|\gamma\|_\infty := \sup_{X, X', X''} |\gamma(X, X', X'')| \leq \sup_{X, X'} |g(X, X')|^2 \leq C_0^2 k_n^2.$$

And,

$$\text{Var}(\gamma) := \text{Var}(\gamma(X, X', X'')) \leq \mathbb{E}\left[\gamma(X, X', X'')^2\right] \leq \mathbb{E}\left[|g(X, X')|^4\right]$$
$$\leq \mathbb{E}\left[|g(X, X')|^3\right] \cdot C_0 k_n \leq C_{\text{moment}} \cdot C_0 k_n,$$

where we use Lemma D.1 for the last bound.

Next, we have

$$\mathbb{E}\left[\gamma(X, X', X'')\right] = \mathbb{E}\left[g(X, X')g(X, X'')\right] = \mathbb{E}\left[\mathbb{E}\left[g(X, X')g(X, X'')|X\right]\right]$$
$$= \mathbb{E}\left[\mathbb{E}\left[g(X, X')|X\right]\mathbb{E}\left[g(X, X'')|X\right]\right] = \mathbb{E}\left[g_1(X)^2\right].$$

Therefore,

$$x^\top \Sigma_{h_1} x = \nu_{g_1}^2 = \text{Var}(g_1(X)) = \mathbb{E}\left[\gamma(X, X', X'')\right] - \mathbb{E}\left[g_1(X)\right]^2 = \mathbb{E}\left[\gamma(X, X', X'')\right] - \mathbb{E}\left[g(X, X')\right]^2.$$

Next, examining the definition of $\widehat{\Sigma}_{h_1}$, we obtain

$$x^\top \widehat{\Sigma}_{h_1} x = \frac{1}{n(n-1)^2}\left[\sum_{i \neq i' \neq i''} \gamma(X_i, X_{i'}, X_{i''}) + \sum_{i \neq i'} g(X_i, X_{i'})^2\right] - \left(\frac{1}{\binom{n}{2}}\sum_{i < i'} g(X_i, X_{i'})\right)^2.$$

Therefore, using the fact that $|g(X, X')| \leq C_0 k_n$ always,

$$\left| x^\top \hat{\Sigma}_{h_1} x - x^\top \Sigma_{h_1} x \right| \leq \left| \frac{1}{\binom{n}{3}} \sum_{i<i'<i''} \gamma(X_i, X_{i'}, X_{i''}) - \mathbb{E}[\gamma(X, X', X'')] \right|$$
$$+ \frac{C_0^2 k_n^2}{n-1} + \left| \left( \frac{1}{\binom{n}{2}} \sum_{i<i'} g(X_i, X_{i'}) \right)^2 - \mathbb{E}[g(X, X')]^2 \right|.$$

Now, using Bernstein's inequality for U-statistics (Peel et al. [2010, Theorem 2]), for any $\delta > 0$,

$$\mathbb{P}\left\{ \left| \frac{1}{\binom{n}{3}} \sum_{i<i'<i''} \gamma(X_i, X_{i'}, X_{i''}) - \mathbb{E}[\gamma(X, X', X'')] \right| > \sqrt{\frac{2\mathsf{Var}(\gamma) \log(2/\delta)}{(n/3)}} + \frac{2||\gamma||_\infty \log(2/\delta)}{3(n/3)} \right\} \leq \delta.$$

Therefore, with probability at least $1 - \frac{1}{72 p_n}$,

$$\left| \frac{1}{\binom{n}{3}} \sum_{i<i'<i''} \gamma(X_i, X_{i'}, X_{i''}) - \mathbb{E}[\gamma(X, X', X'')] \right| \leq$$
$$\sqrt{\frac{2 C_{\mathsf{moment}} \cdot C_0 k_n \log(2 \cdot 72 p_n)}{(n/3)}} + \frac{2 C_0^2 k_n^2 \log(2 \cdot 72 p_n)}{3(n/3)} \leq \sqrt{\frac{k_n^2 \log(p_n)}{n}} \cdot C',$$

where

$$C' = \sqrt{6 C_{\mathsf{moment}} C_0 \log_2(144)} + 2 C_0^2 \log_2(144),$$

and we use the assumption $n \geq k_n^2 \log(p_n)$ and $p_n \geq 2$. And, again using Bernstein's inequality for U-statistics, and using the fact that $|g(X, X')| \leq C_0 k_n$ always, with probability at least $1 - \frac{1}{72 p_n}$,

$$\left| \frac{1}{\binom{n}{2}} \sum_{i<i'} g(X_i, X_{i'}) - \mathbb{E}[g(X, X')] \right| \leq$$
$$\sqrt{\frac{2 C_0^2 k_n^2 \log(2 \cdot 72 p_n)}{(n/2)}} + \frac{2 C_0 k_n \log(2 \cdot 72 p_n)}{3(n/2)} \leq \sqrt{\frac{k_n^2 \log(p_n)}{n}} \cdot C'',$$

where

$$C'' = \sqrt{\frac{2 C_0^2 \log_2(144)}{(1/2)}} + \frac{2 C_0 \log_2(144)}{3/2},$$

and we use the assumption $n \geq k_n^2 \log(p_n)$ and $p_n \geq 2$. Therefore,

$$\left| \left( \frac{1}{\binom{n}{2}} \sum_{i<i'} g(X_i, X_{i'}) \right)^2 - \mathbb{E}[g(X, X')]^2 \right| \leq \left| \frac{1}{\binom{n}{2}} \sum_{i<i'} g(X_i, X_{i'}) - \mathbb{E}[g(X, X')] \right|^2 +$$
$$\left| \frac{1}{\binom{n}{2}} \sum_{i<i'} g(X_i, X_{i'}) - \mathbb{E}[g(X, X')] \right| \cdot 2|\mathbb{E}[g(X, X')]| \leq C''' \sqrt{\frac{k_n^2 \log(p_n)}{n}},$$

where we set

$$C''' = C''^2 + 2 C'' \cdot C_{\mathsf{moment}}^{1/3}$$

and again use $n \geq k_n^2 \log(p_n)$, and apply Lemma D.1 to bound $|\mathbb{E}[g(X, X')]|$. Combining everything, this proves that, with probability at least $1 - \frac{1}{36 p_n}$,

$$\left| x^\top \hat{\Sigma}_{h_1} x - x^\top \Sigma_{h_1} x \right| \leq \sqrt{\frac{k_n^2 \log(p_n)}{n}} \cdot C' + \frac{C_0^2 k_n^2}{n-1} + C''' \sqrt{\frac{k_n^2 \log(p_n)}{n}}.$$

Setting
$$C_{\text{studentized}} = C' + C''' + 2C_0^2$$
and using the fact that $n \geq 2$ and $n \geq k_n^2 \log(p_n)$, we have
$$\left|x^\top \widehat{\Sigma}_{h_1} x - x^\top \Sigma_{h_1} x\right| \leq C_{\text{studentized}} \sqrt{\frac{k_n^2 \log(p_n)}{n}} \ .$$

$\square$

**Lemma F.8.** *With probability at least $1 - \frac{1}{9p_n}$,*
$$||\widehat{\Sigma}_{h_1} - \Sigma_{h_1}||_\infty \leq 100 \sqrt{\frac{\log(p_n)}{n}} \ .$$

*Proof of Lemma F.8.* From our definitions, we see that
$$\Sigma_{h_1} = \mathsf{Var}(h_1(X)) = \mathbb{E}\left[h_1(X)h_1(X)^\top\right] - \mathbb{E}[h_1(X)]\mathbb{E}[h_1(X)]^\top \ ,$$
and
$$\widehat{\Sigma}_{h_1} = \frac{1}{n}\sum_i \widehat{h}_1(X_i)\widehat{h}_1(X_i)^\top - \left(\frac{1}{n}\sum_i \widehat{h}_1(X_i)\right)\left(\frac{1}{n}\sum_i \widehat{h}_1(X_i)\right)^\top \ .$$

First, we bound $||\frac{1}{n}\sum_i \widehat{h}_1(X_i)\widehat{h}_1(X_i)^\top - \mathbb{E}\left[h_1(X)h_1(X)^\top\right]||_\infty$. We have

$$\left\|\frac{1}{n}\sum_i \widehat{h}_1(X_i)\widehat{h}_1(X_i)^\top - \mathbb{E}\left[h_1(X)h_1(X)^\top\right]\right\|_\infty \leq$$
$$\left\|\frac{1}{n}\sum_i \widehat{h}_1(X_i)\widehat{h}_1(X_i)^\top - \frac{1}{n}\sum_i h_1(X_i)h_1(X_i)^\top\right\|_\infty +$$
$$\left\|\frac{1}{n}\sum_i h_1(X_i)h_1(X_i)^\top - \mathbb{E}\left[h_1(X)h_1(X)^\top\right]\right\|_\infty \ . \quad \text{(F.20)}$$

We handle these two terms separately. First, we bound $||\frac{1}{n}\sum_i \widehat{h}_1(X_i)\widehat{h}_1(X_i)^\top - \frac{1}{n}\sum_i h_1(X_i)h_1(X_i)^\top||_\infty$. For convenience we define $A := \frac{1}{n}\sum_i \widehat{h}_1(X_i)\widehat{h}_1(X_i)^\top$ and $B := \frac{1}{n}\sum_i h_1(X_i)h_1(X_i)^\top$. Since $A$ and $B$ are both positive semidefinite matrices with ones on the diagonal, we have

$$||A - B||_\infty = \frac{1}{2} \max_{j,k \in [p_n^2]} \left|f_{jk}^\top (A - B) f_{jk}\right| \ , \quad \text{(F.21)}$$

where $f_{jk} \in \mathbb{R}^{p_n^2}$ is the vector with $(f_{jk})_j = 1$, $(f_{jk})_k = -1$, and zeros elsewhere. Next we have

$$\left|f_{jk}^\top (A - B) f_{jk}\right| = \left|\sqrt{f_{jk}^\top A f_{jk}} - \sqrt{f_{jk}^\top B f_{jk}}\right| \cdot \left(\sqrt{f_{jk}^\top A f_{jk}} + \sqrt{f_{jk}^\top B f_{jk}}\right)$$
$$\leq 4\left|\sqrt{f_{jk}^\top A f_{jk}} - \sqrt{f_{jk}^\top B f_{jk}}\right| = \frac{4}{\sqrt{n}}\left|\sqrt{\sum_i (\widehat{h}_1(X_i)^\top f_{jk})^2} - \sqrt{\sum_i (h_1(X_i)^\top f_{jk})^2}\right|$$
$$\leq \frac{4}{\sqrt{n}} \sqrt{\sum_i \left((\widehat{h}_1(X_i) - h_1(X_i))^\top f_{jk}\right)^2} \ ,$$

where the first inequality follows from the fact that $||f_{jk}||_1 \leq 2$ while $||A||_\infty, ||B||_\infty \leq 1$, and the second inequality follows from the triangle inequality. Next, for each $i$ and each $j, k$, observe that

$$\widehat{h}_1(X_i)^\top f_{jk} = \frac{1}{n-1} \sum_{i' \neq i} \left(\text{sign}(X_i - X_{i'}) \otimes \text{sign}(X_i - X_{i'})\right)^\top f_{jk},$$

which after conditioning on $X_i$, is a mean of $(n-1)$ i.i.d. variables, each taking values in $[-2, 2]$ since $||f_{jk}||_1 \leq 2$. Furthermore, conditioning on $X_i$, we have $\mathbb{E}[\widehat{h}_1(X_i)] = h_1(X_i)$. Therefore, applying Hoeffding's lemma [see, for example, Lemma 2.6 in Massart, 2007], for each $i, j, k$, for any $t \in \mathbb{R}$,

$$\mathbb{E}\left[\exp\left\{t \cdot (\widehat{h}_1(X_i) - h_1(X_i))^\top f_{jk}\right\}\right] \leq \exp\left\{\frac{2t^2}{n-1}\right\}. \tag{F.22}$$

Applying Lemma F.10 (stated below), then,

$$\mathbb{P}\left\{\frac{1}{n}\sum_i \left((\widehat{h}_1(X_i) - h_1(X_i))^\top f_{jk}\right)^2 > \frac{80}{n-1} \cdot (1 + \log(27 p_n^5))\right\} \leq \frac{1}{27 p_n^5}.$$

Taking a union bound over all $j, k \in [p_n^2]$, and returning to (F.21), we then have

$$\mathbb{P}\left\{\left\|\frac{1}{n}\sum_i \widehat{h}_1(X_i)\widehat{h}_1(X_i)^\top - \frac{1}{n}\sum_i h_1(X_i)h_1(X_i)^\top\right\|_\infty > 2\sqrt{\frac{80}{n-1} \cdot (1 + \log(27 p_n^5))}\right\} \leq \frac{1}{27 p_n}.$$

Next we turn to the second term in (F.20). Since $||h_1(X)||_\infty \leq 1$ always, we see that for each $j, k \in [p_n]$,

$$\left(\frac{1}{n}\sum_i h_1(X_i)h_1(X_i)^\top\right)_{jk}$$

is a mean of $n$ i.i.d. terms, each taking values in $[-1, 1]$. Applying Hoeffding's inequality, for each $j, k$,

$$\mathbb{P}\left\{\left|\left(\frac{1}{n}\sum_i h_1(X_i)h_1(X_i)^\top - \mathbb{E}[h(X)h(X)^\top]\right)_{jk}\right| \geq t\right\} \leq 2e^{-nt^2/2}$$

for any $t \geq 0$. Setting $t = \sqrt{\frac{2\log(54 p_n^3)}{n}}$, and taking a union bound, we see that

$$\mathbb{P}\left\{\left\|\frac{1}{n}\sum_i h_1(X_i)h_1(X_i)^\top - \mathbb{E}[h(X)h(X)^\top]\right\|_\infty \geq \sqrt{\frac{2\log(54 p_n^3)}{n}}\right\}$$

$$\leq 2p_n^2 \cdot e^{-n\left(\sqrt{\frac{54\log(p_n^3)}{n}}\right)^2/2} = \frac{1}{27 p_n}.$$

Returning to (F.20), then, with probability at least $1 - \frac{2}{27 p_n}$,

$$\left\|\frac{1}{n}\sum_i \widehat{h}_1(X_i)\widehat{h}_1(X_i)^\top - \mathbb{E}\left[h_1(X)h_1(X)^\top\right]\right\|_\infty \leq 2\sqrt{\frac{80}{n-1} \cdot (1 + \log(27 p_n^3))} + \sqrt{\frac{2\log(54 p_n^3)}{n}}. \tag{F.23}$$

This proves a bound on (F.20).

Next, to complete the proof, we bound

$$\left\|\left(\frac{1}{n}\sum_i \widehat{h}_1(X_i)\right)\left(\frac{1}{n}\sum_i \widehat{h}_1(X_i)\right)^\top - \mathbb{E}\left[h_1(X)\right]\mathbb{E}\left[h_1(X)\right]^\top\right\|_\infty.$$

We have

$$\left(\frac{1}{n}\sum_i \widehat{h}_1(X_i)\right)\left(\frac{1}{n}\sum_i \widehat{h}_1(X_i)\right)^\top - \mathbb{E}\left[h_1(X)\right]\mathbb{E}\left[h_1(X)\right]^\top$$

$$= \left(\frac{1}{n}\sum_i \widehat{h}_1(X_i)\right)\left(\frac{1}{n}\sum_i \widehat{h}_1(X_i) - \mathbb{E}\left[h_1(X)\right]\right)^\top - \left(\frac{1}{n}\sum_i \widehat{h}_1(X_i) - \mathbb{E}\left[h_1(X)\right]\right)\mathbb{E}\left[h_1(X)\right]^\top$$

and, since $||\frac{1}{n}\sum_i \widehat{h}_1(X_i)||_\infty, ||\mathbb{E}\left[h_1(X)\right]||_\infty \leq 1$, we therefore have

$$\left\|\left(\frac{1}{n}\sum_i \widehat{h}_1(X_i)\right)\left(\frac{1}{n}\sum_i \widehat{h}_1(X_i)\right)^\top - \mathbb{E}\left[h_1(X)\right]\mathbb{E}\left[h_1(X)\right]^\top\right\|_\infty \leq 2||\left\|\frac{1}{n}\sum_i \widehat{h}_1(X_i) - \mathbb{E}\left[h_1(X)\right]\right\|_\infty.$$

For each sign $s \in \{\pm 1\}$, for each $j \in [p_n^2]$, writing $\mathbf{e}_j$ to denote the $j$th basis vector in $\mathbb{R}^{p_n^2}$, we have

$$\mathbb{E}\left[\exp\left\{t \cdot s \cdot \mathbf{e}_j^\top \left(\frac{1}{n}\sum_i \widehat{h}_1(X_i) - \mathbb{E}\left[h_1(X)\right]\right)\right\}\right]$$
$$\leq \frac{1}{n}\sum_i \mathbb{E}\left[\exp\left\{t \cdot s \cdot \mathbf{e}_j^\top \left(\widehat{h}_1(X_i) - \mathbb{E}\left[h_1(X)\right]\right)\right\}\right] \leq \exp\left\{\frac{t^2}{2(n-1)}\right\},$$

where the first inequality follows from the convexity of $x \mapsto e^x$, while the second applies Hoeffding's lemma, as in (F.22) above. Then,

$$\mathbb{P}\left\{s \cdot \mathbf{e}_j^\top \left(\widehat{h}_1(X_i) - \mathbb{E}\left[h_1(X)\right]\right) > \sqrt{\frac{2\log(54p_n^3)}{n-1}}\right\} \leq \frac{1}{54p^3},$$

and therefore taking a union bound over each $s \in \{\pm 1\}$ and each $j \in [p_n^2]$,

$$\mathbb{P}\left\{\left\|\frac{1}{n}\sum_i \widehat{h}_1(X_i) - \mathbb{E}\left[h_1(X)\right]\right\|_\infty > \sqrt{\frac{2\log(54p_n^3)}{n-1}}\right\} \leq \frac{1}{27p_n}.$$

Therefore, combining this with (F.23), with probability at least $1 - \frac{1}{9p_n}$,

$$||\widehat{\Sigma}_{h_1} - \Sigma_{h_1}||_\infty \leq 2\sqrt{\frac{80}{n-1} \cdot (1+\log(27p_n^3))} + \sqrt{\frac{2\log(54p_n^3)}{n}} + 2\sqrt{\frac{2\log(54p_n^3)}{n-1}} \leq 100\sqrt{\frac{\log(p_n)}{n}},$$

where the last step uses the fact that $n, p_n \geq 2$. $\square$

**Lemma F.9.** *Let $\Sigma_{h_1}$ be defined as in Assumption 3.4. For every $z \in \mathbb{R}^{p_n^2}$,*

$$z^\top \Sigma_{h_1} z \leq \lambda_{\max}(\Sigma) \cdot ||\mathsf{mat}(z)||_{\ell_1/\ell_2}^2,$$

*where $||\mathsf{mat}(z)||_{\ell_1/\ell_2}$ is defined as in the statement of Lemma F.6.*

*Proof of Lemma F.9.* Since the statement is deterministic, we can treat $M = \mathsf{mat}(z) \in \mathbb{R}^{p_n \times p_n}$ as fixed.

Then $z = \mathsf{vec}(M)$ and

$$\begin{aligned}
\mathsf{vec}(M)^\top \Sigma_{h_1} \mathsf{vec}(M) &= \mathsf{Var}\left(\mathsf{vec}(M)^\top h_1(X)\right) \\
&= \mathsf{Var}\left(\mathsf{vec}(M)^\top \mathbb{E}[h(X, X') \mid X]\right) \\
&= \mathsf{Var}\left(\mathbb{E}[\mathsf{vec}(M)^\top h(X, X') \mid X]\right) \\
&\leq \mathsf{Var}\left(\mathsf{vec}(M)^\top h(X, X')\right) \quad \text{(by the law of total variance)} \\
&\leq \mathbb{E}\left[(\mathsf{vec}(M)^\top h(X, X'))^2\right] \\
&= \mathbb{E}\left[(\mathsf{vec}(M)^\top \left(\mathsf{sign}(X - X') \otimes \mathsf{sign}(X - X'))\right)^2\right] \\
&= \mathbb{E}\left[\left(\mathsf{sign}(X - X')^\top M \, \mathsf{sign}(X - X')\right)^2\right] \\
&= \mathbb{E}\left[\left(\sum_j \mathsf{sign}(X - X')^\top M_j \cdot \mathsf{sign}(X_j - X'_j)\right)^2\right] \quad \text{(where } M_j \text{ is the } j\text{th column of } M) \\
&\leq \mathbb{E}\left[\left(\sum_j |\mathsf{sign}(X - X')^\top M_j|\right)^2\right] \\
&= \sum_{jk} \mathbb{E}\left[|\mathsf{sign}(X - X')^\top M_j| \cdot |\mathsf{sign}(X - X')^\top M_k|\right] \\
&\leq \sum_{jk} \sqrt{\mathbb{E}\left[|\mathsf{sign}(X - X')^\top M_j|^2\right]} \cdot \sqrt{\mathbb{E}\left[|\mathsf{sign}(X - X')^\top M_k|^2\right]} \\
&= \sum_{jk} \sqrt{M_j^\top \mathbb{E}\left[\mathsf{sign}(X - X') \mathsf{sign}(X - X')^\top\right] M_j} \cdot \sqrt{M_k^\top \mathbb{E}\left[\mathsf{sign}(X - X') \mathsf{sign}(X - X')^\top\right] M_k} \\
&= \sum_{jk} \sqrt{M_j^\top T M_j} \cdot \sqrt{M_k^\top T M_k} \\
&\leq \sum_{jk} \sqrt{\|M_j\|_2^2 \cdot \lambda_{\max}(T)} \cdot \sqrt{\|M_k\|_2^2 \cdot \lambda_{\max}(T)} \\
&= \lambda_{\max}(T) \cdot \left(\sum_j \|M_j\|_2\right)^2.
\end{aligned}$$

Finally, by Wegkamp and Zhao [2013, Theorem 2.3], $\lambda_{\max}(T) \leq \lambda_{\max}(\Sigma)$. □

**Lemma F.10.** *Let $v \in \mathbb{R}^p$ be a fixed vector and let $Z_1, \ldots, Z_n \in [-1, 1]^p$ be random vectors, not necessarily independent, such that $v^\top(Z_i - \mathbb{E}[Z_i])$ is $C$-subgaussian for each $i$, that is,*

$$\mathbb{E}[\exp\{tv^\top(Z_i - \mathbb{E}[Z_i])\}] \leq \exp(Ct^2/2).$$

*Then for any $\delta \in (0, 1)$, with probability at least $1 - \delta$,*

$$\frac{1}{n} \sum_i \left(v^\top(Z_i - \mathbb{E}[Z_i])\right)^2 \leq 20C(1 + \log(1/\delta)).$$

*Proof of Lemma F.10.* For each $i$, by assumption,

$$\mathbb{E}\left[\exp\left\{t \cdot \frac{1}{\sqrt{C}} v^\top(Z_i - \mathbb{E}[Z_i])\right\}\right] \leq \exp\left\{\frac{t^2}{2}\right\}.$$

By Vershynin [2012, Lemma 5.5] (and tracking constants carefully in this Lemma), for each $i$,

$$\mathbb{E}\left[\exp\left\{\frac{1}{20C} \cdot \left(v^\top(Z_i - \mathbb{E}[Z_i])\right)^2\right\}\right] \leq e.$$

By the convexity of $x \mapsto e^x$, then,

$$\mathbb{E}\left[\exp\left\{\frac{1}{20C} \cdot \frac{1}{n}\sum_i \left(v^\top(Z_i - \mathbb{E}[Z_i])\right)^2\right\}\right] \leq \frac{1}{n}\sum_i \mathbb{E}\left[\exp\left\{\frac{1}{20C} \cdot \left(v^\top(Z_i - \mathbb{E}[Z_i])\right)^2\right\}\right] \leq e \ .$$

Therefore, we have

$$\mathbb{P}\left\{\frac{1}{n}\sum_i \left(v^\top(Z_i - \mathbb{E}[Z_i])\right)^2 > t\right\} \leq \mathbb{E}\left[\exp\left\{\frac{1}{20C}\frac{1}{n}\sum_i \left(v^\top(Z_i - \mathbb{E}[Z_i])\right)^2 - \frac{1}{20C}t\right\}\right]$$
$$\leq \exp\left\{1 - \frac{1}{20C}t\right\} \ .$$

Setting $t = 20C(1 + \log(1/\delta))$, then, we have proved the desired result.

$\square$